\documentclass[12pt]{article}
\usepackage{amsmath}
\usepackage{amsfonts}
\usepackage{latexsym}
\usepackage{graphicx,psfrag,epsf}
\usepackage{enumerate}
\usepackage{color}
\usepackage{colortbl}
\usepackage{natbib}
\usepackage{multirow}
\usepackage{url} % not crucial - just used below for the URL 

\newcommand{\blind}{0}

    \numberwithin{equation}{section}
    \newtheorem{thrm}{Theorem}[section]             
    \newtheorem{prop}[thrm]{Proposition}

    \newtheorem{lmma}[thrm]{Lemma}
    \newtheorem{remk}[thrm]{Remark}

    \def\calG{{\mathcal G}}

    \def\bbr{{\mathbb R}}
    \def\bbe{\mathbb E}
    \def\bbp{{\mathbb P}}
    
    \def\bbn{{\mathbb N}}

    \def\m{{m}}

    \def\cp{\;{\stackrel{\bbp}{\longrightarrow}}\;}

    \definecolor{Red}{rgb}{0,0,0}%{1,0,0}%{1.00, 0.00, 0.00}
    
    \definecolor{DRed}{rgb}{0,0,0}%{0.0, 0.00, 0.00}
    
    \definecolor{Blue}{rgb}{0,0,1}%{0,0,1}%{0.00, 0.00, 0.00}
    
    \definecolor{PaleGrey}{rgb}{0,0,0}%{.6, .6, .6}

\begin{document}

\def\spacingset#1{\renewcommand{\baselinestretch}%
{#1}\small\normalsize} \spacingset{1}

%%%%%%%%%%%%%%%%%%%%%%%%%%%%%%%%%%%%%%%%%%%%%%%%%%%%%%%%%%%%%%%%%%%%%%%%%%%%%%

\if0\blind
{
  \title{\bf Estimation of a noisy subordinated Brownian Motion via {two-scales} power variations}
  \author{Jos\'e E. Figueroa-L\'opez\thanks{
  The first author's research is partially supported by the NSF {grants DMS-1149692 and DMS-1613016}.}\hspace{.2cm}\\
    Department of Mathematics, Washington University in St. Louis\\
    and \\
    Kiseop Lee \\
    Department of Statistics, Purdue University}
  \maketitle
} \fi

\if1\blind
{
  \bigskip
  \bigskip
  \bigskip
  \begin{center}
    {\LARGE\bf Estimation of a noisy subordinated Brownian Motion via {two-scales} power variations}
\end{center}
  \medskip
} \fi

\bigskip
\begin{abstract}
High frequency based estimation methods for a semiparametric pure-jump subordinated Brownian motion exposed to a small additive microstructure noise are developed {building on} the {two-scales} realized variations approach {originally} developed by \cite{ZhaMykAit:2005} for the estimation of {the integrated variance of a} continuous It\^o process. The proposed estimators are shown to be robust against the noise and, {surprisingly,} to attain better rates of convergence than {their precursors, method of moment estimators}, \emph{even in the absence of microstructure noise}. Our main results give approximate optimal values for the number $K$ of regular sparse  subsamples to be used, which is {an important} tune-up parameter of the method. Finally, a data-driven {plug-in} procedure is devised to implement the proposed estimators with the optimal $K$-value. The {developed estimators exhibit superior performance as illustrated by Monte Carlo simulations and a real high-frequency data application}.
\end{abstract}

\noindent%
\vfill

\newpage
\spacingset{1.45} % DON'T change the spacing!
\section{Introduction}
\label{sec:intro}

In this paper, we develop estimation methods for a {semiparametric} subordinated Brownian motion (SBM),  whose sampling observations have been contaminated by a small additive noise {along the lines of the framework of \cite{ZhaMykAit:2005}}. 
In addition to a ``volatility" parameter $\sigma$, which {controls} the variance of {the increments of the process at regular time intervals}, a SBM is endowed with an additional parameter, {hereafter} denoted by $\kappa$, which accounts for the tail heaviness of the increments' distribution. {Therefore,} $\kappa$ determines the proneness of the process to produce extreme increment observations. Such a measure is clearly of critical relevance in many applications such as {to model extreme events in insurance and risk management and optimal asset allocation in finance}.  
The models considered here are pure-jump L\'evy models and $\sigma$ is not the volatility of a continuous {It\^o process}. {Nevertheless, given that $\sigma^{2}$ is proportional to the variance of {the increments of the process}, it is natural to refer to $\sigma$ as the volatility parameter of the model.}

As in the context of a regression model, the additive noise, typically called microstructure noise, can be seen as a modeling artifact to {account for} any deviations between the {observed process} and the SBM {model}.  {However,} {in some circumstances, the noise can be {link to some specific physical mechanism} such as in the case of bid/aks bounce effects {in tick by tick trading} (cf. \cite{Roll}).} At low frequencies the microstructure noise is {typically} negligible {(compared to the SBM's increments),} but at {high-frequencies} the noise is significant and heavily tilts any estimates that {do} not account for it.  The aim is then to develop inference methods that are robust against potential microstructure noises.

The literature of statistical estimation methods under microstructure noise has grown extensively during the last decade. {See \cite{AitJacodBook} for a recent in depth survey on the topic and, also, \cite{AitMykZh05}, \cite{ZhaMykAit:2005}, \cite{HanLun}, \cite{Bandi}, \cite{MyklandZhang2012} {for a few seminal works in the area}.  Most of these works have focused on the estimation of the integrated variance of a semimartingale model. However,  the problem of translating 
some of the proposed methods into estimation methods for semiparametric models contaminated by additive noise, as it is the case in the present work, has received much less attention in the literature, in particular, when it comes to the estimation of a kurtosis type parameter.} The performances of {some classical parametric methods in the estimation} {of} some popular parametric L\'evy models have been analyzed in a few works such as \cite{Seneta}, \cite{Zeng:2007}, \cite{Behr}, and \cite{FLLLM}, {but none of them have incorporated microstructure noise.}

{To motivate our estimation procedure, we start by considering Method of Moment Estimators (MME) for $\sigma^{2}$ and $\kappa$, in the absence of microstructure noise. Throughout the remainder of the introduction, these estimators are respectively denoted by $\hat\sigma_{n,T}^{2}$ and $\hat\kappa_{n,T}$, where $n$ and $T$ denote the number of observations and the sampling horizon, respectively}.  {MMEs and related estimators} are widely used in high-frequency data analysis due to their simplicity, computational efficiency, and known robustness against potential correlation between observations.  
{In order to establish asymptotic benchmarks for the convergence rates of our proposed estimators, we characterize the asymptotic behavior of the MME estimators, both in the absence and presence of microstructure noise,} when {$\delta_n=T/n$, the time span between observations, shrink to $0$} (infill asymptotics) and $T\to\infty$ (long-run asymptotics). We identify the order $O(T^{-1})$, as the {rate} of convergence of the estimators {under the} absence of noise. Hence, {a desirable objective} is to develop estimators that are able to achieve at least this rate of convergence in the presence of microstructure noise. An asymptotic analysis of the estimators in the presence of noise allows to {show that} $\hat\sigma^{2}_{n,T}\to\infty$ and $\hat\kappa_{n,T}\to{}0$, as {$n\to{}\infty$, both of which are stylized {empirical properties of} high-frequency {financial} observations (see Section \ref{Sec:EmStudy} below). {Furthermore, it is shown} that $\delta_{n}\hat\sigma^{2}_{n,T}$ and $\delta_{n}^{-1}\hat\kappa_{n,T}$ converge to the second moment and the excess kurtosis of the microstructure noise, respectively. 

In order to develop estimators that are robust {against} a microstructure noise component, we borrow ideas from \cite{ZhaMykAit:2005}'s seminal approach based on combining the realized quadratic variations at two-scales or frequencies. More concretely, {there are three main steps in this approach}. First, {the high-frequency sampling observations are divided} in $K$ groups of observations taken at a lower frequency (sparse subsampling). Second, the relevant estimators {(say, realized quadratic variations)} are applied to each group and the resulting $K$ point estimates are averaged. Finally, a bias correction step is necessary for which one typically uses the estimators at the highest possible frequency. 

A fundamental problem in the approach described in the previous paragraph is {how to tune up} the number of subgroups, $K$, {which strongly affects the performance of the estimators}. We propose a method to find approximate optimal values for $K$ {under a white microstructure noise setting}. For the estimator of $\sigma^{2}$, it is found that the optimal $K$ takes the form
\begin{equation}\label{Exp2OK0}
	K^{*}_{\sigma}:=n^{\frac{2}{3}}\left(\frac{6\left(\bbe\varepsilon^{4}+(\bbe\varepsilon^{2})^{2}\right)}{T^{2}\sigma^{4}}\right)^{\frac{1}{3}},
\end{equation}
where $\varepsilon$ represents the additive microstructure noise associated to one observation of the SBM. {Interestingly, the optimal value (\ref{Exp2OK0}) is consistent, but different from that proposed by \cite{ZhaMykAit:2005}  in the context of a continuous It\^o semimartingale\footnote{{The optimal value of $K$ proposed in  \cite{ZhaMykAit:2005} (see Eq.~(58) and (63) therein) lacks the term $\bbe\varepsilon^{4}$ in the numerator.}}}. It is also found that the mean-squared error (MSE) of the resulting estimator (using $K$ as above) attains a rate of convergence $C_{\sigma}\left(\bbe\varepsilon^{4}+(\bbe \varepsilon^{2})^{2}\right)^{\frac{1}{3}}n^{-\frac{1}{3}}T^{-\frac{2}{3}}$ (up to a constant $C_{\sigma}$), which, since $T/n\to{}0$, shows the surprising fact that the estimator converges at a rate of $o(T^{-1})$, which is faster than the rate attained by the MMEs in the absence of noise. For the estimation of $\kappa$, it is found that the optimal $K$ takes the form 
\begin{equation}\label{Exp1OK3}
	K^{*}_{\kappa}=n^{\frac{4}{5}}\left(\frac{5\,{\rm Var}\left((\varepsilon_{2}-\varepsilon_{1})^{4}\right)}{3^{3}2^{4}T^{4}\sigma^{8}}\right)^{\frac{1}{5}},
\end{equation}
while the mean-squared error of the resulting estimator converges at the rate of 
\[
C_{\kappa} {\rm Var}\left((\varepsilon_{2}-\varepsilon_{1})^{4}\right)^{\frac{3}{5}}n^{-\frac{3}{5}}T^{-\frac{2}{5}},
\]
 up to constant $C_{\kappa}$. Here, $\varepsilon_{1}$ and $\varepsilon_{2}$ represent the microstructure noise {corresponding to} two different observations of the SBM. In particular, we again infer that the resulting estimator attains a better MSE performance than the plain MME in the absence of noise. 

In order to implement the estimators with the corresponding optimal choices of $K^{*}$, we propose an iterative procedure in which an initial reasonable guess for $\sigma^{2}$ is used to find $K^{*}$, which in turn is used to improve the initial guess of $\sigma$, and so forth. The  resulting estimators {exhibit superior finite-sample performance both on simulated and real high-frequency stock data. In particular, we found that} the estimators {are quite stable as} the sampling frequency increases, {when compared to} their MME counterparts, {which, as mentioned above, converge to either $0$ or $\infty$ for $\sigma$ or $\kappa$, respectively.}

The rest of the paper is organized as follows. In Section \ref{SctModels}, we give the model and the estimation framework. Section \ref{MMESec} introduces the method of moment estimators. Their in-fill and long-run asymptotic behavior are analyzed in Section \ref{MNSct}. Section \ref{MMECorrected} introduces {the} estimators for $\sigma$ and $\kappa$ that are robust to a microstructure noise component together with bias corrected versions of these with optimal selection of $K$. Section \ref{Numerics} shows the finite-sample performance of the proposed estimators via simulations as well as their empirical robustness using real high-frequency transaction data. Finally, the proofs of the paper are deferred to the Appendix.

\section{The model and the sampling scheme}\label{SctModels}

In this section, we introduce the model used throughout the paper. We {consider a subordinated Brownian motion of the form  
\begin{equation}\label{VGM}
	X_{t}=\sigma W(\tau_{t}) + \theta \tau_{t} + bt,
\end{equation}
where $\sigma,\kappa>0,\theta,b\in\bbr$,} $W:=\{W(t)\}_{t\geq{}0}$ is a {standard Brownian motion}, and $\{\tau_{t}\}_{t\geq{}0}:=\{\tau(t;\kappa)\}_{t\geq{}0}$ is an independent subordinator (i.e., a non-decreasing L\'evy process) satisfying the following conditions:
\begin{equation}\label{MnVarRndClock}
	{\rm (i)}\;\bbe\tau_{t}=t,\qquad {\rm (ii)}\; {\rm Var}(\tau_{t})=\kappa t,\qquad%\text{ and }\quad 
	{\rm (iii)}\;\bbe\tau_{1}^{j}<\infty, \quad{ j=1,\dots,{}8.}
\end{equation} 
{The first condition is needed for identifiability purposes, while the second one allows to interpret $\kappa$ {as a measure} of the excess kurtosis. The} condition (\ref{MnVarRndClock}-iii) is imposed so that $X_{t}$ admits finite moments of {sufficiently large order}. {In financial applications, $X$ is often interpreted as the log-return process $X_{t}=\log(S_{t}/S_{0})$ of a risky asset with price process  $\{S_{t}\}_{t\geq{}0}$. In that case, $\tau$ plays the role of} a random clock aimed at incorporating variations in business {``activity"} through time. 
It is well known that the process $X$ is a L\'evy process (see, e.g,  \cite{Sato}).  Hereafter, $\nu$ will denote the L\'evy measure of $X$, {which controls the jump behavior of the process in that $\nu((x,x+dx))$ measures the expected number of jumps with size near $x$ per unit time.}

{Two prototypical} examples of (\ref{VGM}) are the Variance Gamma (VG) and the Normal Inverse Gaussian (NIG) L\'evy processes, which were proposed by \cite{Madan1} and \cite{Barndorff:1998}, respectively. In the VG model, {$\tau(t;\kappa)$} is Gamma distributed with scale parameter $\beta:=\kappa$ and shape parameter $\alpha:=t/\kappa$, while in the NIG model {$\tau(t;\kappa)$} follows an Inverse Gaussian distribution with mean $\mu=1$ and shape parameter $\lambda=1/(t\kappa)$. 

{As seen from the formulas for their moments (see (\ref{MEVGM}) below), the} model's parameters have the following interpretation:
\begin{enumerate}
	\item $\sigma$ dictates the overall variability of {the process' increments or, in financial terms, the} log returns of the asset; 
	in the ``symmetric" case ($\theta=0$), $\sigma^{2}$ is the variance of log returns {divided by the time span of the returns}; 
	\item $\kappa$ controls the kurtosis or {the tail's} heaviness of the log return distribution; in the symmetric case ($\theta=0$), $\kappa$ is the excess kurtosis of log returns multiplied by the time span {of the returns}; 
	\item $b$ is a drift component in the calendar time; 
	\item $\theta$ is a drift component in the business time and controls the skewness of log returns; 
\end{enumerate}

Throughout the paper, we {also} assume that the log return process $\{X_{t}\}_{t\geq{}0}$ is sampled during a {time} interval $[0,T]$ at evenly spaced times:
\begin{equation}\label{SmplSchm}
	{t_{i,n}}=t_{i}:=i\delta_{n}, \quad i=1,\dots,n,\quad \text { where }\quad\delta_{n}:=\frac{T}{n}.
\end{equation}
This sampling scheme is sometimes called \emph{calendar time sampling} (c.f. \cite{Oomen:2006}). Under the assumption of independence and stationarity of increments, we have at our disposal a random sample 
\begin{equation}\label{SD}
	\Delta_{i}^{n}X:=X_{i\delta_{n}}-X_{(i-1)\delta_{n}},\quad i=1,\dots, n,
\end{equation}
of size $n$ of the {distribution} of $X_{\delta_{n}}$.

In real markets, {high-frequency} log returns exhibit {certain} stylized features{,} which cannot {be accurately explained} by efficient models such as (\ref{SD}). There are different approaches to model these features, widely termed as \emph{microstructure noise}. Microstructure noises {may} come from different sources, such as clustering noises, non-clustering noises {such as bid/ask bounce effects}, and roundoff errors  {(cf. \cite{Campbell}, \cite{Zeng:2003})}.  
In what follows, we adopt a popular approach due to \cite{ZhaMykAit:2005}, where the {net} effect of the market microstructure is incorporated as an additive noise to the observed log-return process:
\begin{equation}\label{OPDfn}
	{\widetilde{X}_t:=\widetilde{X}(t):=X_t+\varepsilon_t},
\end{equation}
where $\{\varepsilon_{t}\}_{t\geq{}0}$ is assumed to be a centered process, independent of $X$. In particular, under this setup, the log return observations at a frequency $\delta_{n}$ are given by 
\begin{equation}\label{MILMN20}
	\Delta_{i}^{n}\widetilde{X}:=\widetilde{X}_{i \delta_n} - \widetilde{X}_{(i-1) \delta_n} 
=
	 \Delta_{i}^{n} X+\tilde{\varepsilon}_{i,{\delta_{n}}},%\qquad i=1,\dots,n.
\end{equation}
{where $\tilde\varepsilon_{i,\delta}:=\varepsilon_{i\delta}-\varepsilon_{(i-1)\delta}$ can be interpreted as the contribution of the microstructure noise to the observed increment $\Delta_{i}^{n}\widetilde{X}$.} In the simplest case, the noise $\{\varepsilon_{t}\}_{t\geq{}0}$ is a white noise; i.e., the variables $\{\varepsilon_{t}\}_{t\geq{}0}$ are independent identically distributed with mean $0$. 

It is well known {(and not surprising)} that standard statistical methods do not perform well when applied to high-frequency observations if {the microstructure noise is} {not taken into account}. A standing problem is then to derive {inference} methods that are robust against a wide range of microstructure noises. In Section \ref{MMECorrected}, we proposed an approach to address the latter problem, borrowing ideas from the {seminal} {two-scales} correction technique of \cite{ZhaMykAit:2005} applied to Method of Moment Estimators {(MME)}. Before that, we first introduce {the considered MMEs} and carry on a simple infill asymptotic analysis of the estimators both in the absence and presence of the microstructure noise.

\section{Method of Moment Estimators}\label{MMESec}

The Method of Moment Estimators (MME) are widely used to deal with high-frequency data due to their simplicity, computational efficiency, and known robustness against potential correlation between observations. For the general subordinated Brownian model (\ref{MnVarRndClock})-(\ref{VGM}), the central moments {can easily be computed} in closed forms as
\begin{align}
	\mu_{1}(X_{\delta})&:=\bbe(X_{\delta})=(\theta+b)\delta,\quad 
	\mu_{2}(X_{\delta}):={\rm Var}(X_{\delta})=(\sigma^{2} +\theta^{2}\kappa)\delta,\nonumber\\
	\mu_{3}(X_{\delta})&:=\bbe(X_{\delta}-\bbe X_{\delta})^{3}=\left(3\sigma^{2} \theta \kappa + \theta^{3}c_{3}(\tau_{1})\right)\delta,\label{MEVGM} \\
	\mu_{4}(X_{\delta})&:=\bbe(X_{\delta}-\bbe X_{\delta})^{4}=\left(3\sigma^{4}\kappa +6\sigma^{2}\theta^{2}c_{3}(\tau_{1})+\theta^{4}c_{4}(\tau_{1})\right)\delta+3\mu_{2}(X_{\delta})^{2},\nonumber
\end{align}
where, hereafter, 
\[
	c_{k}(Y):=\frac{1}{i^{{k}}}\left.\frac{d^{k}}{du^{k}}\ln\bbe\left(e^{iuY}\right)\right|_{u=0},
\] 
represents the $k$-th cumulant of a r.v. $Y$. For the VG model, $\left(c_{3}(\tau_{1}),c_{4}(\tau_{1})\right)=(2\kappa^{2},6\kappa^{3})$, while for the NIG model, $\left(c_{3}(\tau_{1}),c_{4}(\tau_{1})\right)=(3\kappa^{2},15\kappa^{3})$. 

Throughout, we assume that $\theta=0$ {or, more generally, that $\theta$ is negligible compare to $\sigma$ (see Remark \ref{skewness0} below for further discussion about this assumption).
The assumption that $\theta=0$ allows us to propose} tractable expressions for the MME of the parameters $\sigma^{2}$ and $\kappa$ {as follows:}
\begin{align}\label{SME1}	
	{\tilde\sigma_{n}}^{2}(X):=\frac{1}{\delta_{n}}\hat{\mu}_{2,n}(X),\qquad
	{\tilde\kappa_{n}(X)}:={\frac{\delta_{n}}{3}\frac{\hat\mu_{4,n}(X)}{\hat\mu_{2,n}^{2}(X)}-\delta_{n}},
\end{align}
%where 
%$\widehat{\rm Var}_{n}(X)$ and $\widehat{\rm Krt}_{n}(X)$ denote the sample variance and kurtosis of the increments of $X$. Concretely, 
%\begin{equation}\label{SVVSK}
%	\widehat{\rm Var}_{n}(X):=\hat{\mu}_{2,n}(X),\qquad
%	\widehat{\rm Krt}_{n}(X):=\frac{\hat\mu_{4,n}(X)}{\hat\mu_{2,n}^{2}(X)}-3.
%\end{equation}
where {hereafter} $\hat\mu_{k,n}(X)$ represents the sample central moment of $k^{th}$ order as defined by
\begin{equation}\label{DSM}
	\hat\mu_{k,n}(X):=\frac{1}{n}\sum_{i=1}^{n}\left(\Delta_{i}^{n}X-\overline{\Delta^{n}X}\right)^{k},\quad k\geq 2,
	\quad \overline{\Delta^{n}X}:=\frac{1}{n}\sum_{i=1}^{n}\Delta_{i}^{n}X=\frac{1}{n}\log\frac{S_{T}}{S_{0}}.
\end{equation}
We can further simplify the above statistics by omitting the terms of order $O(\delta_{n})=O(1/n)$ {(in particular, we leave out the term $\delta_{n}$ in (\ref{SME1}) and $\overline{\Delta^{n}X}$ in sample moments of (\ref{DSM})):}
\begin{align}\label{SME1b}
	{\hat\sigma_{n}^{2}(X)}
%	:=\frac{1}{\delta_{n}}\hat{m}_{2,n}(X)=\frac{1}{T}\sum_{i=1}^{n}\left(\Delta_{i}^{n}X\right)^{2}
	:=\frac{1}{T}\left[X,X\right]_{2},\qquad 
	{\hat{\kappa}_{n}(X)}
%	:=\frac{\delta_{n}}{3}\frac{\hat{m}_{4,n}(X)}{\hat{m}_{2,n}^{2}(X)}
	:=
	\frac{\delta_{n}}{3}\frac{\frac{1}{n}\sum_{i=1}^{n}\left(\Delta_{i}^{n}X\right)^{4}}{\left(\frac{1}{n}\sum_{i=1}^{n}\left(\Delta_{i}^{n}X\right)^{2}\right)^{2}}=\frac{1}{3}\frac{T^{-1}\widehat{\left[X,X\right]}_{4}}{\left(T^{-1}\widehat{\left[X,X\right]}_{2}\right)^{2}},
\end{align}
where above we have expressed the estimators in terms of the realized variations of order $2$ and $4$, {which hereafter are defined by}
\[
	{\widehat{\left[X,X\right]}_{2}:=\sum_{i=1}^{n}\left(\Delta_{i}^{n}X\right)^{2},
	\qquad 
	\widehat{\left[X,X\right]}_{4}= \sum_{i=1}^{n}\left(\Delta_{i}^{n}X\right)^{4}}.
\]
\begin{remk}\label{skewness0}
{In the case that $|\theta|<<\sigma$ (i.e., $|\theta|$ is negligible relative to $\sigma$), we can see the estimators (\ref{SME1})-(\ref{SME1b}) as approximate Method of Moment Estimators. The assumption of $\theta\approx{}0$ has been suggested by some empirical literature (e.g., \cite{Seneta}, who in turns cites \cite{Hurst}). Using MME and MLE and intraday high-frequency data, this was also validated by \cite{FLLLM} for NIG and VG models. In the latter framework, we can perform a simple experiment to assess this assumption. From the formulas for $\mu_{2}$ and $\mu_{3}$ in (\ref{MEVGM}) as well as the formula for $c_{3}(\tau_{1})$, we have that 
\[
	\frac{|\mu_{3}(X_{\delta})|}{2\mu_{2}(X_{\delta})}\geq{}|\theta| \kappa\geq{}\theta^{2}\kappa,
\]
assuming that, as it is usually the case, $|\theta|\leq{}1$. Therefore, 
\[
	\frac{\sigma^{2}}{\theta^{2}\kappa}\geq{}\frac{2\mu_{2}(X_{\delta})^{2}}{\delta|\mu_{3}(X_{\delta})|}-1.
\]
The following table reports the values of $\frac{2\hat{\mu}_{2}(X)^{2}}{\delta|\hat{\mu}_{3}(X)|}-1$ for a few stocks. Thus, for instance, the value of 44 for 1 minute INTEL data suggests that $\sigma^{2}$ is at least 44 times larger than $\theta^{2}\kappa$ and thus, we can assume that $\mu_{2}(X_{\delta})\approx \sigma^{2}\delta$. One can do a similar analysis to justify that $\mu_{4}(X_{\delta})\approx 3\sigma^{4}\kappa\delta$.}
\begin{table}[ht]
\footnotesize
%\tiny
\centering\addtolength{\tabcolsep}{-2pt}{
\begin{tabular}{|c|c|c|c|c|c|c|c|}
\hline $\delta$ & 5 sec	&	10 sec	&	30 sec	&	1 min 	&	5 min 	&	10 min &	30 min 
\tabularnewline \hline INTEL & $144$   & $82.4$    & $57$   & $44$    & $26.7$    & $24.7$ & $13$
\tabularnewline CVX &  $3146.8$  & $3023.8$    & $8706.9$   & $212.9$    & $251.0$    & $1231.5$    & $175.3$
\tabularnewline CSCO & $587.5$  &  $255.8$  &  $94.1$  &  $77.5$  &  $67.1$  &  $52.3$   & $37.6$
\tabularnewline PFE & $47.8$  & $24.2$ &  $10.7$  & $7.89$ &  $7.67$  & $7.63$  & $8.15$
\tabularnewline \hline
\end{tabular}
}
\vspace{.2 cm}
\caption{Computation  of  $\frac{2\hat{\mu}_{2}(X)^{2}}{\delta|\hat{\mu}_{3}(X)|}-1$ for different stocks {based on high-frequency data during the year of 2005 ($T=252$ days)}.}
\label{Table1}	
\end{table}
\end{remk}

\subsection{Simple infill properties in the absence of noise}
We now proceed to show some ``in-fill" ($n\to\infty$ with fixed $T$) asymptotic properties of the estimators  in (\ref{SME1})-(\ref{SME1b}). 
As above, in the sequel we assume that $\theta=0$ and neglect $O(\delta_{n})=O(1/n)$ terms. In that case, it is easy to see that 
\begin{align}\label{ABSMNMSE}
	{\bbe\hat\sigma_{n}^{2}=\bbe\tilde\sigma_{n}^{2}=\sigma^{2}+O\left(\frac{1}{n}\right),\qquad
	{\rm Var}\left(\hat\sigma_{n}^{2}\right)={\rm Var}\left(\tilde\sigma_{n}^{2}\right)=\frac{3\sigma^{2}\kappa}{T}+O\left(\frac{1}{n}\right)}.
\end{align}
From the above formulas, we conclude the (not surprising) fact that, {on a finite time horizon},  $\hat\sigma_{n}^{2}$ is not a mean-squared consistent estimator for $\sigma^{2}$, when the sampling frequency increases, {but the MSE is of order $O(1/T)$, as $T\to{}\infty$.}%However, the standard error of $\hat{\sigma}_{n}$ decreases inversely proportional to the time horizon $T$.

An analysis of the bias and variance of $\hat\kappa_{n}$ and $\tilde{\kappa}_{n}$ is more complicated due to the non-linearity of the sample kurtosis. However, we can deduce some interesting features of its infill asymptotic behavior.  {First, we have
\begin{align}\label{DfCkrt}
\underset{n\to\infty}{{\rm lim^{P}}}\;\hat{\kappa}_{n}=\underset{n\to\infty}{{\rm lim^{P}}}\;\tilde\kappa_{n}=\frac{1}{3} \frac{\frac{1}{T}\sum_{t\leq{}T}\left(\Delta X_{t}\right)^{4}}{\left(\frac{1}{T}\sum_{t\leq{}T}\left(\Delta X_{t}\right)^{4}\right)^{2}}
	=: \hat\kappa^{(T)},
\end{align}
where above $\Delta X_{t}=X_{t}-X_{t^{-}}$ is the jump size of $X$ at time $t$ and the summations are over the random countable set of times $t$ for which $\Delta X_{t}\neq{}0$. The limit (\ref{DfCkrt}) follows from the well-known formula $\sum_{i=1}^{n} \left(X_{i\delta_{n}}-X_{(i-1)\delta_{n}}\right)^{k}\cp \sum_{t\leq{}T}(\Delta X_{t})^{k}$, as $n\to\infty$, valid for any $k\geq{}2$ and a \emph{pure-jump} L\'evy process $X$. Furthermore, the convergence of the corresponding moments also holds true since 
$0\leq\delta_{n}\hat\mu_{4,n}/\hat\mu_{2,n}^{2}\leq \delta_{n} n=T<\infty$, and, thus, 
\begin{align}\label{FSA}
	\lim_{n\to{}\infty}\bbe \hat{\kappa}_{n}=\lim_{n\to\infty}\bbe \tilde{\kappa}_{n}= \bbe \hat\kappa^{(T)}\quad
	{\rm and}\quad 
	\lim_{n\to{}\infty}{\rm Var}\left(\hat{\kappa}_{n}\right)=\lim_{n\to\infty}{\rm Var}\left(\tilde{\kappa}_{n}\right)= {\rm Var}\left(\hat\kappa^{(T)}\right).
\end{align}
The following result, whose proof is given in the Appendix, expands the expectation and variance of $\hat\kappa^{(T)}$ above and shows that the MSE of $\hat\kappa^{(T)}$ is $O(T^{-1})$, as $T\to\infty$.}
\begin{prop}\label{BEK}
	{Let $X$ be a general L\'evy process with L\'evy measure $\nu$.} Let $c_{i}:=c_{i}(X_{1})$ be the $i^{th}$ cumulant of $X_{1}$, $\kappa:=c_{4}/3c_{2}^{2}$, and suppose that $\int |x|^{i}\nu(dx)<\infty$ for any $i\geq{}2$. Then, as $T\to{}\infty$, 
	\begin{align}\label{BsEstKapCnt}
	&\bbe\, {\hat\kappa}^{(T)}={\kappa}
	+\frac{3c_{4}^{2}-2c_{6}c_{2}}{3c_{2}^{4}}T^{-1}+O(T^{-2}),\\
	&\bbe\left( {\hat\kappa}^{(T)}-{\kappa}\right)^{2}=
	\frac{c_{8}c_{2}-4c_{4}c_{6}+4c_{4}^{2}c_{2}}{9c_{2}^{5}}T^{-1}+O(T^{-2}).\label{MSEEstKapCnt}
\end{align}
\end{prop}

\subsection{Properties of {the MME} under microstructure noise}\label{MNSct}

In this part we {characterize the effects of a microstructure noise component into the asymptotic properties of the MME  introduced above. The results for the case of the volatility estimators are classical and their proofs are given only for the sake of completeness. The results for the estimators of the kurtosis parameter $\kappa$ are not hard to get either but are less known.}

We adopt the setup introduced at the end of Section \ref{SctModels}, under which the observed log-returns are given by
\begin{equation}\label{MILMN}
{\Delta_{i}^{n}\widetilde{X} := \widetilde{X}_{i \delta_n} - \widetilde{X}_{(i-1) \delta_n}  = \left(X_{i \delta_n}-X_{(i-1) \delta_n}\right)+\left(\varepsilon_{i \delta_n}-\varepsilon_{(i-1) \delta_n}\right)  =:\Delta_{i}^{n} X+\tilde{\varepsilon}_{i,n}}.
\end{equation}
Furthermore, throughout we assume that, for each $n$, $(\tilde{\varepsilon}_{i,n})_{i\geq{}1}$ {satisfies the following mild assumption, for any positive integer $k\geq{}1$:
\begin{align}\label{SAN1}
	\frac{1}{n} \sum_{i=1}^{n} (\tilde{\varepsilon}_{i,n})^{k}\cp m_{k}(\tilde{\varepsilon}), \quad (n\to\infty), \quad \text{for some }{m}_{k}(\tilde\varepsilon)\in\bbr.
\end{align}
Obviously, the previous assumption 
covers the microstructure white-noise case, 
where $(\varepsilon_{t})_{t\geq{}0}$ are i.i.d.,
in which case} ${m_{k}(\tilde\varepsilon):=\bbe\left(\left(\tilde{\varepsilon}_{1,n}\right)^{k}\right)}$.
{Note that $\tilde{\varepsilon}$ is not required to be independent of the process $X$ and, furthermore, we only need for $X$ to be a pure-jump semimartingale.}

Let us first {{describe} the infill asymptotic behavior of the estimators for $\sigma^{2}$, {introduced in}  {(\ref{SME1})-(\ref{SME1b})}, but based on the noisy observations:}
\begin{equation}\label{NEFV}
	{\tilde{\sigma}_{n}^2(\widetilde{X})}:=\frac{1}{\delta_{n}n} \sum_{i=1}^{n}(\Delta_{i}^{n}\widetilde{X} -	\overline{\Delta^{n}\widetilde{X} })^{2},\qquad 
	{\hat{\sigma}_{n}^2(\widetilde{X})}:=\frac{1}{T}{\widehat{\big[\widetilde{X},\widetilde{X}\big]}_{2}}=\frac{1}{\delta_{n}n} \sum_{i=1}^{n}(\Delta_{i}^{n}\widetilde{X})^{2}.
\end{equation}
{For future reference, let us state the following {simple} result that follows from applying Cauchy's inequality, the condition (\ref{SAN1}), and the fact that $\sum_{i=1}^{n}|\Delta_{i}^{n} X|^{2m}\stackrel{\mathbb{P}}{\to}\,\sum_{s\leq{}T}|\Delta X_{s}|^{2m}$}. 
\begin{lmma}\label{L1NMN}
For {arbitrary integers $m\geq{}1$ and $k\geq{}0$,}
\begin{equation}\label{ECET}
	\frac{1}{n}  \sum_{i=1}^{n} (\Delta_i^n X)^m (\tilde{\varepsilon}_{i,n})^k \cp 0,\quad \text{as } n\to\infty.
\end{equation}
\end{lmma}

{We are now ready to analyze the asymptotic behavior of the {estimators in} (\ref{NEFV}). {The following result gives the in-fill asymptotic behavior of $\hat{\sigma}_{n}^2(\widetilde{X})$ and $\tilde{\sigma}_{n}^2(\widetilde{X})$.}
 \begin{prop}\label{PrpABSNO}
	Both estimators  $\hat{\sigma}_{n}^2(\widetilde{X})$ and $\tilde{\sigma}_{n}^2(\widetilde{X})$ admit the decomposition 
	\[
		\hat{\sigma}_{n}^2(\widetilde{X})=A_{n}+B_{n},\quad \tilde{\sigma}_{n}^2(\widetilde{X})=\widetilde{A}_{n}+\widetilde{B}_{n}
	\]
	where the r.v.'s above are such that%, as $n\to\infty$, 
	\[
		\underset{n\to\infty}{{\rm lim^{P}}}A_{n}=\underset{n\to\infty}{{\rm lim^{P}}} \widetilde{A}_{n}= \frac{1}{T}\sum_{s\leq{}T}(\Delta X_{s})^{2}, \quad
		{\underset{n\to\infty}{{\rm lim^{P}}}\delta_{n}{B}_{n}=m_{2}(\tilde\varepsilon)},\quad
		{\underset{n\to\infty}{{\rm lim^{P}}}\delta_{n}\widetilde{B}_{n}=m_{2}(\tilde\varepsilon)-(m_{1}(\tilde\varepsilon))^{2}}.
	\]
\end{prop}

\noindent\textbf{Proof}. {We only give the proof for {$\tilde{\sigma}_{n}^2:=\tilde{\sigma}_{n}^2(\widetilde{X})$}. The proof for {$\hat{\sigma}_{n}^2(\widetilde{X})$} is identical.} 
{First} note that
\begin{align*}
{\tilde{\sigma}_{n}^2} 
                &= 
                \frac{1}{n\delta_{n}} \sum_{i=1}^{n}(\Delta_{i}^{n}X-\overline{\Delta^{n}X})^2  + \frac{1}{n\delta_{n}}\sum_{i=1}^{n} 
                {(\tilde{\varepsilon}_{i,n}-\overline{\tilde{\varepsilon}_{n}})^2}+
                \frac{2}{n\delta_{n}} \sum_{i=1}^{n} (\Delta_{i}^{n}X-\overline{\Delta^{n}X})(\tilde{\varepsilon}_{i,n}-\overline{\tilde{\varepsilon}_{n}})\\
                &=:{\widetilde{A}_{n}+\widetilde{B}_{n,1}+\widetilde{B}_{n,2}}.
\end{align*}
The term ${\widetilde{A}_n}$ converges to $T^{-1}\sum_{s\leq{}T}(\Delta X_{s})^{2}$, as $n \to \infty$, since $\sum_{i=1}^{n}(\Delta_{i}^{n}X)^2\to\sum_{s\leq{}T}(\Delta X_{s})^2$ and  $\overline{\Delta^{n}X}=O_{P}(1/n)$. 
 Clearly, (\ref{SAN1}) implies that 
$\delta_{n}{\widetilde{B}_{n,1}}=n^{-1}\sum_{i=1}^{n} (\tilde{\varepsilon}_{i,n})^{2}-\left(\overline{\tilde{\varepsilon}_{n}}\right)^{2}$ converges to $m_{2}(\tilde\varepsilon)-(m_{1}\left(\tilde\varepsilon)\right)^{2}$, in probability, when $n\to\infty$.
Also, using Lemma \ref{L1NMN}, 
{$\delta_{n}{\widetilde{B}_{n,2}}=
	\frac{2}{n}  \sum_{i=1}^{n} (\Delta_{i}^{n}X) (\tilde{\varepsilon}_{i,n}) - 2 \overline{\Delta^n X} \,\overline{\tilde{\varepsilon}_{n}}$} goes to $0$ in probability. 
	 \hfill$\Box$

{Next,} let us consider the estimators for $\kappa$ introduced in {(\ref{SME1})-(\ref{SME1b})}, but applied to the noisy process $\widetilde{X}$:
\[
	\tilde{\kappa}_{n}(\widetilde{X}) = 
	\frac{\delta_{n}}{3}\left( \frac{\hat{\mu}_{4,n}(\widetilde{X})}{\hat{\mu}_{2,n}^2(\widetilde{X})}-3\right),\quad\quad
	\hat{\kappa}_{n}(\widetilde{X})
	:=
	\frac{T}{3}\frac{{\widehat{\big[\widetilde{X},\widetilde{X}\big]}_{4}}}{{\widehat{\big[\widetilde{X},\widetilde{X}\big]^{2}_{2}}}}.
\]
The following result states that, for {large $n$, the above estimators behave asymptotically as $\delta_n C$, for some constant $C$, depending on the ergodic properties of the microstructure noise.}
\begin{prop}\label{PrpABKNO}
	There exist non-zero constants $C$ and $\widetilde{C}$ such that, as $n\to\infty$, 
	\begin{equation}\label{AsymtkEsts}
		\frac{1}{\delta_{n}}\hat\kappa_{n}(\widetilde{X})\cp C, \quad\quad\quad \frac{1}{\delta_{n}}\tilde\kappa_{n}(\widetilde{X})\cp \widetilde{C}.%:=\frac{m_{4}(\tilde\varepsilon)}{3m_{2}^{2}(\tilde\varepsilon)}
	\end{equation}
\end{prop} 
 \noindent\textbf{Proof}.
 We only give the proof for ${\tilde{\kappa}_{n}:=\tilde{\kappa}_{n}(\widetilde{X})}$. The proof for {$\hat{\kappa}_{n}(\widetilde{X})$} is similar. {First, 
observe} that
\begin{align*}
	\hat{\mu}_{2,n}(\widetilde{X})
 &= \frac{1}{n} \sum_{i=1}^{n} (\Delta_{i}^{n} X -\overline{\Delta^n X})^2 
 + \frac{1}{n} \sum_{i=1}^{n} (\tilde{\varepsilon}_{i,n} - \bar{\tilde{\varepsilon}}_{n})^2+ \frac{2}{n} \sum_{i=1}^{n} (\Delta_{i}^{n} X -\overline{\Delta^n X})(\tilde{\varepsilon}_{i,n} - \bar{\tilde{\varepsilon}}_{n}).  
\end{align*}
By Lemma \ref{L1NMN}, the first and third terms on the last expression above tend to $0$ in probability, while the second term converges to ${\widetilde{C}_{0}}:=\m_{2}(\tilde{\varepsilon})-(m_{1}(\tilde{\varepsilon}))^{2}$ by (\ref{SAN1}).
Similarly,
\begin{align*}
\hat{\mu}_{4,n}(\widetilde{X})
 &= %\frac{1}{n} \sum_{i=1}^{n} (\Delta_{i}^{n} X -\overline{\Delta^n X})^4
  {\frac{1}{n} \sum_{\ell=0}^{3}\binom{4}{\ell}\sum_{i=1}^{n} (\Delta_{i}^{n} X -\overline{\Delta^n X})^{4-\ell}(\tilde{\varepsilon}_{i,n} - \bar{\tilde{\varepsilon}}_{n})^{\ell}
+ \frac{1}{n} \sum_{i=1}^{n} (\tilde{\varepsilon}_{i,n} - \bar{\tilde{\varepsilon}}_{n})^4}, 
\end{align*}
and, again, by Lemma \ref{L1NMN}, all the terms {in the first summation above tend to $0$ in probability, while the second term therein converges to 
	$\widetilde{C}_{1}:=m_{4}(\tilde{\varepsilon})-4m_{3}(\tilde{\varepsilon})m_{1}(\tilde{\varepsilon})+6m_{2}(\tilde{\varepsilon})m_{1}^{2}(\tilde{\varepsilon})-3m_{1}^{4}(\tilde{\varepsilon})$,
	in light of our assumption (\ref{SAN1}).
Therefore, the second limit in (\ref{AsymtkEsts}) follows with $\widetilde{C}:=\widetilde{C}_{1}/3\widetilde{C}_{0}^{2}-1$.} 
\hfill$\Box$

\begin{remk}
	As a consequence of the proof, it follows that, if {$m_{1}(\tilde{\varepsilon})=0$}, then 
	\[
		C=\widetilde{C}=\frac{m_{4}(\tilde\varepsilon)}{3\left(m_{2}(\tilde\varepsilon)\right)^{2}}.
	\]
	In particular, if the microstructure noise {$(\varepsilon_{t})_{t\geq{}0}$} in (\ref{OPDfn}) is white-noise, then the constant coincides with the excess kurtosis, {$\bbe\tilde{\varepsilon}^{4}/3\left(\bbe\tilde{\varepsilon}^{2}\right)^{2}$,} of the random variable $\tilde{\varepsilon}:=\varepsilon_{2}-\varepsilon_{1}$.
\end{remk}

\section{{Robust} Method of Moments Estimators}\label{MMECorrected}

In this section, we adapt the so-called {two-scales} bias correction technique of \cite{ZhaMykAit:2005} to develop estimators for $\sigma^{2}$ and $\kappa$ that are robust against microstructure noises.  Roughly, their approach consists of three main ingredients: sparse subsampling, averaging, and bias correction.  
Let us first introduce some needed notation. Let  $\bar{\mathcal{G}}_{n}:=\{t_{0},t_{1},\dots,t_{n}\}$ be {the} complete set of available sampling times {as described in (\ref{SmplSchm})}. For a subsample ${\mathcal{G}}=\{t_{i_{1}},\dots,t_{i_{m}}\}$ with $i_{1}\leq\dots\leq i_{m}$ and {a} natural $\ell\in\bbn$, we define the $\ell^{th}$-order realized variation of the process $\widetilde{X}$ over $\mathcal{G}$ as 
\[
	[\widetilde{X},\widetilde{X}]_{\ell}^{\mathcal{G}}=\sum_{j=0}^{m-1}\left|\widetilde{X}(t_{i_{j+1}})-\widetilde{X}(t_{i_{j}})\right|^{\ell}.
\]
Next, we partition the  grid $\bar{\mathcal{G}}_{n}$ into $K$ mutually exclusive regular sub grids as follows:
\[
	\mathcal{G}^{(i)}_{n}:=\mathcal{G}^{(i)}_{n,K}:=\{t_{i-1},t_{i-1+K},t_{i-1+2K}, \dots, t_{i-1+n_{i}K}\}, \qquad i=1,\dots, K,
\]
with $n_{i}:=n_{i,K}:=[(n-i+1)/K]$. As in \cite{ZhaMykAit:2005}, the key idea to improve the estimators introduced in (\ref{SME1b}) consists of averaging the relevant realized variations over the different sparse sub grids $\mathcal{G}^{(i)}_{n}$, instead of using only one realized variation over the complete set $\bar{\calG}_{n}$. Hence, for instance, for estimating $\sigma^{2}$, we shall consider the estimator 
\begin{equation}\label{FrstEstSigma}
	\hat{\sigma}^{2}_{n}:=\hat\sigma^{2}_{n,K}:=\frac{1}{K}\sum_{i=1}^{K}{\frac{1}{T_{i,K}}}[\widetilde{X},\widetilde{X}]_{2}^{\mathcal{G}_{n}^{(i)}},
\end{equation}
{where $T_{i,K}:=t_{i-1+n_{i}K}-t_{i-1}=K\delta_{n}n_{i}$. The estimator (\ref{FrstEstSigma})} is constructed by averaging estimators of the form $\hat{\sigma}^2(\widetilde{X})$ in (\ref{NEFV}) over sparse sub-grids.
The above estimator corresponds to the so-called ``second-best  estimator" in \cite{ZhaMykAit:2005}. {This estimator} can be improved in two ways. First, by correcting the bias of the estimator and, second, by choosing the number of sub grids, $K$, in an ``optimal" way. We analyze these two approaches in the subsequent two subsections. 

At this point it is convenient to recall that we are assuming the subordinated Brownian motion model (\ref{VGM}) with $\theta=0$. For simplicity, we also assume that $b=0$, which won't affect much what follows since we are considering high-frequency type estimators and, thus, the contribution of the drift is negligible in that case. {Regarding the microstructure noise, we assume that the noise process $\{{\varepsilon}_{t}\}_{t\geq{}0}$ appearing in Eq.~(\ref{OPDfn}) is a centered stationary process with finite moments of arbitrary order, independent of $X$. Furthermore, we assume that, for any $\ell\in\bbn$,\begin{equation}\label{MCMNT}
	\lim_{\delta_{1},\delta_{2}\to0}\frac{\bbe\left[\tilde\varepsilon_{\delta_{2}}^{\ell}\right]-\bbe\left[\tilde\varepsilon_{\delta_{1}}^{\ell}\right]}{\delta_{2}-\delta_{1}}=0;
\end{equation}
where hereafter $\tilde{\varepsilon}_{\delta}$ denotes a random variable with the same distribution as $\tilde{\varepsilon}_{t,\delta}:=\varepsilon_{t+\delta}-\varepsilon_{t}$, which does not depend on $t$. Note that (\ref{MCMNT}) implies the existence of a constant  ${m}_{\ell}(\tilde{\varepsilon})\in\bbr$ such that 
\begin{equation}\label{MCMNTC}
	\lim_{\delta\to{}0}\bbe\left[\tilde\varepsilon_{\delta}^{\ell}\right]= {m}_{\ell}(\tilde\varepsilon).
\end{equation}
The simplest case is the white noise, when the variables $\{\varepsilon_{t}\}_{t\geq{}0}$ are independent identically distributed. In that case,  $\{\tilde{\varepsilon}_{i\delta,\delta}\}_{i\geq{}1}$ follows a stationary Moving Average (MA) process with $\bbe\left(\tilde{\varepsilon}_{i\delta,{\delta}}\right)=0$ and $\bbe\big(\tilde{\varepsilon}_{i\delta,{\delta}}^{\,2}\big)=2\bbe\left(\varepsilon_{1}^{2}\right)$.}

\subsection{Bias corrected estimators}
{In order to deduce the bias correction, we first adopt the white noise case, where $\{\varepsilon_{t}\}_{t\geq{}0}$ are i.i.d. In that case, the distribution of $\tilde{\varepsilon}_{\delta}$ does not depend on $\delta$. A random variable with this distribution is denoted $\tilde{\varepsilon}$.}
We start by devising bias correction techniques for the estimator (\ref{FrstEstSigma}).  Clearly, from (\ref{MEVGM}) and the independence of the noise $\tilde{\varepsilon}$ and {the process} $X$, we have:
\begin{align}
	\bbe\left({\hat{\sigma}^{2}_{n,K}}\right)
%	=\bbe\left(\frac{1}{K}\sum_{i=1}^{K}T^{-1}[\widetilde{X},\widetilde{X}]_{2}^{\mathcal{G}_{n}^{(i)}}\right)
	=\sigma^{2}
	+ \bbe \left({\tilde{\varepsilon}^{\,2}}\right)\frac{1}{K}\sum_{i=1}^{K}\frac{n_{i}}{T_{i,K}}
%	+O\left(\frac{KT}{n}\right)
%	\\
	=\sigma^{2}
	+ \bbe \left({\tilde{\varepsilon}^{\,2}}\right){\frac{1}{K\delta_{n}}}.\label{BiasSE0}
\end{align}
The relation (\ref{BiasSE0}) shows that the bias of the estimator diverges to infinity when the time span between observation $\delta_{n}:=T/n$ tends  to $0$. To correct this issue, first note that (\ref{BiasSE0}) also implies that
\begin{equation}\label{CstEstm2}
	\bbe\left(\delta_{n}{\hat{\sigma}^{2}_{n,1}}\right)
	=\sigma^{2}\delta_{n}
	+ \bbe \left({\tilde{\varepsilon}^{\,2}}\right)\,\stackrel{n\to\infty}{\longrightarrow}\,\bbe \left({\tilde{\varepsilon}^{\,2}}\right).
\end{equation}
Hence, a natural ``bias-corrected" estimator would be 
\begin{align}\label{ZMAS0}
	\hat{\tilde{\sigma}}^{2}_{n}:=\hat{\tilde{\sigma}}^{2}_{n,K}:={\hat{\sigma}^{2}_{n,K}}-{\frac{1}{K\delta_{n}}}\hat\mu_{2,n}(\tilde\varepsilon),
\end{align}
where $\hat\mu_{2,n}(\tilde\varepsilon):=\delta_{n}\hat\sigma_{n,1}^{2}$.
However, from (\ref{BiasSE0}) with $K=1$, we have:
\begin{align*}
	\bbe\left(\hat{\tilde{\sigma}}^{2}_{n}\right)=\sigma^{2}
	+ \bbe \left(\tilde{\varepsilon}^{\,2}\right)\frac{n}{KT}-\frac{1}{K}\left(\sigma^{2}
	+ \bbe \left(\tilde{\varepsilon}^{\,2}\right)\frac{n}{T}\right)=\frac{K-1}{K}\sigma^{2},
\end{align*}
which {implies} that $\hat{\tilde{\sigma}}_{n}^{2}$ is not truly unbiased. 
Nevertheless, the above relationship yield the following unbiased estimator for $\sigma^{2}$:
\begin{align}\label{ZMAS0b}
	\hat{\bar{\sigma}}^{2}_{n,K}:=\frac{K}{K-1}\hat{\tilde{\sigma}}_{n,K}=\frac{1}{K-1}\sum_{i=1}^{K}\frac{1}{T_{i,K}}[\widetilde{X},\widetilde{X}]_{2}^{\mathcal{G}_{n}^{(i)}}-\frac{1}{(K-1)T}[\widetilde{X},\widetilde{X}]_{2}^{\bar{\mathcal{G}}_{n}}.
\end{align}
The estimator (\ref{ZMAS0b}) corresponds to the small-sample adjusted ``First-Best Estimator" of \cite{ZhaMykAit:2005}. 
\begin{prop}\label{PropEstmUbsSigma}
Under a centered stationary noise process $\{\varepsilon_{t}\}_{t\geq{}0}$ independent of $X$, 
\begin{align*}
	\bbe\left(\hat{\bar{\sigma}}^{2}_{n,K}\right)
	&=\sigma^{2}+\frac{\bbe\left(\tilde{\varepsilon}^{2}_{K\delta_{n}}\right)-\bbe\left(\tilde{\varepsilon}^{2}_{\delta_{n}}\right)}{\delta_{n}(K-1)}.
\end{align*}
In particular, $\hat{\bar{\sigma}}^{2}_{n,K}$ is an asymptotically unbiased (respectively, unbiased) estimator for $\sigma^{2}$ under the condition (\ref{MCMNT}) (respectively, a white microstructure noise setting). 
\end{prop}

We %next attempt to 
now devise (approximate) bias-corrected estimators for $\kappa$. In order to separate the problem of estimating $\kappa$ and $\sigma^{2}$, in this part we assume that $\sigma$ is known. In practice, we have to replace $\sigma$ with an ``accurate" estimate such as the estimator (\ref{ZMAS0b}).  Let us start by considering the mean of the statistic 
\[
\frac{1}{K}\sum_{i=1}^{K}\frac{1}{T_{i,K}}\left[\widetilde{X},\widetilde{X}\right]_{4}^{\mathcal{G}^{(i)}_{n}},
\]
which is the analog of (\ref{FrstEstSigma}). To this end, we use the fact that 
%\begin{align*}
	$\bbe\left(X_{\delta}+\tilde{\varepsilon}\right)^{4}
%	&=\bbe\left(X_{\delta}^{4}\right)
%	+4\bbe\left(X_{\delta}^{3}\,\tilde{\varepsilon}\right)
%	+6\bbe\left(X_{\delta}^{2}\,\tilde{\varepsilon}^{2}\right)
%	+4\bbe\left(X_{\delta}\,\tilde{\varepsilon}^{3}\right)+4\bbe\left(\tilde{\varepsilon}^{4}\right)\\
	=3\sigma^{4}\kappa \delta+6\sigma^{2}\bbe\left(\tilde{\varepsilon}^{\,2}\right) \delta+\bbe\left(\tilde{\varepsilon}^{\, 4}\right)+3\sigma^{4}\delta^{2}$, which is an easy consequence of (\ref{MEVGM}) and  the independence of the noise $\tilde\varepsilon$ and ${X}$. In that case, we have 
\begin{align}\label{QuartVarMean0}
	\bbe\left(\frac{1}{K}\sum_{i=1}^{K}\frac{1}{T_{i,K}}[\widetilde{X},\widetilde{X}]_{4}^{\mathcal{G}^{(i)}_{n}}\right)=3\sigma^{4}\kappa+6\sigma^{2}\bbe\left(\tilde{\varepsilon}^{\, 2}\right)+\frac{1}{K\delta_{n}}\bbe\left(\tilde{\varepsilon}^{\,4}\right)+3\sigma^{4}K\delta_{n},
\end{align}
This identifies the estimator 
\begin{equation}\label{FrstEstKappa}
	\hat{\kappa}_{n,K}:=\frac{1}{3\sigma^{4}K}\sum_{i=1}^{K}\frac{1}{T_{i,K}}[\widetilde{X},\widetilde{X}]_{4}^{\mathcal{G}^{(i)}_{n}}-K\delta_{n},
\end{equation}
as an unbiased estimator for $\kappa$ in the absence of microstructure noise. However, as with the estimate of $\sigma$, the bias of the above estimate blows up when $\delta_{n}\to{}0$ due to the third term in (\ref{QuartVarMean0}). To correct this issue we need an estimate for $\bbe\left(\tilde\varepsilon^{4}\right)$, which can be inferred from the following limit
\begin{align}\label{QuartVarMean0b}
	\lim_{n\to\infty}\bbe\left(\frac{\delta_{n}}{T}[\widetilde{X},\widetilde{X}]_{4}^{\bar{\mathcal{G}}_{n}}\right)=\bbe\left(\tilde{\varepsilon}^{\,4}\right),
\end{align}
which is an easy consequence of (\ref{QuartVarMean0}) with $K=1$. 
Together with (\ref{CstEstm2}), these two suggests the following estimate:
\begin{equation}\label{EstmAlmstUbsKap}
	\hat{\tilde{\kappa}}_{n}:=\hat\kappa_{n,K}
%	\frac{1}{3\sigma^{4}K}\sum_{i=1}^{K}\frac{1}{T_{i,K}}[\widetilde{X},\widetilde{X}]_{4}^{\mathcal{G}^{(i)}_{n}}
	-\frac{2}{\sigma^{2}}\hat\mu_{2,n}(\tilde{\varepsilon})
	-\frac{1}{3\sigma^{4}K\delta_{n}}\hat\mu_{4,{n}}(\tilde{\varepsilon})
%	-K\delta_{n}
	,
%	\\
%	&=\frac{1}{K}\sum_{i=1}^{K}T^{-1}[\widetilde{X},\widetilde{X}]_{2}^{\mathcal{G}_{n}^{(i)}}-\frac{n-K+1}{KTn}[\widetilde{X},\widetilde{X}]_{2}^{\bar{\mathcal{G}}_{n}}.
\end{equation}
where 
\begin{equation}\label{FouthMnEst}
	\hat{\mu}_{2,n}(\tilde\varepsilon):=\frac{\delta_{n}}{T}[\widetilde{X},\widetilde{X}]_{2}^{\bar{\mathcal{G}}_{n}}, \quad \hat{\mu}_{4,n}(\tilde\varepsilon):=\frac{\delta_{n}}{T}[\widetilde{X},\widetilde{X}]_{4}^{\bar{\mathcal{G}}_{n}}.
\end{equation}
However, as with the estimator $\hat{\tilde{\sigma}}_{n}^{2}$ above, the above estimator is only asymptotically unbiased for large $n$ and $K$. The following result provides an unbiased estimator for $\kappa$ based on the realized variations of the process on two scales. The proof follows from (\ref{BiasSE0}) and (\ref{QuartVarMean0}) and is omitted. 
\begin{prop}\label{PropEstmUbsKap}
	Let 
	\begin{align}\label{EstmUbsKap}
	\hat{\bar{\kappa}}_{n}&
	:=%\hat{\bar{\kappa}}_{n,2}-\frac{T(n-K-1)(K-1)}{n(n+1)}\\
%	&=\hat{\bar{\kappa}}_{n,1}-\frac{2}{\sigma^{2}}\hat{m}_{2,n}(\tilde\varepsilon)-\frac{T(n-K-1)(K-1)}{n(n+1)}\\
%	&=\frac{nK}{(n+1)(K-1)}\hat{\bar{\kappa}}_{n,0}-\frac{2}{\sigma^{2}}\hat{m}_{2,n}(\tilde\varepsilon)-\frac{T(n-K-1)(K-1)}{n(n+1)}\\
%	&=\frac{nK}{(n+1)(K-1)}\left(\frac{1}{3\sigma^{4}K}\sum_{i=1}^{K}T^{-1}[\widetilde{X},\widetilde{X}]_{4}^{\mathcal{G}^{(i)}_{n}}-\frac{n-K+1}{3\sigma^{4}TK}\hat{m}_{4,n}(\tilde\varepsilon)\right)-\frac{2}{\sigma^{2}}\hat{m}_{2,n}(\tilde\varepsilon)-\frac{T(n-K-1)(K-1)}{n(n+1)}\\
	\frac{1}{3\sigma^{4}(K-1)}\sum_{i=1}^{K}\frac{1}{T_{i,K}}[\widetilde{X},\widetilde{X}]_{4}^{\mathcal{G}^{(i)}_{n}}-\frac{1}{3\sigma^{4}(K-1)T}[\widetilde{X},\widetilde{X}]_{4}^{\bar{\mathcal{G}}_{n}}\\
	&\quad-\frac{2}{n\sigma^{2}}[\widetilde{X},\widetilde{X}]_{2}^{\bar{\mathcal{G}}_{n}}-(K-1)\delta_{n}.\nonumber
\end{align}
Then, under a white microstructure noise independent of $X$, $\hat{\bar{\kappa}}_{n}$ is an unbiased estimator for $\kappa$. Furthermore, for a general centered stationary noise process, we have
\begin{align*}
	\bbe\left(\hat{\bar{\kappa}}_{n}\right)&=\kappa+\frac{2}{\sigma^{2}}\frac{K}{K-1}\left(\bbe\left(\tilde{\varepsilon}_{K\delta_{n}}^{\, 2}\right)-\bbe\left(\tilde{\varepsilon}_{\delta_{n}}^{\, 2}\right)\right)+\frac{1}{3\sigma^{4}}\frac{\bbe\left(\tilde{\varepsilon}_{K\delta_{n}}^{\, 4}\right)-\bbe\left(\tilde{\varepsilon}_{\delta_{n}}^{\, 4}\right)}{\delta_{n}(K-1)},
\end{align*}
which shows that $\hat{\bar{\kappa}}_{n}$ is asymptotically unbiased under condition (\ref{MCMNT}).
\end{prop}

\subsection{Optimal selection of $K$}
{In this part, given a specified function $b(K,n,T)$, ${O_{u}}(b(K,n,T))$ means that there exists a constant $c$, independent of $K$, $n$, and $T$, such that ${|{O_{u}}(b(K,n,T))|}\leq{} {c b(K,n,T)}$, for all $K$, $n$, and $T$. We also assume {the white-noise case where the microstructure noise $\{\varepsilon_{t}\}_{t\geq{}0}$} are centered i.i.d. r.v.'s.}

An important issue when using the {two-scales} procedure {described in the previous section} is the selection of the number of subclasses, $K$. A natural approach to deal with this issue consists of minimizing the variance of the relevant estimators over all $K$. This procedure will yield an optimal $K^{*}$ for the number of subclasses. %Once such a $K^{*}$ is selected, one could then correct the bias of the estimator using the techniques of the previous section.  
Let us first illustrate {this approach for} the estimator $\hat\sigma_{n,K}^{2}$ given in (\ref{FrstEstSigma}). The next result, whose proof is given in  Appendix \ref{PrfSec5}, gives the variance of $\hat{\sigma}_{n,K}^{2}$.
\begin{thrm}\label{VarSigmEst}
	The estimator (\ref{FrstEstSigma}) is such that 
	 \begin{align}\label{VarSE}
	 {\rm Var}\left(\hat\sigma^{2}_{n,K}\right)&=\frac{4\sigma^{4}K}{3n}+{\frac{4n\bbe\varepsilon^{4}}{K^{2}T^{2}}}+\frac{4\sigma^{4}}{3n}+\frac{3\sigma^{4}\kappa}{T}
+\frac{2\sigma^{4}}{3Kn}+\frac{8\sigma^{2}\bbe(\varepsilon^{2})}{KT}
%\frac{8n\left(3(\bbe\varepsilon^{2})^{2}+\bbe\varepsilon^{4}\right)}{KT^{2}}
\\
&\quad+{O_{u}}\left(\frac{K^{2}}{n^{2}}\right)+{O_{u}}\left(\frac{K}{Tn}\right)+{{O_{u}}\left(\frac{1}{KT^{2}}\right)}.\nonumber
	\end{align}
\end{thrm}
\begin{remk}
		As a consequence of (\ref{VarSE}), for a fixed arbitrary $K$ and a high-frequency/long-horizon sampling setup {($T_{n}\to\infty$ and $\delta_{n}=T_{n}/n\to{}0$)}, a sufficient asymptotic relationship between $T$ and $\delta_{n}$ for the estimator $\hat{\sigma}^{2}_{n,K}$ to be mean square consistent is that  {$\delta_{n}T_{n}\to{}\infty$}. {If} $K$ is chosen depending on $n$ and $T$, {as we intend to do next,} the feasible values $K:=K_{n,T}$ must be such that $K_{n,T}/n\to{}0$ and $n/{(K^{2}_{n,T}T^{2})}\to{}0$ {as} $T\to\infty$ and $\delta_{n}=T/n\to{}0$.% and .}
\end{remk}
Now, we are ready to propose an approximately ``optimal" $K^{*}$. To that end, let us first recall from (\ref{BiasSE0}) that the bias of the estimator is 
\begin{align}\label{BiasSE}
	{\rm Bias}\left(\hat\sigma^{2}_{n,K}\right)
	=2 \bbe\varepsilon^{2}{\frac{n}{TK}}.
\end{align}
Together (\ref{VarSE})-(\ref{BiasSE}) implies that
\begin{align}
	{\rm MSE}\left(\hat\sigma_{n,K}^{2}\right)&=\frac{4\sigma^{4}K}{3n}+\frac{4\sigma^{4}}{3n}+\frac{3\sigma^{4}\kappa}{T}
+\frac{2\sigma^{4}}{3Kn}+\frac{8\sigma^{2}\bbe(\varepsilon^{2})}{KT}+{\frac{4n\bbe\varepsilon^{4}}{K^{2}T^{2}}}+\frac{4 n^{2}\left(\bbe\varepsilon^{2}\right)^{2}}{T^{2}K^{2}}\label{MSETilSig}
\\
&\quad+O_{u}\left(\frac{K^{2}}{n^{2}}\right)+O_{u}\left(\frac{K}{Tn}\right)+{O_{u}\left(\frac{1}{KT^{2}}\right)}.\nonumber
\end{align}
Our goal is to minimize the MSE {with respect to $K$} when $n$ is large. Note that the only term that is increasing in $K$ is $4\sigma^{4}K/3n$, while out of the terms decreasing in $K$, {the term $4 n^{2}\left(\bbe\varepsilon^{2}\right)^{2}/T^{2}K^{2}$ is the dominant (when $n$ is large)}. It is then reasonable to consider only these two terms leading to the ``approximation":
\begin{equation}\label{MSE Exp1}
	{\rm MSE}\left(\hat\sigma_{K}^{2}\right)\approx \frac{4\sigma^{4}K}{3n}+\frac{4 n^{2}}{T^{2}K^{2}}(\bbe\varepsilon^{2})^{2}
%	+\frac{8 n}{K T^{2}}\left(\bbe\varepsilon^{4}+3\left(\bbe \varepsilon^{2}\right)^{2}\right)
	=:{\rm MSE}_{1}\left(\hat\sigma_{K}^{2}\right).
\end{equation}
The right-hand side in the above expression attains its minimum at the value:
\begin{equation}\label{Exp1OK}
	{K^{*}_{{1}}}=n\left(\frac{6(\bbe\varepsilon^{2})^{2}}{T^{2}\sigma^{4}}\right)^{\frac{1}{3}}.
\end{equation}
Interestingly enough, the value above coincides with the optimal $K^{*}$ proposed in \cite{ZhaMykAit:2005} (see Eq.~(8) therein). Plugging  (\ref{Exp1OK}) in (\ref{MSETilSig}) and, {since} $\delta=T/n\to{}0$, it follows that %(\ref{MSE Exp1}), we get 
\begin{equation}\label{MSEEST1}
	MSE\left(\hat\sigma^{2}_{{K^{*}_{1}}}\right)%=\left({\rm Bias}\left(\tilde\sigma^{2}_{K^{*}}\right)\right)^{2}
%	+{\rm Var}\left(\tilde{\sigma}_{K^{*}}^{2}\right)
	=2^{\frac{4}{3}}3^{\frac{1}{3}}\left(\bbe \varepsilon^{2}\right)^{\frac{2}{3}}\sigma^{\frac{8}{3}}T^{-\frac{2}{3}}+3\kappa\sigma^{4}T^{-1}+o(T^{-1}).%c_{2}(\bbe\varepsilon^{2})^{-\frac{2}{3}}\left(\bbe\varepsilon^{4}+3(\bbe\varepsilon^{2})^{2}\right)\sigma^{\frac{4}{3}}T^{-\frac{4}{3}}
\end{equation}
In particular, the above expression shows that, in the presence of {a microstructure noise component}, the rate of convergence reduces from $O(T^{-1})$ to only $O(T^{-2/3})$ {and, furthermore, that the convergence is worst when $\sigma$, $\bbe\varepsilon^{2}$, and $\kappa$ are larger.}

The following result {gives an estimate of} the variance of the unbiased estimator (\ref{ZMAS0b}). Its proof is given in Appendix \ref{PrfSec5}.
\begin{prop}\label{VarSigmEst2}
	The estimator (\ref{ZMAS0b}) is such that 
	 \begin{align}\label{VarSE2}
	 {\rm Var}\left(\hat{\bar{\sigma}}^{2}_{n,K}\right)&={\frac{4\sigma^{4}K}{3n}+\frac{4n\left({\bbe\varepsilon^{4}}+(\bbe\varepsilon^{2})^{2}\right)}{T^{2}K^{2}}+{O_{u}}\left(\frac{1}{n}\right)+{O_{u}}\left(\frac{n}{K^{3}T^{2}}\right)+{O_{u}}\left(\frac{1}{TK}\right).}
	\end{align}
\end{prop}
As before, the previous result suggests to fix $K$ so that {to minimize} the first two leading terms in (\ref{VarSE2}). Such a minimum  is given by 
\begin{equation}\label{Exp2OK}
	K^{*}_{{2}}=n^{\frac{2}{3}}\left(\frac{6\left(\bbe\varepsilon^{4}+(\bbe\varepsilon^{2})^{2}\right)}{T^{2}\sigma^{4}}\right)^{\frac{1}{3}},
\end{equation}
which is similar\footnote{{The optimal value of $K$ proposed in  \cite{ZhaMykAit:2005} lacks the term $\bbe\varepsilon^{4}$ in the numerator.}} (but not identical) to the analog optimal $K^{*}$ proposed in \cite{ZhaMykAit:2005} (see Eq.~(58) \& (63) therein). {After plugging {$K_{2}^{*}$} in (\ref{VarSE2}), the resultant estimator attains the MSE:
\begin{equation}\label{MSEEST2}
	MSE\left(\hat{\bar{\sigma}}^{2}_{K^{*}_{2}}\right)=2^{\frac{4}{3}}3^{\frac{1}{3}}\left(\bbe\varepsilon^{4}+(\bbe \varepsilon^{2})^{2}\right)^{\frac{1}{3}}\sigma^{\frac{8}{3}}n^{-\frac{1}{3}}T^{-\frac{2}{3}}+o(T^{-1}).
\end{equation}
Interestingly enough, since $T/n\to{}0$, the estimator $\hat{\bar{\sigma}}^{2}_{K^{*}_{1}}$ attains the order $o(T^{-1})$, which was not achievable by the estimators $\hat\sigma_{K}^{2}$, even in the absence of microstructure noise, nor by the standard estimators introduced in Section \ref{MMESec} (see (\ref{ABSMNMSE})).}

Now, we proceed to study the optimal selection problem of $K$ for the estimator (\ref{FrstEstKappa}) for $\kappa$. As with $\hat\sigma^{2}_{n,K}$, we first need to {analyze} the variance of the estimator.

\begin{thrm}\label{VarKappaEst}
	The estimator (\ref{FrstEstKappa}) is such that 
	 \begin{equation}\label{VarKappEq}
	 {\rm Var}\left(\hat\kappa_{n,K}\right)=\frac{{64}}{5}\frac{T^{2}K^{3}}{n^{3}}+
{O_{u}}\left(\frac{T^{2}K^{2}}{n^{3}}\right).
	\end{equation}
\end{thrm}

{We are now ready to} propose a method to choose a value of $K$ that approximately minimizes the MSE of the estimator $\hat\kappa_{n,K}$. Let us first recall from (\ref{QuartVarMean0}) that the bias of the estimator $\hat\kappa_{n,K}$ is 
\begin{align}\label{BiasKurt0}
	{\rm Bias}\left(\hat\kappa_{n,K}\right)=\bbe\left(\hat\kappa_{n,K}\right)-\kappa=\bbe\left(\tilde{\varepsilon}^{\,4}\right){\frac{n}{TK\sigma^{4}}}+2\frac{\bbe\left(\tilde{\varepsilon}^{\, 2}\right)}{\sigma^{2}}.
	%=\bbe\left(\tilde{\varepsilon}^{\,4}\right)\frac{n}{TK\sigma^{4}}+\text{{h.o.t.}}
\end{align}
Together, (\ref{VarKappEq})-(\ref{BiasKurt0}) imply that
\begin{align}
	{\rm MSE}\left(\hat\kappa_{n,K}\right)&=\frac{{64}}{5}\frac{T^{2}K^{3}}{n^{3}}+\frac{n^{2}\left(\bbe\tilde{\varepsilon}^{4}\right)^{2}}{T^{2}K^{2}\sigma^{8}}+\text{{h.o.t.},}\label{MSEKappa0}
\end{align}
{where {h.o.t. mean ``higher order terms"}.}
{It is then reasonable to select $K$ so that the leading terms of the MSE are minimized. The aforementioned minimum is reached at}
\begin{equation}\label{Exp1OK2}
	K^{*}_{3}
%	=n\left(\frac{5(9)(\bbe\tilde\varepsilon^{4})^{2}}{({288})(3)T^{4}\sigma^{8}}\right)^{1/5}
	=n\left(\frac{5(\bbe\tilde\varepsilon^{4})^{2}}{{96}T^{4}\sigma^{8}}\right)^{\frac{1}{5}}.
\end{equation}
Plugging  (\ref{Exp1OK2}) in (\ref{MSEKappa0}), it follows that %(\ref{MSE Exp1}), we get 
\begin{align*}
	MSE\left(\hat\kappa_{{K^{*}_{3}}}\right)
	={(4)5^{\frac{3}{5}}3^{-\frac{3}{5}}\left(\bbe \tilde\varepsilon^{4}\right)^{\frac{6}{5}}\sigma^{-\frac{24}{5}}T^{-\frac{2}{5}}+o\left(T^{-\frac{2}{5}}\right)},
%	=O\left(T^{-\frac{2}{5}}\right),%2^{\frac{8}{5}}5^{-\frac{2}{5}}3^{\frac{1}{5}}\left(\bbe \tilde\varepsilon^{4}\right)^{\frac{6}{3}}\sigma^{-\frac{24}{5}}T^{-\frac{2}{5}}+o\left(T^{-\frac{2}{5}}\right),
\end{align*}
whose rate of convergence to $0$ is slower than the rate of {$O\left(T^{-2/3}\right)$} attained by the estimator $\hat\sigma_{K^{*}}^{2}$. 

Finally, we consider the unbiased estimator for $\kappa$ introduced in Proposition \ref{PropEstmUbsKap}. {As above, h.o.t. refers to higher order terms}.
\begin{thrm}\label{ThrmVarKappaEst2}
	The estimator (\ref{EstmUbsKap}) is such that 
	 \begin{align}\label{VarKappaSE2}
	 	{\rm Var}\left(\hat{\bar{\kappa}}_{n,K}\right)=\frac{64}{5}\frac{T^{2}K^{3}}{n^{3}}+\frac{2n}{9\sigma^{8}T^{2}K^{2}}e(\varepsilon)+{{\rm h.o.t.}},
	\end{align}
	where $e(\varepsilon)={\rm Var}\left((\varepsilon_{2}-\varepsilon_{1})^{4}\right)$.
\end{thrm}
The two terms on the right-hand side of (\ref{VarKappaSE2}) reach their minimum value at 
\begin{equation}\label{Exp1OK3}
	K^{*}_{4}=n^{\frac{4}{5}}\left(\frac{5e(\varepsilon)}{(27)(16)T^{4}\sigma^{8}}\right)^{\frac{1}{5}}.
\end{equation}
{After plugging $K^{*}_{4}$ in (\ref{VarKappaSE2}), we obtain that 
\[
	MSE\left({\hat{\bar{\kappa}}_{{K^{*}_{4}}}}\right)
	={2^{\frac{28}{5}}5^{-\frac{2}{5}}3^{-\frac{9}{5}}e(\varepsilon)^{\frac{3}{5}}\sigma^{-\frac{24}{5}}n^{-\frac{3}{5}}T^{-\frac{2}{5}}+o\left(T^{-1}\right)},
\]
which again, since $T/n\to{}0$, implies that $MSE\left({\hat{\bar{\kappa}}_{{K^{*}_{4}}}}\right)=o(T^{-1})$. The aforementioned result should be compared to (\ref{MSEEstKapCnt}), which essentially says that the estimator {$\hat{\bar{\kappa}}_{K^{*}_{4}}$} has better efficiency than the continuous-time based estimator $\hat\kappa^{(T)}$, obtained by making $n\to\infty$ in the estimators $\hat\kappa_{n}$ and $\tilde\kappa_{n}$ (see (\ref{DfCkrt})).}
It is worth pointing out here that one can devise a consistent estimator for $e(\varepsilon)$ {using} the relationships
\begin{equation}\label{Estfouthadnsecmmnt}
	{\rm (i)}\;\frac{1}{n}[\widetilde{X},\widetilde{X}]_{4}^{\bar{\mathcal{G}}_{n}}\cp \bbe\left(\varepsilon_{2}-\varepsilon_{1}\right)^{4}, \quad 
	{\rm (ii)}\;
	\frac{1}{n}[\widetilde{X},\widetilde{X}]_{8}^{\bar{\mathcal{G}}_{n}}\cp \bbe\left(\varepsilon_{2}-\varepsilon_{1}\right)^{8}. 
\end{equation}

\begin{remk}
{It is natural to wonder if some types of central limit theorems are feasible for  the estimators considered here. In spite of the fact that we are considering a L\'evy model, whose increments are independent, the estimators cannot be written in terms of a row-wise independent triangular array. For instance, consider the estimator $\hat\sigma_{n,K}$ for $\sigma$ introduced in (\ref{FrstEstSigma}) and, for simplicity, assume that $T_{i,K}=T$, which asymptotically is satisfied, and absence of microstructure noise. It can be shown that 
\begin{equation*}
	\frac{1}{K}\sum_{i=1}^{K}[{X},{X}]_{2}^{\mathcal{G}_{n}^{(i)}}=
\frac{1}{K} \sum_{i= 0}^{n - K} (X_{t_{i + K}} - X_{t_{i}})^2,
\end{equation*}
whose terms are correlated.}

\end{remk}

\section{Numerical Performance and Empirical Evidence}\label{Numerics}

{In this section, we propose an iterative method to implement the estimators described in the previous section, with the corresponding optimal choices of $K^{*}$. The main issue arises from the fact that in order to accurately estimate $\sigma$, we need to choose $K$ as in (\ref{Exp2OK}) (or (\ref{Exp1OK})), which precisely depends on what we want to estimate, $\sigma$. So, we propose to start with an initial reasonable guess for $\sigma^{2}$ to find $K^{*}$, which in turn is then used to improve the initial guess of $\sigma$, and so forth. The finite-sample and empirical performance of the resulting estimators are illustrated by simulation and a real high-frequency data application.}
For briefness, {in what follows we will} make use of the following notation 
\[
 	\mathcal{K}^{*}_{1}(m_{2},\sigma):=n\left(\frac{6m_{2}^{2}}{T^{2}\sigma^{4}}\right)^{\frac{1}{3}},\qquad 
	\mathcal{K}^{*}_{2}(m_{2},m_{4},\sigma):=n^{\frac{2}{3}}\left(\frac{6\left(m_{4}+m_{2}^{2}\right)}{T^{2}\sigma^{4}}\right)^{\frac{1}{3}}.
\]
{For the simulation portion of this section,} we consider a Variance Gamma (VG) model with white Gaussian microstructure noise. The variance of the noise  $\varepsilon_{t}$ is denoted by $\varrho^{2}$ so that the noise of the $i^{th}$ increment, $\tilde{\varepsilon}_{i,n}$, is $\mathcal{N}(0,2\varrho^{2})$. Other parameters are set as:
%\[
	$\sigma=0.02$, $\kappa=0.3$, and $\varrho=0.005$. 
%\]
The time unit here is a day. In particular, the above value of $\sigma$ corresponds to an annualized volatility of $0.02\sqrt{252}=0.31$.

\subsection{Estimators for $\sigma$}\label{Sect:EstSigma}
	
We compare the finite sample performance of the following estimators:
\begin{enumerate}
	\item The estimator $\hat\sigma_{n,K}^{2}$ given in (\ref{FrstEstSigma}) with $K$ determined by a suitable estimate of the optimal value {$K^{*}_{1}$} given in (\ref{Exp1OK}), {as described next}.
%	\begin{equation}\label{Exp1OKb}
%	K^{*}=n\left(\frac{6(\bbe\varepsilon^{2})^{2}}{T^{2}\sigma^{4}}\right)^{1/3}.
%\end{equation}
As shown in Proposition \ref{PrpABSNO} {and (\ref{CstEstm2}), a consistent and unbiased} estimator 
%for $m_{2}(\tilde\varepsilon)=\bbe\tilde\varepsilon^{2}=2\bbe\varepsilon^{2}=2\varrho^{2}$ is {provided} by $\hat{m}_{2,n}(\tilde\varepsilon):=\delta_{n}\tilde{\sigma}_{n}^{2}=[\widetilde{X},\widetilde{X}]_{2}^{\bar{\mathcal{G}}_{n}}/n$, which suggests the following consistent estimate 
for {$\bbe\varepsilon^{2}=\bbe\tilde\varepsilon^{2}/2$} is given:
\begin{equation}\label{estVarRho}
	\hat\varrho^{2}:=\widehat{\bbe\varepsilon^{2}}:=\frac{1}{2n}[\widetilde{X},\widetilde{X}]_{2}^{\bar{\mathcal{G}}_{n}}.
\end{equation}
The only missing ingredient for estimating (\ref{Exp1OK}) is an initial preliminary estimate of $\sigma^{2}$, which we will then proceed to improve via $\hat\sigma_{n,K^{*}}^{2}$. Concretely,  we propose the following procedure. First, we evaluate the estimate
%\begin{align}\label{EstValKOpt}
${\hat{K}^{*}_{1}}:=\mathcal{K}_{1}^{*}(\hat\varrho^{2},\sigma_{0})$, %	{\hat{K}^{*}_{1}}:=n\frac{6^{1/3}\left(\widehat{\bbe\varepsilon^{2}}\right)^{2/3}}{T^{2/3}\sigma^{4/3}_{0}},
%\end{align}
 where $\sigma_{0}$ is {an initial} ``reasonable" value for the volatility. 
Second, we estimate $\sigma$ via $\hat\sigma'_{1}:=\hat\sigma_{n,\hat{K}^{*}_{1}}$. Next, we use $\hat\sigma'_{1}$ to improve our estimate of $K^{*}$ by setting 
%\begin{align}\label{EstValKOptb}
	$\hat{\hat{K}}^{*}_{{1}}:=\mathcal{K}_{1}^{*}(\hat\varrho^{2},\hat\sigma'_{1})$.
%	\hat{\hat{K}}^{*}_{{1}}:=n\frac{6^{1/3}\left(\widehat{\bbe\varepsilon^{2}}\right)^{2/3}}{T^{2/3}\left(\hat\sigma'_{1}\right)^{4/3}}
%\end{align}
Finally, we {set} $\hat\sigma''_{1}:=\hat\sigma_{n,\hat{\hat{K}}^{*}_{{1}}}$

\item We consider the bias-corrected estimator $\hat{\bar{\sigma}}^{2}_{n,K}$ introduced in (\ref{ZMAS0b}), with a value of $K$ given by $\hat{K}^{*}_{1}$ as defined in the point 1 above. We denote {this estimator} $\hat{\sigma}'_{2}$. We also analyze an iterative procedure similar to that in item 1, but using $\hat{{\sigma}}'_{2}$. Concretely, we set  {$\hat{\sigma}''_{2}=\hat{\bar{\sigma}}_{n,\hat{\bar{K}}^{*}_{{1}}}$,
where $\hat{\bar{K}}^{*}_{{1}}:=\mathcal{K}_{1}^{*}(\hat\varrho^{2},\hat\sigma'_{2})$.} 

\item Finally, we also consider the estimator $\hat{\bar{\sigma}}^{2}_{n,K}$ introduced in (\ref{ZMAS0b}) but using {an estimate of} the optimal value $K^{*}_{{2}}$ {as defined in Eq.~(\ref{Exp2OK})}. Concretely, we set {$\hat\sigma'_{3}=\hat{\bar{\sigma}}_{n,\hat{K}^{*}_{{2}}}$}  with 
%\[
	$\hat{K}^{*}_{{2}}:=\mathcal{K}^{*}_{2}(\hat\varrho^{2},\hat\varpi,\sigma_{0})$,
%	n^{2/3}\frac{6\Big[\widehat{\bbe\varepsilon^{4}}+\left(\widehat{\bbe\varepsilon^{2}}\right)^{2}\Big]^{1/3}}{T^{2/3}\sigma^{4/3}_{0}}
%\]
{where $\sigma_{0}$ is an initial reasonable value for $\sigma$ and} $\hat\varpi$ is a consistent estimator for ${\bbe\varepsilon^{4}}$. Next, we improve the estimate of $\hat\sigma'_{3}$ by setting 
\begin{equation}\label{BstEstKapp}
	\hat{\sigma}''_{3}:={\hat{\bar{\sigma}}_{n,\hat{\hat{K}}^{*}_{{2}}}},\quad\text{ with }\quad \hat{\hat{K}}^{*}_{{2}}:=\mathcal{K}_{2}^{*}(\hat\varrho,\hat\varpi,\hat\sigma_{3}').
%	n^{2/3}\left(\frac{6\left[\widehat{\bbe\varepsilon^{4}}+\left(\widehat{\bbe\varepsilon^{2}}\right)^{2}\right]}{T^{2}(\hat\sigma'_{3})^{4}}\right)^{1/3}
\end{equation}
To estimate ${\bbe\varepsilon^{4}}$, we use (\ref{FouthMnEst}). Concretely, as shown in the proof of Proposition \ref{PrpABKNO} {and also in Eq.~(\ref{QuartVarMean0b})}, the statistics $\hat{m}_{4,n}(\tilde\varepsilon):=[\widetilde{X},\widetilde{X}]_{4}^{\bar{\mathcal{G}}_{n}}/n$ converges to $\bbe\left(\tilde\varepsilon^{4}\right)=2\bbe\varepsilon^{4}+6\left(\bbe\varepsilon^{2}\right)^{2}$.
%\begin{equation}
%	\hat{m}_{4,n}(\tilde\varepsilon):=\frac{1}{n}[\widetilde{X},\widetilde{X}]_{4}^{\bar{\mathcal{G}}_{n}}\to m_{4}(\tilde\varepsilon):=\bbe\left(\tilde\varepsilon^{4}\right)=2\bbe\varepsilon^{4}+6\left(\bbe\varepsilon^{2}\right)^{2}.
%\end{equation}
Therefore, a consistent estimate for  ${\bbe\varepsilon^{4}}$ is given by 
\[
	\hat\varpi:=\widehat{\bbe\varepsilon^{4}}:=\frac{1}{2n}[\widetilde{X},\widetilde{X}]_{4}^{\bar{\mathcal{G}}_{n}}-{3}\left(\widehat{\bbe\varepsilon^{2}}\right)^{2}.
\]

\end{enumerate} 

The sample mean, standard deviation, and mean-squared error (MSE) based on $1000$ simulations are presented in the {Table \ref{Table2}}. Here, we take $T=252\text{ days}$ and $\sigma_{0}\approx{}0.063$, which corresponds to an annualized volatility of $1$. As expected, the estimator $\hat\sigma'_{1}$ exhibits a {noticeable bias and that this bias is corrected} by $\hat\sigma'_{2}$. However, $\hat\sigma''_{3}$ is much more superior to other considered estimators, {which is consistent with the asymptotic results for the mean-squared errors described in Eqs.~(\ref{MSEEST1}) and (\ref{MSEEST2}).}

\begin{table}[ht]
{\scriptsize%footnotesize
	\begin{tabular}{|c|c|c|c|c|c|c|c|}
	\hline
	 	$\delta_{n}$ & & $\hat\sigma'_{1}$ & $\hat\sigma''_{1}$ &  $\hat\sigma'_{2}$ & $\hat\sigma''_{2}$  & $\hat\sigma'_{3}$ & $\hat\sigma''_{3}$  
\\	\hline 
\multirow{3}{*}{$\text{5 min}$}& Mean & 
	0.02274333 & 0.02066226 & 0.01998258 & 0.01988843 & 0.01999695 & 0.01999614\\
	& Std Dev & 
	0.0006854182 & 0.0011434344 & 0.0007945224 & 0.0012479476 & 0.0008839566 &0.0007044640\\
	& MSE & 
	7.995654e-06 & 1.746024e-06 & 6.315694e-07 & 1.569822e-06 & 7.813885e-07 & 4.962843e-07
\\	\hline 
\multirow{3}{*}{$\text{1 min}$}&	Mean & 
	0.02288498 & 0.02066931 & 0.01995456 & 0.01984824 & 0.01997237 & 0.02000242\\
	& Std Dev & 
	0.0006482329 & 0.0010605652  & 0.0007468549 & 0.0011609025 & 0.0007887707 & 0.0006469303\\
	& MSE & 
	8.743311e-06 & 1.572774e-06 & 5.598574e-07 & 1.370725e-06 & 6.229225e-07 & 4.185247e-07\\
	\hline 
	\multirow{3}{*}{$\text{30 sec}$} &
	Mean & 
	0.02293765 & 0.02075251 & 0.01998865 & 0.01993685 & 0.02000009 & 0.02001709\\
	& Std Dev & 
	0.0006537998 & 0.0010611910 & 0.0007515176 & 0.0011497640 & 0.0007185258 & 0.0006364266\\
	& MSE & 
	9.057229e-06 & 1.692391e-06 & 5.649076e-07 & 1.325945e-06  & 5.162794e-07 & 4.053310e-07\\
		\hline 
\multirow{3}{*}{$\text{1 sec}$}&
	Mean & 
	0.02296041 & 0.02076158 & 0.01998938 & 0.01994110 & 0.02000240 & 0.02000628\\
	& Std Dev & 
	0.0006346972 & 0.0010546469 & 0.0007285086 & 0.0011415267 & 0.0006393828 & 0.0005973219\\
	& MSE & 
	9.166839e-06 & 1.692287e-06 & 5.308377e-07 & 1.306553e-06 & 4.088161e-07 & 3.568328e-07\\
	\hline
	\end{tabular}
}
\vspace{.2 cm}
\caption{Sample means, standard deviations, and mean-squared errors for different estimators of $\sigma=0.02$ based on $1000$ simulations.}
\label{Table2}	
\end{table}

\subsection{Estimators for $\kappa$}
We compare the finite sample performance of the following three estimators, which are respectively denoted by $\hat\kappa_{1},\hat\kappa_{2},\hat\kappa_{3}$.
\begin{enumerate}
	\item The estimator $\hat\kappa_{n,K}$ given in (\ref{FrstEstKappa}) with $\sigma$ replaced with the estimate $\hat\sigma''_{3}$ in Eq.~(\ref{BstEstKapp}) and $K$ determined by {an estimate of the optimal value $K^{*}_{3}$ given in (\ref{Exp1OK2}) obtained by replacing $\sigma$ and $\bbe\tilde{\varepsilon}^{4}$ with $\hat\sigma''_{3}$ and Eq.~(\ref{Estfouthadnsecmmnt}-i), respectively.}
%	\begin{equation}\label{Exp1OKb}
%	K^{*}=n\left(\frac{6(\bbe\varepsilon^{2})^{2}}{T^{2}\sigma^{4}}\right)^{1/3}.
%\end{equation}

\item The unbiased estimator $\hat{\bar{\kappa}}_{n}$ defined in (\ref{EstmUbsKap}) with the same value of $K$ as the previous item. As before, we replace $\sigma$ by the estimator  $\hat\sigma''_{3}$.

\item Again, the unbiased estimator $\hat{\bar{\kappa}}_{n}$ in (\ref{EstmUbsKap}) replacing $\sigma$ with  $\hat\sigma''_{3}$, but now the value of $K$ is given by (\ref{Exp1OK3}). We replace $\sigma$ therein with $\hat\sigma''_{3}$, while to estimate $e(\varepsilon)={\rm Var}\left((\varepsilon_{2}-\varepsilon_{1})^{4}\right)$, we exploit the limits in (\ref{Estfouthadnsecmmnt}).

\end{enumerate} 

The sample mean, standard deviation, and mean-squared error (MSE) based on $1000$ simulations are presented in Table \ref{Table3}. Here, we take $T=252$ days and $\sigma_{0}=0.063$. As expected, the estimator {$\hat\kappa_{3}$} has {much} better performance than any other estimator therein. 
  
\begin{center}
\begin{table}[ht]
{\scriptsize
	\begin{tabular}{|c||c|c|c||c|c|c||}
	\hline
	 	& $\hat\kappa_{1}$ & $\hat\kappa_{2}$ &  $\hat\kappa_{3}$  &   $\hat\kappa_{1}$ & $\hat\kappa_{2}$ &  $\hat\kappa_{3}$  
\\	
\hline 
 &\multicolumn{3}{c||}{
	$\delta_{n}=\text{5 min}$ } & \multicolumn{3}{c||}{
	$\delta_{n}=\text{1 min}$ }\\
	\hline
	Mean & 
	0.57771957 & 0.29982420 & 0.29967835  &  0.57428966 & 0.29189326 & 0.29686684\\%& 0.26176571\\
	Std Dev & 
	0.1783289311 & 0.1832631941 & 0.0979104650 &0.1571320926 & 0.1599275870 & 0.0758019358\\ %& 0.1793404144 \\
	MSE & 
	1.089294e-01 & 3.358543e-02 & 9.586563e-03 & 9.992531e-02 & 2.564255e-02 & 5.755750e-03\\ %& 3.362485e-02 
\hline 
& \multicolumn{3}{c||}{
	$\delta_{n}=\text{30 sec}$ } &  \multicolumn{3}{c||}{
	$\delta_{n}=\text{1 sec}$ }\\
	\hline
	Mean & 
	0.58111784 & 0.29929056 & 0.29677713 &0.57371817 & 0.29046728 & 0.29455234 \\%& 0.29531806\\
	Std Dev & 
	0.161799873 & 0.163678990 & 0.069347518 &0.162874998 & 0.165066890 & 0.066836990\\%& 0.163256840\\
	MSE & 
	1.052064e-01 & 2.679132e-02 & 4.819465e-03 & 1.014499e-01 & 2.733795e-02 & 4.496860e-03\\	\hline
	\end{tabular}
}
%\vspace{.2 cm}
\caption{Sample means, standard deviations, and mean-squared errors for different estimator of $\kappa=0.3$ based on $1000$ simulations.}
\label{Table3}	
\end{table}
\end{center}

\subsection{{Rate of Convergence Analysis}}
{In this section we study the rates of convergence of the standard errors of the unbiased estimators $\hat{\bar{\sigma}}^{2}_{n,K}$ and  $\hat{\bar{\kappa}}_{n,K}$ as defined by Eqs.~(\ref{ZMAS0b}) and  (\ref{EstmUbsKap}), when $K$ is chosen according to the optimal values (\ref{Exp2OK}) and (\ref{Exp1OK3}), respectively. In particular, we want to assess our claim that the convergence rates of the estimator's variances are faster than $T^{-1}$. To this end, we plot $\log(\widehat{{\rm Var}}\left(\hat{\bar{\sigma}}_{n,K^{*}_{2},T}\right))$  against $\log(T)$  for $T$'s ranging from 2 months to 2 years and eight intraday sampling frequencies $\delta_{n}$ (see left panel in Figure \ref{Plot1}). We also show the best linear fit for each plot. Here, $\widehat{{\rm Var}}\left(\hat{\bar{\sigma}}_{n,K^{*}_{2},T}\right)$ represents the sample variance of the estimator $\hat{\bar{\sigma}}_{n,K^{*}_{2},T}$ computed by Monte Carlo using 200 simulations. In Table \ref{Table6}, we also report the 95\% confidence intervals for the slopes of the best linear fits (second column in the table). It is apparent that the linear fit is very good, which indicates that ${{\rm Var}}\left(\hat{\bar{\sigma}}_{n,K^{*}_{2},T}\right)\propto T^{-\beta}$, for large $T$ and some $\beta<0$, and furthermore, the slope's estimates indicate that the convergence rate of  ${{\rm Var}}\left(\hat{\bar{\sigma}}_{n,K^{*}_{2},T}\right)$ is slightly better than $T^{-1}$ (the average rate is $T^{-1.03}$). We also perform the same analysis for the estimator $\hat{\sigma}_{3}''$, as described in Section \ref{Sect:EstSigma}, which is designed to be a data-drive proxy of the oracle estimator $\hat{\bar{\sigma}}_{n,K^{*}_{2},T}$. The results are show in the right panel of  Figure \ref{Plot1} and the third column of Table \ref{Table6}. The average convergence rate of ${{\rm Var}}\left(\hat{{\sigma}}_{3}''\right)$ is $T^{-1.045}$. Note that the CI's indicate that the slope is significantly different than $-1$ in almost all cases. We carry out the same analyses for the estimators for $\kappa$. The graphs of  $\log(\widehat{{\rm Var}}\left(\hat{\bar{\kappa}}_{n,K^{*}_{4},T}\right)$ and $\log(\widehat{{\rm Var}}\left(\hat{{\kappa}}_{3}\right)$ against $\log(T)$ are shown in Figure \ref{Plot2}. The CI's for the slope of the best linear fits are shown in Table \ref{Table6} (last two columns). The average convergence rate of the variance of  $\hat{\bar{\kappa}}_{n,K^{*}_{4},T}$ is $T^{-1.15}$, while the average convergence rate of the variance of $\hat{\kappa}_{3}$ is $T^{-1.18}$.}

\begin{center}
\begin{table}[ht]
{\footnotesize
	\begin{tabular}{|c|c|c|c|c|}
	\hline
	 $\delta_{n}$ & $\log\left(\widehat{{\rm Var}}\left(\hat{\bar{\sigma}}_{n,K^{*}_{2}}\right)\right)$  & $\log\left(\widehat{\rm Var}\left(\hat\sigma_{3}''\right)\right)$ &  $\log\left(\widehat{{\rm Var}}\left(\hat{\bar{\kappa}}_{n,K^{*}_{4}}\right)\right)$ & $\log\left(\widehat{\rm Var}\left(\hat\kappa_{3}\right)\right)$ \\
	\hline
	5 sec & $-1.036\pm 0.025$ & $-1.032\pm0.027$  & $-1.234\pm 0.122$ & $-1.219\pm 0.127$ \\
	\hline
	10 sec & $-1.053\pm 0.026$ & $-1.040\pm0.026$  & $-1.272\pm 0.151$ & $-1.219\pm 0.171$ \\
	\hline
	30 sec& $-1.031\pm 0.025$ & $-1.058\pm 0.026$  & $-1.22\pm 0.138$ & $-1.197\pm 0.132$ \\
	\hline
	1 min& $-1.043\pm 0.032$ & $-1.031\pm0.032$ & $-1.315\pm 0.158$ & $-1.196\pm 0.178$ \\
	\hline
	10 min & $-1.001\pm 0.026$ & $-0.998\pm 0.024$  & $-1.086\pm 0.199$ & $-1.229\pm 0.15$ \\
	\hline
	20 min & 	$-1.045\pm0.030$ & $-1.073\pm 0.026$  & $-1.056\pm 0.099$ & $-1.268\pm 0.187$ \\
	\hline
	30 min &  $-1.028\pm0.036$ & $-1.076\pm 0.019$  & $-0.931\pm 0.177$ & $-1.056\pm 0.232$ \\
	\hline
	1 hr & 	$-1.041\pm 0.020$  & $- 1.053\pm 0.023$ & $-1.124\pm 0.105$  & $-1.072\pm 0.133$ \\
	\hline
	\end{tabular}
	\vspace{.2 cm}
\caption{95\% CI's for the slope of the linear regression fit of $\log(\widehat{{\rm Var}}\left(\text{sigma Estimator}\right))$ against $\log(T)$ for $T\in \{2\text{m}, 3\text{m},\dots, 24\text{m}\}$, and $\log(\widehat{{\rm Var}}\left(\text{kappa Estimator}\right))$ against $\log(T)$ for $T\in \{12\text{m}, 13\text{m},\dots, 24\text{m}\}$.}
\label{Table6}
}
\end{table}	
\end{center}

\begin{figure}[htp]
    {\par \centering    
    \includegraphics[width=8.0cm,height=7.5cm]{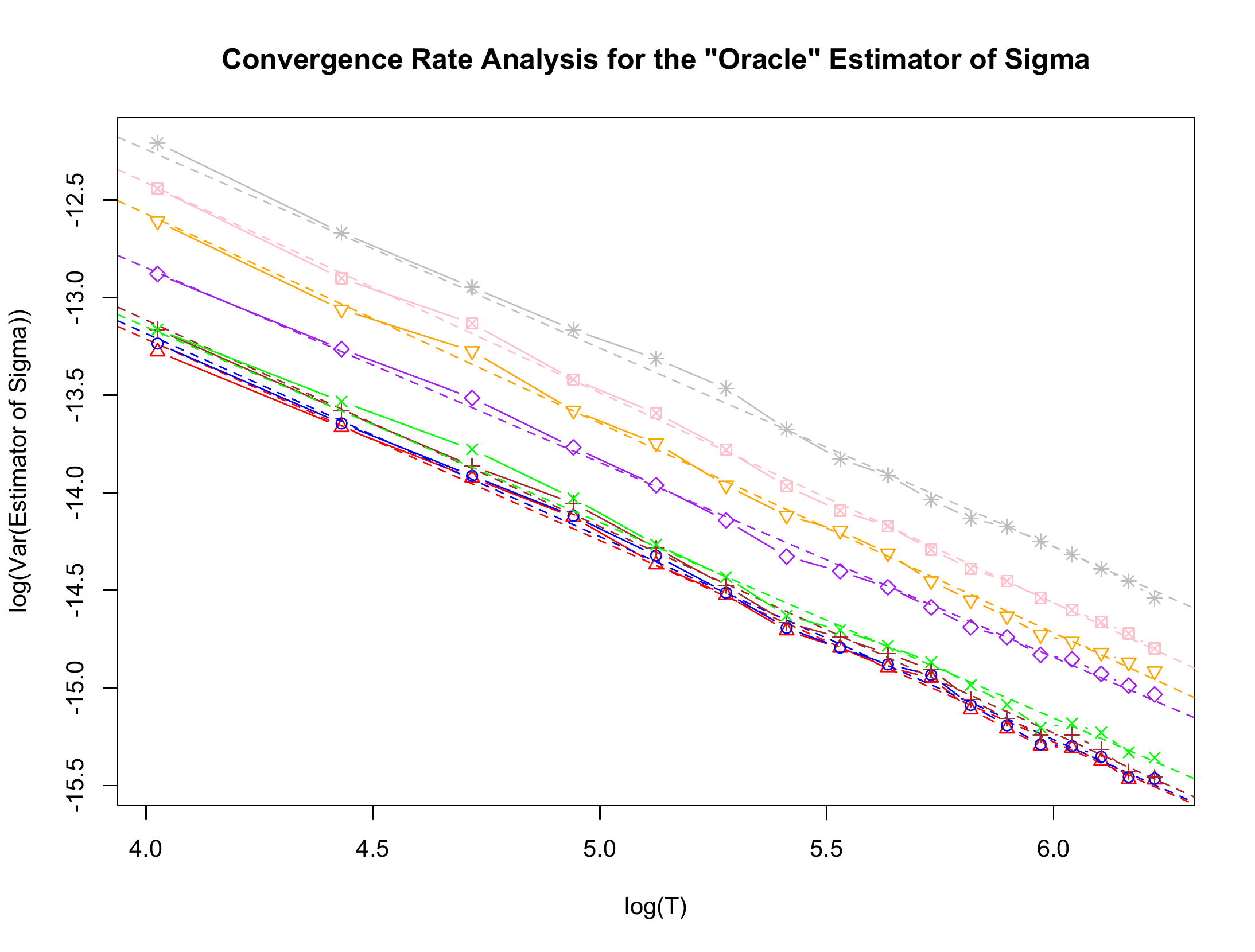}
    \includegraphics[width=8.0cm,height=7.5cm]{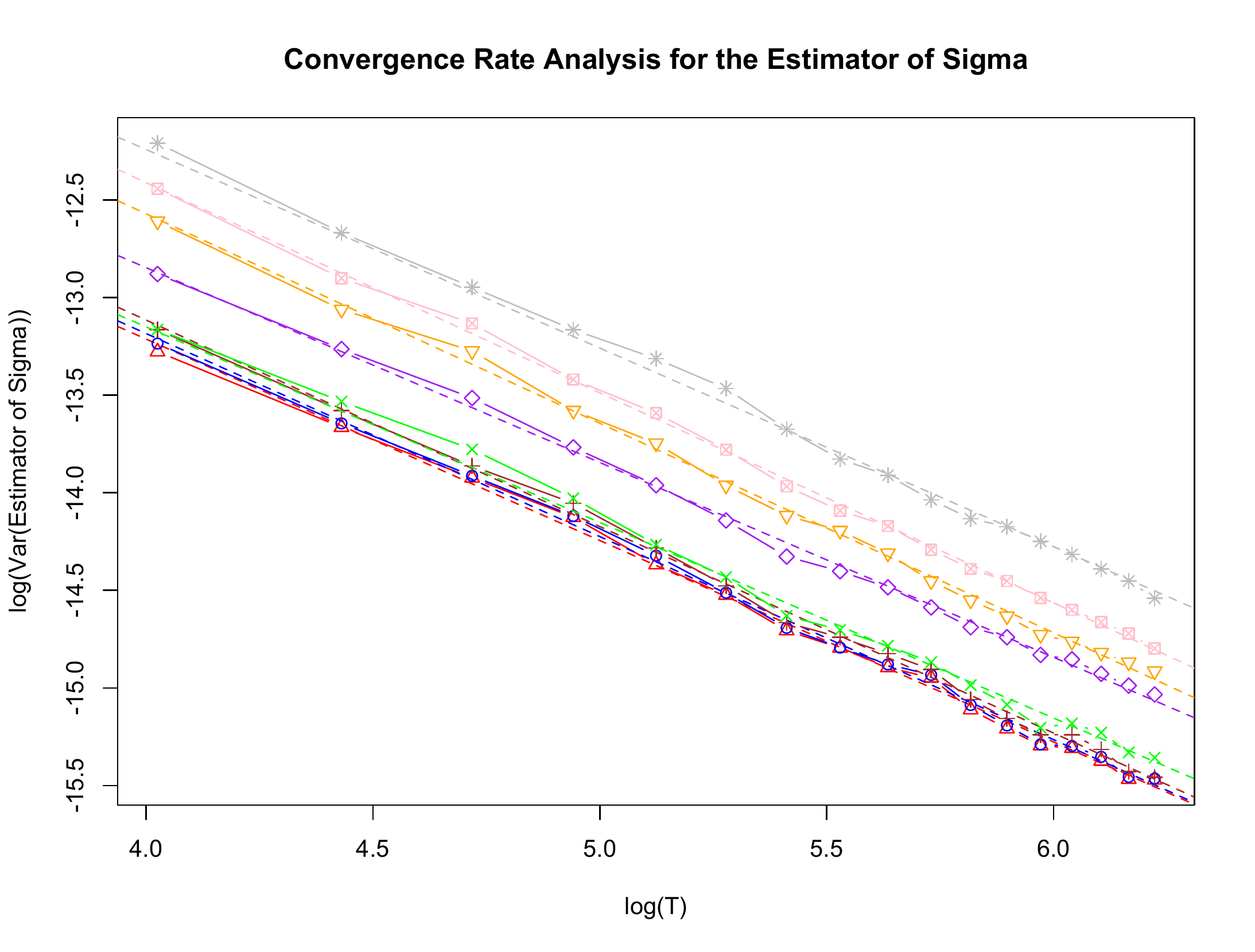}
    \par}\vspace{-.8 cm}
    \caption{Regression Analysis of $\log(\widehat{{\rm Var}}\left(\hat{\bar{\sigma}}_{n,K^{*}_{2},T}\right)$ against $\log(T)$ (left panel) and $\log\left(\widehat{{\rm Var}}\left(\hat{\sigma}_{3}''\right)\right)$ (right panel) for $T\in \{2\text{ m}, 3\text{ m},\dots, 24\text{ m}\}$, and $\delta_{n}=5 \text{ sec}$ (Red), $\delta_{n}=10  \text{ sec}$ (Blue), $\delta_{n}=30  \text{ sec}$ (Brown), $\delta_{n}=1  \text{ min}$ (Green), $\delta_{n}=10  \text{ min}$ (Purple), $\delta_{n}=20  \text{ min}$ (Orange), $\delta_{n}=30  \text{ min}$ (Pink), and $\delta_{n}=1  \text{ hr}$ (Grey). The sample variance is computed based on 200 simulations.} 
\label{Plot1}
\end{figure}

\begin{figure}[htp]
    {\par \centering    
    \includegraphics[width=8.0cm,height=7.5cm]{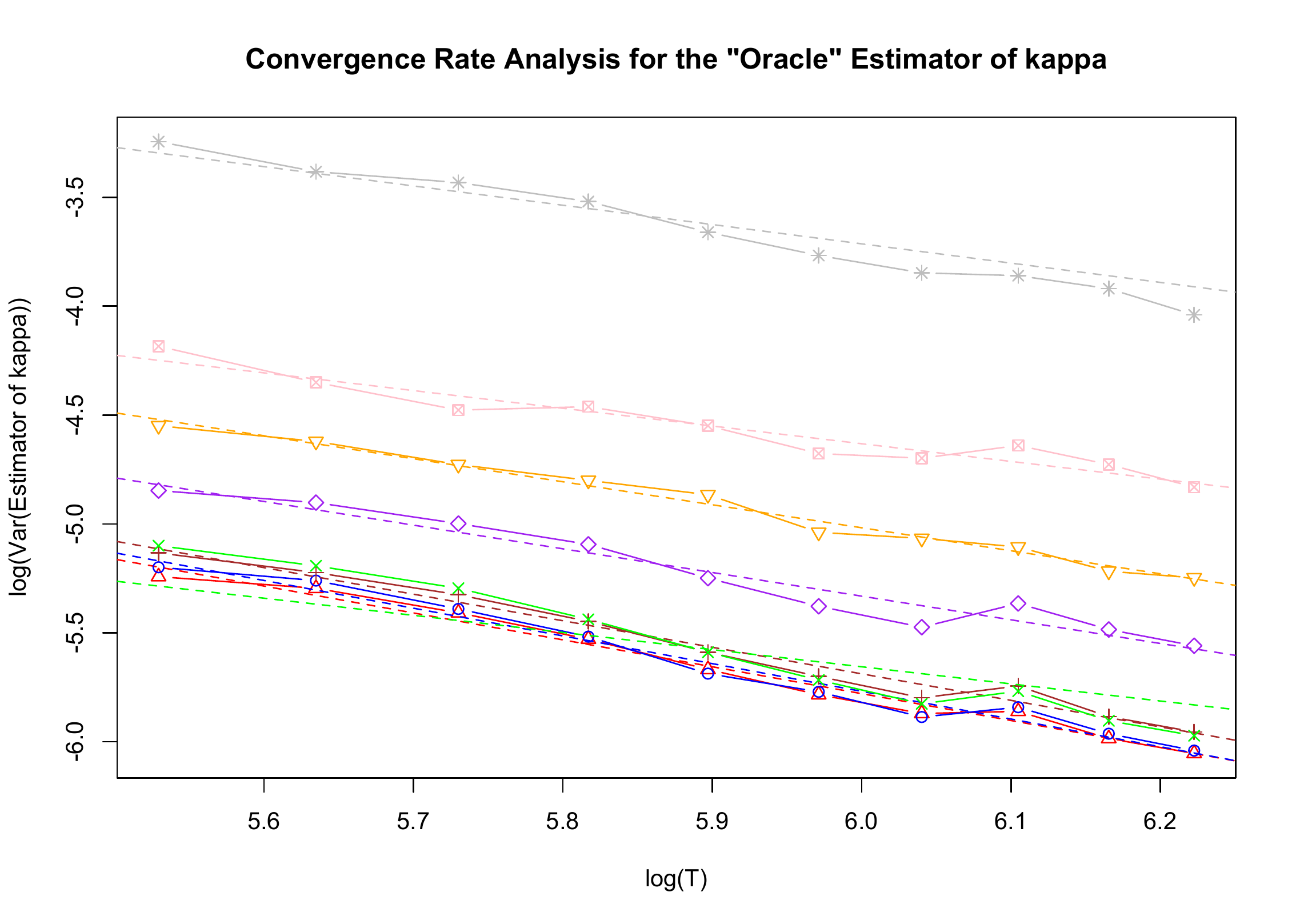}
    \includegraphics[width=8.0cm,height=7.5cm]{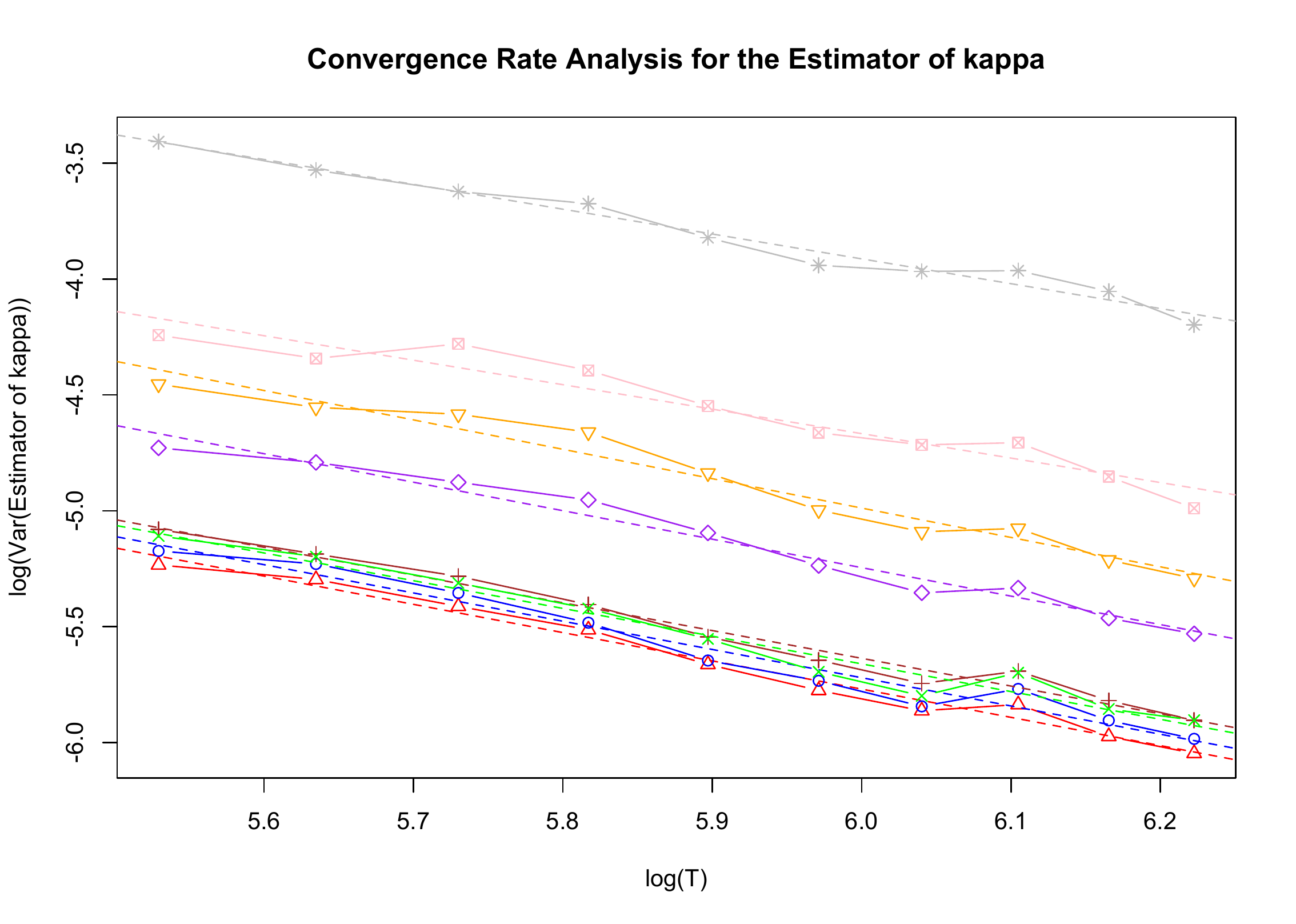}
    \par}\vspace{-.8 cm}
    \caption{Regression Analysis of $\log(\widehat{{\rm Var}}\left(\hat{\bar{\kappa}}_{n,K^{*}_{2},T}\right)$ against $\log(T)$ (left panel) and $\log\left(\widehat{{\rm Var}}\left(\hat{\kappa}_{3}\right)\right)$ (right panel) for $T\in \{2\text{ m}, 3\text{ m},\dots, 24\text{ m}\}$, and $\delta_{n}=5 \text{ sec}$ (Red), $\delta_{n}=10  \text{ sec}$ (Blue), $\delta_{n}=30  \text{ sec}$ (Brown), $\delta_{n}=1  \text{ min}$ (Green), $\delta_{n}=10  \text{ min}$ (Purple), $\delta_{n}=20  \text{ min}$ (Orange), $\delta_{n}=30  \text{ min}$ (Pink), and $\delta_{n}=1  \text{ hr}$ (Grey). The sample variance is computed based on 200 simulations.} 
\label{Plot2}
\end{figure}

\subsection{Empirical study}\label{Sec:EmStudy}
We now proceed to analyze the performance of the proposed estimators when applied to real data. As it was explained above and was theoretically verified by Propositions \ref{PrpABSNO}-\ref{PrpABKNO}, traditional estimators are not stable as the sampling frequency increases. Indeed, $\hat\sigma_{n}$ and $\tilde\sigma_{n}$ both diverge to $\infty$ while $\hat\kappa_{n}$ and $\tilde{\kappa}_{n}$  converge to $0$, as $n\to\infty$. The {objective} is to verify that {the proposed} estimators do not exhibit {the aforementioned} behaviors at very high-frequencies. 

We consider high-frequency stock data for several stocks during 2005, which were obtained from the NYSE TAQ database of Wharton's WRDS system. For briefness and illustration purposes, we only show Intel (INTC)  and Pfeizer (PFE). For these, we compute the estimator $\hat{\varrho}$ defined in (\ref{estVarRho}), the estimator $\hat\sigma_{n,K}$ defined in (\ref{FrstEstSigma}) with $K=1$, the estimator $\hat{\bar{\sigma}}_{n,K}$ defined in (\ref{ZMAS0b}) with $K=\hat{\hat{K}}^{*}_{1}$ as given in (\ref{BstEstKapp}), the estimator $\hat\kappa_{n,K}$ defined in (\ref{FrstEstKappa}) with $K=1$, and finally the estimator $\hat{\bar{\kappa}}_{n,K}$ defined in (\ref{EstmUbsKap}) with $K=\hat{K}_{4}^{*}$ as given in (\ref{Exp1OK3}). In the case of $\hat\kappa_{n,1}$, we used $\sigma=\hat\sigma_{n,1}$. Both $\hat\sigma_{n,1}$ and $\hat\kappa_{n,1}$ represent the estimators without any technique to alleviate the effect of the microstructure noise. As one can see in Tables \ref{Table4}-\ref{Table5}, the estimators $\hat{\bar{\sigma}}$ and $\hat{\bar{\kappa}}$ do not exhibit the drawbacks of the estimators $\hat\sigma$ and $\hat\kappa$ at high frequencies. As a {conclusion} of the {empirical} results therein, {we deduce that} Intel's stock exhibits an annualized volatility $\sigma$ of about $0.014*\sqrt{252}=0.22$ per year, while its excess kurtosis increases with $1/\delta$ at a rate of {about $0.5$} {(see item 2 above Eq.~(\ref{SD}) for the interpretation of $\kappa$)}. {By} comparison, even though the volatility of Pfizer's stock is just slightly larger (about $0.015*\sqrt{252}=0.23$), its excess kurtosis increases at a rate of about $2.3$ with $1/\delta$, showing much more riskiness due to the much heavier tails of its return's distribution.  This example illustrates the importance of considering a parameter which measures the tail heaviness of the return distribution and not only its variance.

\begin{center}
\begin{table}[ht]
{\footnotesize
	\begin{tabular}{|c|c|c|c|c|c|}
	\hline
	 & $\hat{\varrho}$ & $\hat\sigma_{n,1}$ &  $\hat{\bar{\sigma}}_{n,\hat{\hat{K}}^{*}_{1}}$ & $\hat\kappa_{n,1}$ & $\hat{\bar{\kappa}}_{n,\hat{K}^{*}_{4}}$ \\
	\hline
	20 min & 0.002198811 & 0.013732969 & 0.013115165 & 0.772846688 & 0.645084939\\
	\hline
	10 min & 0.001584536 & 0.013995671 & 0.013112833 & 0.589344904 & 0.727208959\\
	\hline
	5 min & 0.001152404 & 0.014394983 & 0.013253727 & 0.495378704 & 0.768302688\\
	\hline
	1 min & 0.0005581856 & 0.0155908617 & 0.0136519981 & 0.3499494734 & 0.7293149570\\
	\hline
	30 sec & 0.0004113675 & 0.0162494093 & 0.0139405766 & 0.2817929514 & 0.6875741045\\
	\hline
	20 sec & 	0.0003483541 & 0.0168528945 & 0.0141596310 & 0.2566280373 & 0.6575495762 \\
	\hline
	10 sec &  0.0002712869 & 0.0185608431 & 0.0145174963 & 0.1831341414 & 0.5921934015 \\
	\hline
	5 sec & 	0.0002174315 & 0.0210381061 & 0.0147818871 & 0.1084570206 & 0.4987667343 \\
	\hline
	\end{tabular}
	\vspace{.2 cm}
\caption{Estimation of the parameters $\sigma$ and $\kappa$ of a subordinated Brownian motion with microstructure noise for INTC (Intel) stock.}
\label{Table4}
}
\end{table}	
\end{center}

\begin{center}
\begin{table}[ht]
{\footnotesize
	\begin{tabular}{|c|c|c|c|c|c|}
	\hline
	 & $\hat{\varrho}$ & $\hat\sigma_{n,1}$ &  $\hat{\bar{\sigma}}_{n,\hat{\hat{K}}^{*}_{1}}$ & $\hat\kappa_{n,1}$ & $\hat{\bar{\kappa}}_{n,\hat{K}^{*}_{4}}$ \\
	\hline
	20 min & 0.002310884 & 0.014432934 & 0.014279133 & 3.552809339 & 3.665645436\\
	\hline
	10 min & 0.001678615 & 0.014826633 & 0.013921679 & 3.330420039 & 4.192632331\\
	\hline
	5 min & 0.001223294 & 0.015280492 & 0.013758805 & 3.395593192 & 4.458814370
	\\
	\hline
	1 min & 0.000581559 & 0.016243711 & 0.014289601 & 2.885849749 & 3.074717720\\
	\hline
	30 sec & 0.0004379718 & 0.0173003060 & 0.0147847384 & 2.1009477905 & 2.5399891978\\
	\hline
	20 sec & 	0.0003733763 & 0.0180634325 & 0.0149589310 & 1.8189209947 & 2.3582752416 \\
	\hline
	10 sec &  0.0003021168 & 0.0206701623 & 0.0150440707 & 1.0395706194 & 2.3194219287 \\
	\hline
	5 sec & 	0.0002547010 & 0.0246442060 & 0.0151395852 & 0.5255478783 & 2.3750789809 \\
	\hline
	\end{tabular}
	\vspace{.2 cm}
\caption{Estimation of the parameters $\sigma$ and $\kappa$ of a subordinated Brownian motion with microstructure noise for PFE (Pfeizer) stock.}
\label{Table5}
}
\end{table}	
\end{center}
\newpage
 
\appendix

\section{Proofs}\label{Prf}

\subsection{Proof of Proposition \ref{BEK}.}
We shall need the following {standard result that can easily be shown using} the moment generating function for Poisson integrals (see, e.g., \cite[Chapter 2]{Cont:2003})}:
\begin{lmma}\label{MPI}
	Suppose that $M$ is a Poisson random measure {on an open domain of $\bbr^{d}$ with mean measure $m$} and let $\bar{M}(f)=\int f(z) (M-m)(dz)$ denote the integral of $f$ with respect the compensated random measure $\bar{M}=M-m$. 
	If $m(|f|^{k}):=\int |f(z)|^{k} m(dz)<\infty$, for $k=1,\dots,5$, then $\bbe \left(\bar{M}(f)^{k}\right)=m(f^{k})$, for $k=2,3$, $\bbe \left(\bar{M}(f)^{4}\right)=3m(f^{2})^{2}+m(f^{4})$, and $\bbe \left(\bar{M}(f)^{5}\right)=10m(f^{2})m(f^{3})+m(f^{5})$. Similarly, $\bbe\left(\bar{M}(g)\bar{M}(f)^{k}\right)=m(g f^{k})$ and $\bbe\left(\bar{M}(g)\bar{M}(f)^{3}\right)=m(g f^{3})+3m(f^{2}) m(gf)$. 
\end{lmma}	
\begin{lmma}\label{MPI2}
	Let $M$ be {the jump} measure of a L\'evy process $X$ with L\'evy measure $\nu$ {(i.e., $M((s,t)\times B):=\#\{u\in(s,t): \Delta X_{u}\in B\}$, for any $s<t$ and $B\in\mathcal{B}(\mathbb{R}^{d})$),} and let $\bar{M}(dt,dx):=M(dt,dx)-dt\nu(dx)$ be the corresponding compensated measure. Also, suppose that $f$ is such that $\int |f(x)|^{k}\nu(dx)<\infty$ for some $k\geq{}2$. Then, {there exists a constant $A_{k}(f)$ such that,} for any $T\geq{}1$, 
	\[
		\bbe \left|\frac{1}{T}\int_{0}^{T}\int f(x)\bar{M}(dt,dx)\right|^{k} \leq A_{k}(f) T^{-k/2}.
	\]
\end{lmma}	
\noindent\textbf{Proof.} 
	{Throughout the proof, let $\bar{M}_{s,t}(f):=\int_{s}^{t}\int f(x)\bar{M}(dt,dx)$ and let $[T]$ be the integer part of $T$. We need the following classical inequality} (see \cite[Lemma 5.3.1]{Bickel}):
	\begin{equation}\label{Bickel}
		\bbe |\bar{Z}_{n}-\mu_{Z}|^{k}\leq C_{k} \bbe |Z_{1}|^{k} n^{-k/2},
	\end{equation}
	where  $\bar{Z}_{n}=\frac{1}{n}\sum_{i=1}^{n}Z_{i}$, $\mu_{Z}=\bbe Z_{1}$, and $\{Z_{i}\}_{i}$ are i.i.d. such that $\bbe |Z_{1}|^{k}<\infty$. First, {note that
%	\[
%		\frac{1}{T}\int_{0}^{T}\int f(x)\bar{M}(dt,dx)=\frac{1}{T}\int_{0}^{[T]}\int f(x)\bar{M}(dt,dx)+\frac{1}{T}\int_{[T]}^{T}\int f(x)\bar{M}(dt,dx).
%	\]
%	Hence, 
	\begin{align*}
		\bbe \left|\frac{1}{T}\bar{M}_{0,T}(f)\right|^{k}
		\leq 2^{k}\bbe \left|\frac{[T]}{T}\frac{1}{[T]}\bar{M}_{0,[T]}\left(f\right)\right|+2^{k} \frac{1}{T^{k}}\bbe \left|\bar{M}_{[T],T}(f)\right|^{k}.
	\end{align*}
	For} the first term on the right-hand side above, we apply (\ref{Bickel}) {with $Z_{i}:=\bar{M}_{i-1,i}(f)$, which are i.i.d. because $M$ is a Poisson random measure.} For the second term, {we apply Burkholder-Davis-Gundy inequality (see \cite{Protter}) to get, 
	\begin{align*}
		\bbe \left|\int_{[T]}^{T}\int f(x)\bar{M}(dt,dx)\right|^{k}
%		&\leq 
%		\bbe \left|\sup_{t\leq{}1}\int_{0}^{t}\int f(x)\bar{M}(dt,dx)\right|^{k}
		\leq B_{k}^{k} \bbe\left| \int_{0}^{1}\int f^{2}(x)M(dt,dx)\right|^{k/2}.
	\end{align*}
	This completes the proof.}
\hfill$\Box$

\noindent
\noindent\textbf{Proof of Proposition \ref{BEK}.}
{Throughout the proof, $M$ denotes the jump measure of the L\'evy process $X$; i.e., $M((s,t)\times B):=\#\{u\in(s,t): \Delta X_{u}\in B\}$, for any $s<t$ and $B\in\mathcal{B}(\mathbb{R})$. In particular, let us note that $M$ is Poisson random measure with mean measure $dt\nu(dx)$ and $\sum_{t\leq{}T}(\Delta X_{t})^{\ell}=\int_{0}^{t}\int x^{\ell}M(dt,dx)$. Let also $\bar{M}(dt,dx):=M(dt,dx)-dt\nu(dx)$ be the corresponding compensated measure. Let us start by noting} the identity
\begin{equation}\label{EI1}
	\frac{1}{(1+x)^{2}}=\sum_{i=0}^{k-1}(-1)^{i} (i+1) x^{i} 
	+\frac{(-1)^{k}x^{k}}{(1+x)^{2}}(k+1+kx),
\end{equation}
and the notation 
\[
	\hat{\mu}_{k}^{(T)}:= \frac{1}{T}\int_{0}^{T}\int x^{k} M(dt,dx),
	\quad \hat{D}_{T}:=\frac{\hat\mu_{2}^{(T)}}{c_{2}(X_{1})}-1.
\]
{In particular, $\hat\kappa^{(T)}=(1/3)\hat\mu_{4}^{(T)}/(\hat\mu_{2}^{(T)})^{2}$. Then,} we have the following decomposition: 
\begin{align*}
	\bbe \hat\kappa^{(T)}&=\frac{1}{3c_{2}^{2}(X_{1})}
	\bbe \left\{\frac{\hat\mu_{4}^{(T)}}{(1+\hat{D}_{T})^{2}}\right\}\\&=\frac{1}{3c_{2}^{2}(X_{1})}
	\bbe \left\{\hat\mu_{4}^{(T)}\left(1-2 \hat{D}_{T}+3\hat{D}_{T}^{2}
	-4\hat{D}_{T}^{3}+5\hat{D}_{T}^{4}-6\hat{D}_{T}^{5}\right)\right\}\\
	&\quad+\frac{1}{3}
	\bbe\left\{\hat\mu_{4}^{(T)}\left(\hat\mu_{2}^{(T)}\right)^{-2}\left(7+6 \hat{D}_{T}\right)\hat{D}_{T}^{6}\right\}\\
	&=:L_{T}+R_{T}.
\end{align*}
Let us first analyze the residual term $R_{T}$ using the following easy consequence of the triangle inequality:
\begin{equation}\label{UI2}
	(\hat\mu_{4}^{(T)})^{1/2}
	=\frac{1}{T^{1/2}}\left(\sum_{s\leq{}T}
	\left(\Delta X_{s}\right)^{4}\right)^{1/2}\leq 
	\frac{1}{T^{1/2}}
	\sum_{s\leq{}T}\left(\Delta X_{s}\right)^{2}
	=T^{1/2}\mu_{2}^{(T)}.
\end{equation}
Thus, since $7+6\hat{D}_{T}=1+6(1+\hat{D}_{T})=1+6\hat\mu_{2}^{(T)}/c_{2}(X_{1})>0$, we have that  
\begin{align*}
	0\leq R_{T}
	&\leq  \frac{7T}{3}
	\bbe\left(\hat{D}_{T}^{6}\right)
	+\frac{6 T}{3}\bbe
	\left(\hat{D}_{T}^{7}\right)\\
	&=  \frac{7T}{3c_{2}^{6}(X_{1})}
	\bbe\left(\hat\mu_{2}^{(T)}-c_{2}(X_{1})\right)^{6}
	+\frac{6 T}{3c_{2}^{7}(X_{1})}\bbe
	\left(\hat\mu_{2}^{(T)}-c_{2}(X_{1})\right)^{7}.
\end{align*}
Using that $\bbe \hat\mu_{2}^{(T)}=c_{2}(X_{1})$ and Lemma \ref{MPI2}, 
\(
	 R_{T}= O(T^{-2}).
\)
Similarly, using Lemma \ref{MPI}, the first four terms of $L_{T}$ (i.e. those multiplying $\hat{D}^{i}_{T}$ up to $i=3$) {are given by}
\[
	\frac{c_{4}(X_{1})}{3 c_{2}^{2}(X_{1})}-\frac{2c_{6}(X_{1})}{3 c_{2}^{3}(X_{1})}T^{-1}
	+\frac{c_{4}^{2}(X_{1})}{c_{2}^{4}(X_{1})}T^{-1}+O(T^{-2}).
\]
The last two term of $L_{T}$ can be seen to be $O(T^{-2})$ from Lemma \ref{MPI2} and Cauchy inequality. Indeed,
\begin{align*}
	\left|\bbe\hat\mu_{4}^{(T)}\hat{D}_{T}^{4}\right|&\leq 
	c_{4}(X_{1})\left|\bbe \hat{D}_{T}^{4}\right|+\frac{1}{c_{2}^{4}(X_{1})}
	\left|\bbe \left(\hat\mu_{4}^{(T)}-c_{4}(X_{1})\right)\left(\hat\mu_{2}^{(T)}-c_{2}(X_{1})\right)^{4}\right|\\
	&\leq K c_{2} T^{-2}+\left(\bbe \left(\hat\mu_{4}^{(T)}-c_{4}(X_{1})\right)^{2}
	\bbe\left(\hat\mu_{2}^{(T)}-c_{2}(X_{1})\right)^{8}\right)^{1/2},
\end{align*}
which is $O(T^{-2})$ in light of Lemma  \ref{MPI2}.
We finally obtain that 
\[
	\bbe \hat{\kappa}_{n}{\;{\stackrel{n\to\infty}{\longrightarrow}}\;} \bbe \hat\kappa^{(T)}=\frac{c_{4}(X_{1})}{3 c_{2}^{2}(X_{1})}
	-\frac{2c_{6}(X_{1})}{3 c_{2}^{3}(X_{1})}T^{-1}
	+\frac{c_{4}^{2}(X_{1})}{c_{2}^{4}(X_{1})}T^{-1}+O(T^{-2}).
\]
%One can similarly verify that 
%\[
%	{\rm Var}\left(\hat\kappa^{(T)}\right)=O(T^{-1}).
%\]
In order to show the bound for the variance, we use again (\ref{EI1}) to get 
\begin{align*}
	 \hat\kappa^{(T)}&=\frac{\hat\mu_{4}^{(T)}}{3c_{2}^{2}(X_{1})}
	\left(1-2 \hat{D}_{T}+3\hat{D}_{T}^{2}-4\hat{D}_{T}^{3}\right)+\frac{1}{3}\frac{\hat\mu_{4}^{(T)}}{\left(\hat\mu_{2}^{(T)}\right)^{2}}\left(5+4 \hat{D}_{T}\right)\hat{D}_{T}^{4}.
\end{align*}
%\begin{align*}
%	 \hat\kappa^{(T)}&=\frac{\hat\mu_{4}^{(T)}}{3c_{2}^{2}(X_{1})}
%	\left(1-2 \hat{D}_{T}+3\hat{D}_{T}^{2}-4\hat{D}_{T}^{3}\right)+\frac{1}{3}\frac{\hat\mu_{4}^{(T)}}{\left(\hat\mu_{2}^{(T)}\right)^{2}}\left(5+4 \hat{D}_{T}\right)\hat{D}_{T}^{4}.
%\end{align*}
Then, 
\begin{align*}
	 \hat\kappa^{(T)}-\frac{c_{4}(X_{1})}{3c_{2}^{2}(X_{1})}&=\frac{1}{3c_{2}^{2}(X_{1})}\left(\hat\mu_{4}^{(T)}-c_{4}(X_{1})\right)
	-\frac{2\hat\mu_{4}^{(T)}}{3c_{2}^{2}(X_{1})} \hat{D}_{T}+\frac{\hat\mu_{4}^{(T)}}{c_{2}^{2}(X_{1})}\hat{D}_{T}^{2}\\
	&\quad-\frac{4\hat\mu_{4}^{(T)}}{3c_{2}^{2}(X_{1})}\hat{D}_{T}^{3}+\frac{1}{3}\frac{\hat\mu_{4}^{(T)}}{\left(\hat\mu_{2}^{(T)}\right)^{2}}\left(5+4 \hat{D}_{T}\right)\hat{D}_{T}^{4}.
\end{align*}
After expanding the squares, taking expectations both sides, and using Cauchy's inequality together with Lemmas \ref{MPI} and \ref{MPI2}, one can check that all the terms are at least $O(T^{-2})$ except possibly the following terms:
\begin{align*}
	 &\frac{1}{9c_{2}^{4}(X_{1})}\bbe\left\{ \left(\hat\mu_{4}^{(T)}-c_{4}(X_{1})\right)^{2}\right\}-
	  \frac{4}{9c_{2}^{4}(X_{1})}\bbe\left\{ \left(\hat\mu_{4}^{(T)}-c_{4}(X_{1})\right)\hat\mu_{4}^{(T)}\hat D_{T}\right\}
	  \\
	   &+\frac{4}{9c_{2}^{4}(X_{1})}\bbe\left\{(\hat\mu_{4}^{(T)})^{2}\hat{D}_{T}^{2}\right\}.
\end{align*}
Subtracting $c_{4}(X_{1})$ from $\hat\mu_{4}^{(T)}$ in the second and third terms above, and using again  Lemmas \ref{MPI} and \ref{MPI2}, we can check that the above expression {indeed coincides with the expression in (\ref{MSEEstKapCnt}).}
\hfill$\Box$

\subsection{Proofs of Section \ref{MMECorrected}.}\label{PrfSec5}
\noindent\textbf{Proof of {Theorem \ref{VarSigmEst}.}} %[\textbf{Proof of {Theorem \ref{VarSigmEst}}}.]
Throughout we write $T_{i}$ for $T_{i,K}$. 
Clearly, 
\begin{equation}\label{VarExpFrm}
	{\rm Var}\left(\hat\sigma^{2}_{n,K}\right)=\frac{2}{K^{2}}\sum_{1\leq i<j\leq K}\frac{1}{{T_{i}T_{j}}}{\rm Cov}\left([\widetilde{X},\widetilde{X}]_{2}^{\mathcal{G}_{n}^{(i)}},[\widetilde{X},\widetilde{X}]_{2}^{\mathcal{G}_{n}^{(j)}}\right)+\frac{1}{K^{2}}\sum_{i=1}^{K}\frac{1}{{T_{i}^{2}}}{\rm Var}\left([\widetilde{X},\widetilde{X}]_{2}^{\mathcal{G}_{n}^{(i)}}\right).
\end{equation}
Each covariance in the first term on the right hand side above is given by 
\begin{align*}
	A_{i,j}&:={\rm Cov}\left([\widetilde{X},\widetilde{X}]_{2}^{\mathcal{G}_{n}^{(i)}},[\widetilde{X},\widetilde{X}]_{2}^{\mathcal{G}_{n}^{(j)}}\right)\\
%	&={\rm Cov}\left(\sum_{q=0}^{n_{i}-1}\left|\widetilde{X}(t_{i-1+(q+1)K})-\widetilde{X}(t_{i-1+qK})\right|^{2},\sum_{r=0}^{n_{j}-1}\left|\widetilde{X}(t_{j-1+(r+1)K})-\widetilde{X}(t_{j-1+rK})\right|^{2}\right)\\
	&=\sum_{q=0}^{n_{i}-1}\sum_{r=0}^{n_{j}-1}{\rm Cov}\left(\left|\widetilde{X}(t_{i-1+(q+1)K})-\widetilde{X}(t_{i-1+qK})\right|^{2},\left|\widetilde{X}(t_{j-1+(r+1)K})-\widetilde{X}(t_{j-1+rK})\right|^{2}\right)\\
%	&=n_{i}{\rm Cov}\left(\left|\widetilde{X}(t_{i-1+K})-\widetilde{X}(t_{i-1})\right|^{2},\left|\widetilde{X}(t_{j-1+K})-\widetilde{X}(t_{j-1})\right|^{2}\right)\\
%	&\quad+(n_{j}-1)
%	{\rm Cov}\left(\left|\widetilde{X}(t_{i-1+2K})-\widetilde{X}(t_{i-1+K})\right|^{2},\left|\widetilde{X}(t_{j-1+K})-\widetilde{X}(t_{j-1})\right|^{2}\right)\\
%	&=(n_{i}-1/2\pm 1/2)C\left(t_{i-1+K}-t_{j-1}\right)+(n_{j}-1/2\mp 1/2)
%	C\left(t_{j-1+K}-t_{i-1+K}\right)\\
	&=n_{i}C\left((K+i-j)\delta_{n}\right)+\left(n_{j}-1\right)
	C\left((j-i)\delta_{n}\right),
\end{align*}
where, for any $u<t<t+\delta<v$, 
\(
	C(\delta):={\rm Cov}\left(\left|\widetilde{X}(t+\delta)-\widetilde{X}(u)\right|^{2},\left|\widetilde{X}(v)-\widetilde{X}(t)\right|^{2}\right),
\)
which can be proved to depend only on $\delta>0$. More specifically, note that $C(\delta)={\rm Cov}\left(\left|S+U\right|^{2},\left|S+V\right|^{2}\right)$, 
%\begin{align*}
%	C(\delta)
%%	&={\rm Cov}\left(\left|\left(X(t+\delta)-{X}(t)\right)+\left({X}(t)-{X}(u)+\varepsilon_{t+\delta}-\varepsilon_{u}\right)\right|^{2},\right.\\
%%	&\qquad\qquad\left.\left|\left({X}(v)-{X}(t+\delta)+\varepsilon_{v}-\varepsilon_{t}\right)+\left({X}(t+\delta)-{X}(t)\right)\right|^{2}\right)\\
%		&={\rm Cov}\left(\left|S+U\right|^{2},\left|S+V\right|^{2}\right)
%\end{align*}
where $S:=X(t+\delta)-{X}(t)$, $U:={X}(t)-{X}(u)+\varepsilon_{t+\delta}-\varepsilon_{u}$, and $V:={X}(v)-{X}(t+\delta)+\varepsilon_{v}-\varepsilon_{t}$. Next, using that independence of $S$, $U$, and $V$, 
\begin{align*}
	C(\delta)&={\rm Var}\left(S^{2}\right)+2{\rm Cov}\left(S^{2},SV\right)+2{\rm Cov}\left(SU,S^{2}\right)+4{\rm Cov}\left(SU,SV\right)\\
	&={\rm Var}\left(S^{2}\right)+2\bbe(V){\rm Cov}\left(S^{2},S\right)+2\bbe(U){\rm Cov}\left(S,S^{2}\right)+4\bbe(U)\bbe(V){\rm Var}\left(S\right).
\end{align*}
Finally, using that $\bbe U=\bbe V=0$ as well as the moment formulas in (\ref{MEVGM}),
%\[
%	C(\delta)%{\rm Cov}\left(\left|\widetilde{X}(t+\delta)-\widetilde{X}(u)\right|^{2},\left|\widetilde{X}(v)-\widetilde{X}(t)\right|^{2}\right)
%	={\rm Cov}\left(\left|\left(\widetilde{X}(t+\delta)-\widetilde{X}(t)\right)+\left(\widetilde{X}(t)-\widetilde{X}(u)\right)\right|^{2},\left|\left(\widetilde{X}(v)-\widetilde{X}(t+\delta)\right)+\left(\widetilde{X}(t+\delta)-\widetilde{X}(t)\right)\right|^{2}\right)
%%	{\rm Cov}((\widetilde{X}(t+\delta-Y_s)^4,(Y_v-Y_u)^4)=Cov((Y_t-Y_u+Y_u-Y_s)^4,(Y_v-Y_t+Y_t-Y_s)^4),
%\]
%which is decomposed into 9 terms after expanding the squares and the covariance. Upon} a direct computation of each term, it follows that
$C(\delta)={\rm Var}\left(S^{2}\right)$ is given by $C(\delta)=2\sigma^{4}\delta^{2}+3\sigma^{4}\kappa \delta$. %+12\bbe\left(\varepsilon^{2}\right)^{2}+4\bbe\left(\varepsilon^{4}\right).
%\)
Using the previous formula together with the fact that {$|\frac{n(n_{j}-1)}{Kn_{i}n_{j}}-1|\leq{}U\frac{K}{n}$ and $|\frac{n}{Kn_{j}}-1|\leq{}U\frac{K}{n}$ for some constant $U$ (independent of $n$, $K$, $i$, $T$, and $j$)}, the first term in (\ref{VarExpFrm}), which we denote $A$, can be computed as follows:
\begin{align*}
	A
%	&=\frac{2}{K^{2}T^{2}}\sum_{1\leq i<j\leq K}{\rm Cov}\left([\widetilde{X},\widetilde{X}]_{2}^{\mathcal{G}_{n}^{(i)}},[\widetilde{X},\widetilde{X}]_{2}^{\mathcal{G}_{n}^{(j)}}\right)\\
	&={\frac{2n}{K^{3}T^{2}}}\sum_{1\leq i<j\leq K}
	\left(2\sigma^{4}(j-i)^{2}\delta_{n}^{2}+3\sigma^{4}\kappa (j-i)\delta_{n}
	%+12\bbe\left(\varepsilon^{2}\right)^{2}+4\bbe\left(\varepsilon^{4}\right)
	\right)\\
	&\quad+{\frac{2n}{K^{3}T^{2}}}\sum_{1\leq i<j\leq K}
	\left(2\sigma^{4}(K+i-j)^{2}\delta_{n}^{2}+3\sigma^{4}\kappa (K+i-j)\delta_{n}
	%+12\bbe\left(\varepsilon^{2}\right)^{2}+4\bbe\left(\varepsilon^{4}\right)
	\right)+\mathcal{R}\\
%	&=\frac{n}{K}\frac{4}{K^{2}T^{2}}\frac{K(K-1)}{2}
%	\left(12\bbe\left(\varepsilon^{2}\right)^{2}+4\bbe\left(\varepsilon^{4}\right)\right)\\
%	&\quad+\frac{n}{K}\frac{2}{K^{2}T^{2}}\frac{K(K-1)}{2}
%	\left(3\sigma^{4}\kappa K\right)\delta_{n}\\	
%	&\quad+\frac{n}{K}\frac{2}{K^{2}T^{2}}\delta_{n}^{2}(2\sigma^{4})\sum_{1\leq i<j\leq K}
%	\left((j-i)^{2}+(K+i-j)^{2}\right)+\mathcal{R}\\
%		&=\frac{n}{K}\frac{4}{K^{2}T^{2}}\frac{K(K-1)}{2}
%	\left(12\bbe\left(\varepsilon^{2}\right)^{2}+4\bbe\left(\varepsilon^{4}\right)\right)\\
%	&\quad+\frac{n}{K}\frac{2}{K^{2}T^{2}}\frac{K(K-1)}{2}
%	\left(3\sigma^{4}\kappa K\right)\delta_{n}\\	
%	&\quad+\frac{n}{K}\frac{2}{K^{2}T^{2}}\delta_{n}^{2}(2\sigma^{4})\frac{K^{2}(K-1)(2K-1)}{6}\\
	&=\frac{n}{K}\frac{K-1}{KT^{2}}\left(
%	24\bbe\left(\varepsilon^{2}\right)^{2}+8\bbe\left(\varepsilon^{4}\right)+
	{3\sigma^{4}\kappa K\delta_{n}+\frac{2}{3}\sigma^{4}K(2K-1)\delta_{n}^{2}}\right)+\mathcal{R}_{1},
\end{align*}
where $\mathcal{R}$ is such that 
%The last term above comes from the following summation:
%\begin{align*}
%	&\sum_{i=1}^{K-1}\sum_{\ell=1}^{K-i}\left(\ell^{2}+(K-\ell)^{2}\right)=\sum_{\ell=1}^{K-1}(K-\ell)\left(\ell^{2}+(K-\ell)^{2}\right)=K\sum_{i=1}^{K-1}\ell^{2}=\frac{K^{2}(K-1)(2K-1)}{6}.
%\end{align*}
%where 
\begin{align}\label{MaxR}
	|\mathcal{R}_{1}|&\leq \frac{{2U}(K-1)}{KT^{2}}\left(
%	24\bbe\left(\varepsilon^{2}\right)^{2}+8\bbe\left(\varepsilon^{4}\right)+
	{3\sigma^{4}\kappa K\delta_{n}+\frac{2}{3}\sigma^{4}K(2K-1)\delta_{n}^{2}}\right).
	\end{align}
Now, we consider the second term in (\ref{VarExpFrm}), which we denote $B$. Each variance term of $B$ can be written as 
%\begin{align*}
%	B=\frac{1}{K^{2}T^{2}}\sum_{i=1}^{K}{\rm Var}\left([\widetilde{X},\widetilde{X}]_{2}^{\mathcal{G}_{n}^{(i)}}\right)
%\end{align*}
\begin{align*}
	B_{i}&:={\rm Var}\left([\widetilde{X},\widetilde{X}]_{2}^{\mathcal{G}_{n}^{(i)}}\right)\\
%	&={\rm Var}\left(\sum_{q=0}^{n_{i}-1}\left|\widetilde{X}(t_{i-1+(q+1)K})-\widetilde{X}(t_{i-1+qK})\right|^{2}\right)\\
%	&=\sum_{q=0}^{n_{i}-1}{\rm Var}\left(\left|\widetilde{X}(t_{i-1+(q+1)K})-\widetilde{X}(t_{i-1+qK})\right|^{2}\right)\\
%	&\quad +2\sum_{0\leq{}q<r\leq n_{i}}{\rm Cov}\left(\left|\widetilde{X}(t_{i-1+(q+1)K})-\widetilde{X}(t_{i-1+qK})\right|^{2},\left|\widetilde{X}(t_{i-1+(r+1)K})-\widetilde{X}(t_{i-1+rK})\right|^{2}\right)\\
		&=\sum_{q=0}^{n_{i}-1}{\rm Var}\left(\left|\widetilde{X}(t_{i-1+(q+1)K})-\widetilde{X}(t_{i-1+qK})\right|^{2}\right)\\
	&\quad +2\sum_{q=0}^{n_{i}-2}{\rm Cov}\left(\left|\widetilde{X}(t_{i-1+(q+1)K})-\widetilde{X}(t_{i-1+qK})\right|^{2},\left|\widetilde{X}(t_{i-1+(q+2)K})-\widetilde{X}(t_{i-1+(q+1)K})\right|^{2}\right).
%	&=(n_{i}-1/2\pm 1/2){\rm Cov}\left(\left|\widetilde{X}(t_{i-1+K})-\widetilde{X}(t_{i-1})\right|^{2},\left|\widetilde{X}(t_{j-1+K})-\widetilde{X}(t_{j-1})\right|^{2}\right)\\
%	&\quad+(n_{j}-1/2\mp 1/2)
%	{\rm Cov}\left(\left|\widetilde{X}(t_{i-1+2K})-\widetilde{X}(t_{i-1+K})\right|^{2},\left|\widetilde{X}(t_{j-1+K})-\widetilde{X}(t_{j-1})\right|^{2}\right)\\
%	&=(n_{i}-1/2\pm 1/2)C\left(t_{i-1+K}-t_{j-1}\right)+(n_{j}-1/2\mp 1/2)
%	C\left(t_{j-1+K}-t_{i-1+K}\right)\\
%	&=\left(n_{j}-\frac{1\pm1}{2}\right)
%	C\left((j-i)\delta_{n}\right)+\left(n_{i}-\frac{1\mp1}{2}\right)C\left((K+i-j)\delta_{n}\right),
\end{align*}
Next, using the relationships
\begin{align*}
&{\rm Var}\left(\left|\widetilde{X}(t+\delta)-\widetilde{X}(t)\right|^{2}\right)
=2\sigma^{4}\delta^{2}+3\sigma^{4}\kappa \delta+8\sigma^{2}\bbe\left(\varepsilon^{2}\right)\delta+2\bbe\left(\varepsilon^{2}\right)^{2}+2\bbe\left(\varepsilon^{4}\right)\\
	&{\rm Cov}\left(\left|\widetilde{X}(t+\delta)-\widetilde{X}(t)\right|^{2},\left|\widetilde{X}(v)-\widetilde{X}(t+\delta)\right|^{2}\right)=\bbe\left(\varepsilon^{4}\right)-\bbe\left(\varepsilon^{2}\right)^{2},
\end{align*}
valid for any $t<t+\delta<v$, we get
\begin{align}
%\nonumber
	B_{i}
%	&=n_{i}\left(2\sigma^{4}\left(K\delta_{n}\right)^{2}+3\sigma^{4}\kappa \left(K\delta_{n}\right)+8\sigma^{2}\bbe\left(\varepsilon^{2}\right)\left(K\delta_{n}\right)+2\bbe\left(\varepsilon^{2}\right)^{2}+2\bbe\left(\varepsilon^{4}\right)\right)+2(n_{i}-1)\left(\bbe\left(\varepsilon^{4}\right)-\bbe\left(\varepsilon^{2}\right)^{2}\right)\\
		&=n_{i}\left(2\sigma^{4}\left(K\delta_{n}\right)^{2}+3\sigma^{4}\kappa \left(K\delta_{n}\right)+8\sigma^{2}\bbe\left(\varepsilon^{2}\right)\left(K\delta_{n}\right)\right)+2(2n_{i}-1)\bbe\left(\varepsilon^{4}\right)+2\bbe\left(\varepsilon^{2}\right)^{2}.\label{CovG}
\end{align}
Therefore, using that {$|1/n_{i}-K/n|\leq{}UK^{2}/n^{2}$ and $|1/n^{2}_{i}-K^{2}/n^{2}|\leq{}U K^{3}/n^{3}$, for a constant $U$ independent of $i$, $K$, $n$, and $T$, we have} 
$%\begin{align*}
	B
%	&=\frac{1}{K^{2}}\sum_{i=1}^{K}\frac{1}{K^{2}\delta_{n}^{2}n_{i}^{2}}\left\{n_{i}\left[2\sigma^{4}\left(K\delta_{n}\right)^{2}+3\sigma^{4}\kappa \left(K\delta_{n}\right)+8\sigma^{2}\bbe\left(\varepsilon^{2}\right)\left(K\delta_{n}\right)\right]+2(2n_{i}-1)\bbe\left(\varepsilon^{4}\right)+2\bbe\left(\varepsilon^{2}\right)^{2}\right\}\\
%	&=\frac{n^{2}}{K^{4}T^{2}}\left(\sum_{i=1}^{K}\frac{1}{n_{i}}\right)\left\{2\sigma^{4}\left(K\delta_{n}\right)^{2}+3\sigma^{4}\kappa \left(K\delta_{n}\right)+8\sigma^{2}\bbe\left(\varepsilon^{2}\right)\left(K\delta_{n}\right)+4\bbe\left(\varepsilon^{4}\right)\right\}\\
%	&\quad -\frac{2n^{2}}{K^{4}T^{2}}\left(\sum_{i=1}^{K}\frac{1}{n_{i}^{2}}\right)\left\{\bbe\left(\varepsilon^{4}\right)-\bbe\left(\varepsilon^{2}\right)^{2}\right\}\\
	=C_{1}-C_{2}+\mathcal{R}_{2},
$ %\end{align*}
where 
\begin{align*}
	C_{1}&=\frac{n}{K^{2}T^{2}}\left(2\sigma^{4}\left(K\delta_{n}\right)^{2}+3\sigma^{4}\kappa \left(K\delta_{n}\right)+8\sigma^{2}\bbe\left(\varepsilon^{2}\right)\left(K\delta_{n}\right)+4\bbe\left(\varepsilon^{4}\right)\right),\\
	C_{2}&=\frac{2}{KT^{2}}\left(\bbe\left(\varepsilon^{4}\right)-\bbe\left(\varepsilon^{2}\right)^{2}\right),
\end{align*}
and $\mathcal{R}_{2}=O_{u}\left((K/n)C_{1}\right)=O_{u}\left((K/n)C_{2}\right)$. 
Putting together $A$ and $B$ above,
\begin{align}\nonumber
	{\rm Var}\left(\hat\sigma^{2}_{n,K}\right)&=\frac{n}{K}\frac{K-1}{KT^{2}}\left(
%	24\bbe\left(\varepsilon^{2}\right)^{2}+8\bbe\left(\varepsilon^{4}\right)+
	{3\sigma^{4}\kappa K\delta_{n}+\frac{2}{3}\sigma^{4}K(2K-1)\delta_{n}^{2}}\right)\\
	&\quad+\frac{n}{K^{2}T^{2}}\left(2\sigma^{4}\left(K\delta_{n}\right)^{2}+3\sigma^{4}\kappa \left(K\delta_{n}\right)+8\sigma^{2}\bbe\left(\varepsilon^{2}\right)\left(K\delta_{n}\right)+4\bbe\left(\varepsilon^{4}\right)\right)\nonumber
	\\
	&\quad-\frac{2}{KT^{2}}\left(\bbe\left(\varepsilon^{4}\right)-\bbe\left(\varepsilon^{2}\right)^{2}\right)+\mathcal{R}_{1}+\mathcal{R}_{2}.
	\label{VarExpFrmb2}
\end{align}
Recalling that $\delta_{n}=T/n$ and using (\ref{MaxR}), we get the expression (\ref{VarSE}).
\hfill$\Box$

\smallskip\smallskip\noindent\textbf{Proof of Proposition \ref{VarSigmEst2}.} % [\textbf{Proof of Proposition \ref{VarSigmEst2}}.]
Let {$a_{K}:=\frac{K}{K-1}$} and {$b_{K}:=\frac{1}{T(K-1)}$}. Clearly, 
\begin{align*}
	{\rm Var}\left(\hat{\bar{\sigma}}^{2}_{n,K}\right)&=a_{K}^{2}{\rm Var}\left(\hat{\sigma}^{2}_{n,K}\right)+b_{K}^{2}{\rm Var}\left([\widetilde{X},\widetilde{X}]_{2}^{\bar{\mathcal{G}}_{n}}\right)-2 a_{K}b_{K}{\rm Cov}\left(\hat\sigma^{2}_{n,K},[\widetilde{X},\widetilde{X}]_{2}^{\bar{\mathcal{G}}_{n}}\right)
\end{align*}
From the expressions in Eqs.~(\ref{CovG})-(\ref{VarExpFrmb2}), we have
\begin{align*}
	{\rm Var}\left([\widetilde{X},\widetilde{X}]_{2}^{\bar{\mathcal{G}}_{n}}\right)&=n\left(2\sigma^{4}\delta_{n}^{2}+3\sigma^{4}\kappa \delta_{n}+8\sigma^{2}\bbe\left(\varepsilon^{2}\right)\delta_{n}\right)+2(2n-1)\bbe\left(\varepsilon^{4}\right)+2\bbe\left(\varepsilon^{2}\right)^{2}\\
	{\rm Var}\left(\hat\sigma^{2}_{n,K}\right)&=\frac{n}{K}\frac{K-1}{KT^{2}}\left(
%	24\bbe\left(\varepsilon^{2}\right)^{2}+8\bbe\left(\varepsilon^{4}\right)+
	{3\sigma^{4}\kappa K\delta_{n}+\frac{2}{3}\sigma^{4}K(2K-1)\delta_{n}^{2}}\right)\\
	&\quad+\frac{n-K+1}{K^{2}T^{2}}\left(2\sigma^{4}\left(K\delta_{n}\right)^{2}+3\sigma^{4}\kappa \left(K\delta_{n}\right)+8\sigma^{2}\bbe\left(\varepsilon^{2}\right)\left(K\delta_{n}\right)+4\bbe\left(\varepsilon^{4}\right)\right)
	\\
	&\quad-\frac{2}{KT^{2}}\left(\bbe\left(\varepsilon^{4}\right)-\bbe\left(\varepsilon^{2}\right)^{2}\right)+\mathcal{R}_{1}+\mathcal{R}_{2}.\nonumber
\end{align*}
To compute the last covariance, let us first note that
\begin{equation}
	{\rm Cov}\left(\hat\sigma^{2}_{n,K},[\widetilde{X},\widetilde{X}]_{2}^{\bar{\mathcal{G}}_{n}}\right)=\frac{1}{K}\sum_{i=1}^{K}{\frac{1}{T_{i}}}{\rm Cov}\left([\widetilde{X},\widetilde{X}]_{2}^{\mathcal{G}_{n}^{(i)}},[\widetilde{X},\widetilde{X}]_{2}^{\bar{\mathcal{G}}_{n}}\right)=:\frac{1}{K}\sum_{i=1}^{K}{\frac{1}{T_{i}}}B_{i}.
\end{equation}
Each covariance term on the right hand side above {can be computed as}
\begin{align*}
	B_{i}%&:={\rm Cov}\left([\widetilde{X},\widetilde{X}]_{2}^{\mathcal{G}_{n}^{(i)}},[\widetilde{X},\widetilde{X}]_{2}^{\bar{\mathcal{G}}_{n}}\right)\\
%	&={\rm Cov}\left(\sum_{q=0}^{n_{i}-1}\left|\widetilde{X}(t_{i-1+(q+1)K})-\widetilde{X}(t_{i-1+qK})\right|^{2},\sum_{r=0}^{n_{j}-1}\left|\widetilde{X}(t_{j-1+(r+1)K})-\widetilde{X}(t_{j-1+rK})\right|^{2}\right)\\
	&=\sum_{q=0}^{n_{i}-1}\sum_{r=0}^{n-1}{\rm Cov}\left(\left|\widetilde{X}(t_{i-1+(q+1)K})-\widetilde{X}(t_{i-1+qK})\right|^{2},\left|\widetilde{X}(t_{r+1})-\widetilde{X}(t_{r})\right|^{2}\right)\\
	&=(n_{i}-e_{i})\sum_{r=0}^{n-1}{\rm Cov}\left(\left|\widetilde{X}(t_{i-1+2K})-\widetilde{X}(t_{i-1+K})\right|^{2},\left|\widetilde{X}(t_{r+1})-\widetilde{X}(t_{r})\right|^{2}\right)\\
	&\quad+e_{i}\sum_{r=0}^{n-1}{\rm Cov}\left(\left|\widetilde{X}(t_{K})-\widetilde{X}(t_{0})\right|^{2},\left|\widetilde{X}(t_{r+1})-\widetilde{X}(t_{r})\right|^{2}\right),
\end{align*}
where above $e_{i}$ denote the number of subintervals in {the set} $\left\{[t_{i-1+qK},t_{i-1+(q+1)K}]\right\}_{{q=0}}^{{n_{i}-1}}$ which intersect the end points $0$ and $T$. {Obviously,} $\sum_{i=1}^{K}e_{i}=2$. Now, we use the following formulas:
\begin{align*}
	&{\rm Cov}\left(|\widetilde{X}(v)-\widetilde{X}(u)|^{2},|\widetilde{X}(v')-\widetilde{X}(u')|^{2}\right)=
	{2\sigma^{4}(v'-u')^{2}+3\kappa\sigma^{4}(v'-u')},\quad u<u'<v'<v\\
	&{\rm Cov}\left(|\widetilde{X}(t)-\widetilde{X}(s)|^{2},|\widetilde{X}(u)-\widetilde{X}(t)|^{2}\right)=
	\bbe\varepsilon^{4}-(\bbe\varepsilon^{2})^{2},\qquad s<t<u.
\end{align*}
We then get
%\begin{align*}
\(
	B_{i}%& =(n_{i}-e_{i})\left\{K\left({2\sigma^{4}\delta_{n}^{2}+3\kappa\sigma^{4}\delta_{n}}\right)+2(\bbe\varepsilon^{4}-(\bbe\varepsilon^{2})^{2})\right\}+e_{i}\left\{K\left({2\sigma^{4}\delta_{n}^{2}+3\kappa\sigma^{4}\delta_{n}}\right)+(\bbe\varepsilon^{4}-(\bbe\varepsilon^{2})^{2})\right\}\\
	=n_{i}\left\{K\left({2\sigma^{4}\delta_{n}^{2}+3\kappa\sigma^{4}\delta_{n}}\right)+2(\bbe\varepsilon^{4}-(\bbe\varepsilon^{2})^{2})\right\}-e_{i}(\bbe\varepsilon^{4}-(\bbe\varepsilon^{2})^{2}).
	\)
%\end{align*}
%where above we used the formula
%\[
%	{\rm Cov}\left(|\widetilde{X}(v)-\widetilde{X}(u)|^{2},|\widetilde{X}(v')-\widetilde{X}(u')|^{2}\right)=
%	{2\sigma^{4}(v'-u')^{2}+3\kappa\sigma^{4}(v'-u')},
%\]
%for any $u<u'<v'<v$ (see note in handwriting for this evaluation). 
Next, 
\begin{align*}
	{\rm Cov}\left(\hat\sigma^{2}_{n,K},[\widetilde{X},\widetilde{X}]_{2}^{\bar{\mathcal{G}}_{n}}\right)
%	&=\frac{1}{K}\sum_{i=1}^{K}\frac{1}{K\delta_{n}n_{i}}
%	n_{i}\left\{K\left({2\sigma^{4}\delta_{n}^{2}+3\kappa\sigma^{4}\delta_{n}}\right)+2(\bbe\varepsilon^{4}-(\bbe\varepsilon^{2})^{2})\right\}\\
%	&\quad-\frac{1}{K}\sum_{i=1}^{K}\frac{e_{i}}{K\delta_{n}n_{i}}(\bbe\varepsilon^{4}-(\bbe\varepsilon^{2})^{2})\\
%	&=\frac{n}{KT}
%	\left\{K\left({2\sigma^{4}\delta_{n}^{2}+3\kappa\sigma^{4}\delta_{n}}\right)+2(\bbe\varepsilon^{4}-(\bbe\varepsilon^{2})^{2})\right\}\\
%	&\quad-\frac{n}{K^{2}T}(\bbe\varepsilon^{4}-(\bbe\varepsilon^{2})^{2})\sum_{i=1}^{K}\frac{e_{i}}{n_{i}}\\
		&=\frac{n}{KT}
	\left\{K\left({2\sigma^{4}\delta_{n}^{2}+3\kappa\sigma^{4}\delta_{n}}\right)+2(\bbe\varepsilon^{4}-(\bbe\varepsilon^{2})^{2})\right\}+O_{u}\left(\frac{1}{KT}\right)
%	\\
%	&\quad-\frac{2}{KT}(\bbe\varepsilon^{4}-(\bbe\varepsilon^{2})^{2}).
\end{align*}	
%using that $\sum_{i=1}^{K}n_{i}=n-K+1$ and $\sum_{i=1}^{K}e_{i}=2$,
%\begin{align*}
%	{\rm Cov}\left(\hat\sigma^{2}_{n,K},[\widetilde{X},\widetilde{X}]_{2}^{\bar{\mathcal{G}}_{n}}\right)
%%	&=\frac{1}{KT}\sum_{i=1}^{K}\left(n_{i}\left\{K\left({2\sigma^{4}\delta_{n}^{2}+3\kappa\sigma^{4}\delta_{n}}\right)+2(\bbe\varepsilon^{4}-(\bbe\varepsilon^{2})^{2})\right\}-e_{i}(\bbe\varepsilon^{4}-(\bbe\varepsilon^{2})^{2})\right)\\
%%	&={\frac{n-K+1}{T}\left(2\sigma^{4}\delta_{n}^{2}+3\kappa\sigma^{4}\delta_{n}\right)}+ 2\frac{n-K+1}{TK}(\bbe\varepsilon^{4}-(\bbe\varepsilon^{2})^{2})-2\frac{1}{TK}(\bbe\varepsilon^{4}-(\bbe\varepsilon^{2})^{2})\\
%	&={\frac{n-K+1}{T}\left(2\sigma^{4}\delta_{n}^{2}+3\kappa\sigma^{4}\delta_{n}\right)}+ 2\frac{n-K}{TK}(\bbe\varepsilon^{4}-(\bbe\varepsilon^{2})^{2}).
%\end{align*}
Putting together the previous relationships, 
\begin{align*}
	{\rm Var}\left(\hat{\bar{\sigma}}^{2}_{n,K}\right)
	&=
	a_{K}^{2}{\rm Var}\left(\hat{\sigma}^{2}_{n,K}\right)+b_{K}^{2}{\rm Var}\left([\widetilde{X},\widetilde{X}]_{2}^{\bar{\mathcal{G}}_{n}}\right)-2 a_{K}b_{K}{\rm Cov}\left(\hat\sigma^{2}_{n,K},[\widetilde{X},\widetilde{X}]_{2}^{\bar{\mathcal{G}}_{n}}\right)\\
	&=\left(\frac{4\sigma^{4}K}{3n}+4\frac{n}{K^{2}T^{2}}\bbe(\varepsilon^{4})\right)+O_{u}\left(\frac{1}{n}\right)+O_{u}\left(\frac{n}{K^{3}T^{2}}\right)\\
	&\quad + \frac{1}{T^{2}K^{2}}\left(4n\bbe\varepsilon^{4}\right)-\frac{4n}{T^{2}K^{2}}(\bbe\varepsilon^{4}-(\bbe\varepsilon^{2})^{2})+O_{u}\left(\frac{1}{TK}\right)\\
	&={\frac{4\sigma^{4}K}{3n}+\frac{4n\left(\bbe\left(\varepsilon^{4}\right)+(\bbe\varepsilon^{2})^{2}\right)}{T^{2}K^{2}}+O_{u}\left(\frac{1}{n}\right)+O_{u}\left(\frac{n}{K^{3}T^{2}}\right)+O_{u}\left(\frac{1}{TK}\right).}
\end{align*}
\hfill$\Box$

\smallskip\noindent\textbf{Proof of Theorem \ref{VarKappaEst}.} %[\textbf{Proof of Theorem \ref{VarKappaEst}}.]
%The proof is similar to that of {Theorem \ref{VarSigmEst}}. 
Let us first write the variance of the estimator as follows:
\begin{align}	\nonumber%\label{VarExpFrmb}
	{\rm Var}\left(\hat\kappa_{n,K}\right)&=\frac{2}{9\sigma^{8}K^{2}}\sum_{1\leq i<j\leq K}\frac{1}{T_{i}T_{j}}{\rm Cov}\left([\widetilde{X},\widetilde{X}]_{4}^{\mathcal{G}_{n}^{(i)}},[\widetilde{X},\widetilde{X}]_{4}^{\mathcal{G}_{n}^{(j)}}\right)+\frac{1}{9\sigma^{8}K^{2}}\sum_{i=1}^{K}\frac{1}{T_{i}^{2}}{\rm Var}\left([\widetilde{X},\widetilde{X}]_{4}^{\mathcal{G}_{n}^{(i)}}\right)
	\\
	&=:A+B.\label{VarExpFrmb}
\end{align}
{Let us first note that we can replace $1/(T_{i}T_{j})=1/(K^{2}\delta_{n}^{2}n_{i}n_{j})$ with $1/T^{2}$ for any $1\leq{}i\leq{}j\leq{}K$, since $|1/(n_{i}n_{j})-K^{2}/n^{2}|\leq{}U K^{3}/n^{3}$, for a constant $U$ independent of $i,j,K,n,T$, and, thus, 
\begin{equation}\label{SmpRel}
	\left|\frac{1}{T_{i}T_{j}}-\frac{1}{T^{2}}\right|\leq{}U\frac{K}{T^{2}n}.
\end{equation}
Next, each covariance in the first term of (\ref{VarExpFrmb}) can be computed as:}
\begin{align*}
	A_{i,j}&:={\rm Cov}\left([\widetilde{X},\widetilde{X}]_{4}^{\mathcal{G}_{n}^{(i)}},[\widetilde{X},\widetilde{X}]_{4}^{\mathcal{G}_{n}^{(j)}}\right)\\
%	&={\rm Cov}\left(\sum_{q=0}^{n_{i}-1}\left|\widetilde{X}(t_{i-1+(q+1)K})-\widetilde{X}(t_{i-1+qK})\right|^{2},\sum_{r=0}^{n_{j}-1}\left|\widetilde{X}(t_{j-1+(r+1)K})-\widetilde{X}(t_{j-1+rK})\right|^{2}\right)\\
	&=\sum_{q=0}^{n_{i}-1}\sum_{r=0}^{n_{j}-1}{\rm Cov}\left(\left|\widetilde{X}(t_{i-1+(q+1)K})-\widetilde{X}(t_{i-1+qK})\right|^{4},\left|\widetilde{X}(t_{j-1+(r+1)K})-\widetilde{X}(t_{j-1+rK})\right|^{4}\right)\\
	&=n_{i}{\rm Cov}\left(\left|\widetilde{X}(t_{i-1+K})-\widetilde{X}(t_{i-1})\right|^{4},\left|\widetilde{X}(t_{j-1+K})-\widetilde{X}(t_{j-1})\right|^{4}\right)\\
	&\quad+(n_{j}-1)
	{\rm Cov}\left(\left|\widetilde{X}(t_{i-1+2K})-\widetilde{X}(t_{i-1+K})\right|^{4},\left|\widetilde{X}(t_{j-1+K})-\widetilde{X}(t_{j-1})\right|^{4}\right)\\
%	&=(n_{i}-1/2\pm 1/2)C\left(t_{i-1+K}-t_{j-1}\right)+(n_{j}-1/2\mp 1/2)
%	C\left(t_{j-1+K}-t_{i-1+K}\right)\\
	&=n_{i}C\left((j-i)\delta_{n},(K+i-j)\delta_{n},(j-i)\delta_{n}\right)\\
		&\quad+\left(n_{j}-1\right)
	C\left((K+i-j)\delta_{n},(j-i)\delta_{n},(K+i-j)\delta_{n}\right),
\end{align*}
where, for any $t,{s}_{1},{s}_{2},{s}_{3}>0$, 
\begin{equation}\label{CovArgm}
	C({s}_{1},{s}_{2},{s}_{3}):={\rm Cov}\left(\left|\widetilde{X}_{t+{s}_{1}+{s}_{2}}-\widetilde{X}_{t}\right|^{4},\left|\widetilde{X}_{t+{s}_{1}+{s}_{2}+{s}_{3}}-\widetilde{X}_{t+{s}_{1}}\right|^{4}\right),
\end{equation}
which again can be proved to be independent of $t$. Concretely, with the notation
$S:=X_{t+s_{1}+s_{2}}-X_{t+s_{1}}$, $U:=X_{t+s_{1}}-X_{t}+\varepsilon_{t+s_{1}+s_{2}}-\varepsilon_{t}$, and $V:=X_{t+s_{1}+s_{2}+s_{3}}-X_{t+s_{1}+s_{2}}+\varepsilon_{t+s_{1}+s_{2}+s_{3}}-\varepsilon_{t+s_{1}}$
\begin{align*}
	C(s_{1},s_{2},s_{3})
		&={\rm Cov}\left(\left|S+U\right|^{4},\left|S+V\right|^{4}\right)\\
	&={\rm Var}\left(S^{4}\right)+6\left[\bbe(U^{2})+\bbe(V^{2})\right]{\rm Cov}\left(S^{4},S^{2}\right)\\
	&\quad+36\bbe(U^{2})\bbe(V^{2}){\rm Var}\left(S^{2}\right)+16\bbe(U^{3})\bbe(V^{3}){\rm Var}\left(S\right)
\end{align*}
where above we used the independence of $S$, $U$, and $V$ as well as the fact that $\bbe U=\bbe V=\bbe S^{k}=0$ for any odd positive integer $k$. Upon computation of the relevant moments of $U$ and $V$, we get
\begin{align}\nonumber
	C(s_{1},s_{2},s_{3})
		&={\rm Var}\left(X_{s_{2}}^{4}\right)+6\left[\sigma^{2}(s_{1}+s_{3})+4\bbe\varepsilon^{2}\right]{\rm Cov}\left(X_{s_{2}}^{4},X_{s_{2}}^{2}\right)\\
	&\quad+6^{2}\left(\sigma^{2}s_{1}+2\bbe\varepsilon\right)\left(\sigma^{2}s_{3}+2\bbe\varepsilon\right){\rm Var}\left(X_{s_{2}}^{2}\right)+4^{2}\left(2\bbe\varepsilon^{3}\right)^{2}{\rm Var}\left(X_{s_{2}}\right).\label{VarForthMom}
\end{align}
Note that 
\[
	\bbe X_{s}^{k}=\bbe\left({(\sigma W_{\tau_{s}})}^{k}\right)=%\bbe\left(\bbe\left(\left.W_{\tau_{s}}^{k}\right|\tau_{s}\right)\right)=\bbe\left(\tau_{s}^{k/2}\bbe\left(\left.W_{1}^{k}\right|\tau_{s}\right)\right)
\sigma^{k}\bbe\left(W_{1}^{k}\right)\bbe\left(\tau_{s}^{k/2}\right)=\sigma^{k}\bbe\left(W_{1}^{k}\right)\left(s^{k/2}+
	\sum_{i=1}^{k/2-1}a_{k,i}s^{i}\right),
\]
for some constant $a_{k,i}$'s. 
We now proceed to analyze each term separately:
\begin{itemize}
	\item The contribution to $A$ due to ${\rm Var}\left(X_{s_{2}}^{4}\right)$ can be written as:
	\begin{align*}
		A^{(1)}&:=\frac{n}{K}\frac{2}{9\sigma^{8}K^{2}}\sum_{1\leq i<j\leq K}{\frac{1}{T_{i}T_{j}}}
		{\rm Var}\left({X}_{(K+i-j)\delta_{n}}^{4}\right)+\frac{n}{K}\frac{2}{9\sigma^{8}K^{2}}\sum_{1\leq i<j\leq K}{\frac{1}{T_{i}T_{j}}}
		{\rm Var}\left({X}_{(j-i)\delta_{n}}^{4}\right).
\end{align*}
Using {(\ref{SmpRel}) and that} ${\rm Var}\left(X_{t}^{4}\right)$ is a polynomial of degree $4$ in $t$ with the highest-degree term being $96\sigma^{8}t^{4}$,
\begin{align*}
			A^{(1)}&=\frac{n}{K}\frac{192\delta_{n}^{4}}{9K^{2}T^{2}}\left(\sum_{1\leq i<j\leq K}(K+i-j)^{4}+\sum_{1\leq i<j\leq K}(j-i)^{4}+O\left(K^{5}\right)\right)\\
%			&=\frac{n}{K}\frac{192\delta_{n}^{4}}{9K^{2}T^{2}}\left(\sum_{i=1}^{K-1}\sum_{\ell=1}^{K-i}(K-\ell)^{4}+\sum_{i=1}^{K-1}\sum_{\ell=1}^{K-i}\ell^{4}+O\left(K^{5}\right)\right)\\
%			&=\frac{n}{K}\frac{192\delta_{n}^{4}}{9K^{2}T^{2}}\left(\sum_{\ell=1}^{K-1}(K-\ell)^{5}+\sum_{\ell=1}^{K-1}(K-\ell)\ell^{4}+O\left(K^{5}\right)\right)\\
%			&=\frac{n}{K}\frac{192\delta_{n}^{4}}{9K^{2}T^{2}}K\sum_{\ell=1}^{K-1}\ell^{4}+O\left(\frac{T^{2}K^{2}}{n^{3}}\right)\\
			&=\frac{192}{5(9)}\frac{T^{2}K^{3}}{n^{3}}+O\left(\frac{T^{2}K^{2}}{n^{3}}\right).
\end{align*}

	\item Let us analyze the contribution to $A$ due to ${\rm Var}\left(X_{s_{2}}^{2}\right)$. {Using again (\ref{SmpRel}) and the variance formula in (\ref{MEVGM})}, the leading term is given by:
	\begin{align*}
		A^{(2)}&
		:=6^{2}\frac{n}{K}\frac{2}{9\sigma^{8}K^{2}}\sum_{1\leq i<j\leq K}\frac{1}{T_{i}T_{j}}
		\left(\sigma^{2}(j-i)\delta_{n}\right)^{2}{\rm Var}\left({X}_{(K+i-j)\delta_{n}}^{2}\right)\\
		&\quad+6^{2}\frac{n}{K}\frac{2}{9\sigma^{8}K^{2}T^{2}}\sum_{1\leq i<j\leq K}
		\left(\sigma^{2}(K+i-j)\delta_{n}\right)^{2}{\rm Var}\left({X}_{(j-i)\delta_{n}}^{2}\right)\\
		&=6^{2}\frac{n}{K}\frac{2}{9\sigma^{8}K^{2}T^{2}}\sum_{1\leq i<j\leq K}
		\left(\sigma^{2}(j-i)\delta_{n}\right)^{2}\left(3\sigma^{4}\kappa(K+i-j)\delta_{n}+2\sigma^{4}(K+i-j)^{2}\delta_{n}^{2}\right)\\
		&\quad+6^{2}\frac{n}{K}\frac{2}{9\sigma^{8}K^{2}T^{2}}\sum_{1\leq i<j\leq K}\left(\sigma^{2}(K+i-j)\delta_{n}\right)^{2}\left(3\sigma^{4}\kappa(j-i)\delta_{n}+2\sigma^{4}(j-i)^{2}\delta_{n}^{2}\right)\\
%				&=6^{2}\frac{n}{K}\frac{4\sigma^{8}\delta_{n}^{4}}{\sigma^{8}K^{2}T^{2}}\left(\sum_{1\leq i<j\leq K}(j-i)^{2}(K+i-j)^{2}+\sum_{1\leq i<j\leq K}(K+i-j)^{2}(j-i)^{2}\right)\\
%	&=\frac{(8)6^{2}}{9}\frac{T^{2}}{K^{3}n^{3}}\sum_{1\leq i<j\leq K}(j-i)^{2}(K+i-j)^{2}+O\left(\frac{T^{2}K^{2}}{n^{3}}\right)\\%=(8)6^{2}\frac{T^{2}}{K^{3}n^{3}}\frac{K^{6}}{(6)(5)(2)}\\
	&=\frac{(6)(4)(13)}{5(9)}\frac{T^{2}K^{3}}{n^{3}}+O\left(\frac{T^{2}K^{2}}{n^{3}}\right)
\end{align*}
	\item The contribution to $A$ due to ${\rm Cov}\left(X_{s_{2}}^{4},X_{s_{2}}^{2}\right)$ has the following leading term:
	\begin{align*}
		A^{(3)}&:=6\frac{n}{K}\frac{2}{9\sigma^{8}K^{2}T^{2}}\sum_{1\leq i<j\leq K}
		\left(2\sigma^{2}(j-i)\delta_{n}\right){\rm Cov}\left({X}_{(K+i-j)\delta_{n}}^{2},{X}_{(K+i-j)\delta_{n}}^{4}\right)\\
		&\quad+6\frac{n}{K}\frac{2}{9\sigma^{8}K^{2}T^{2}}\sum_{1\leq i<j\leq K}
		\left(2\sigma^{2}(K+i-j)\delta_{n}\right){\rm Cov}\left({X}_{(j-i)\delta_{n}}^{2},{X}_{(j-i)\delta_{n}}^{4}\right)\\
%		&=\frac{12^{2}}{9}\sigma^{8}\frac{n}{K}\frac{2\delta_{n}^{4}}{\sigma^{8}K^{2}T^{2}}\left(\sum_{1\leq i<j\leq K}
%		(j-i)(K+i-j)^{3}+\sum_{1\leq i<j\leq K}(K+i-j)(j-i)^{3}+O\left(K^{5}\right)\right)\\
			&=\frac{12^{2}(2)}{(5)(4)(9)}\frac{T^{2}K^{3}}{n^{3}}+O\left(\frac{T^{2}K^{2}}{n^{3}}\right)
\end{align*}
where above we used that ${\rm Cov}\left(X_{s}^{2},X_{s}^{4}\right)=\bbe X_{s}^{6}-\bbe(X_{s}^{2})\bbe(X_{s}^{4})%=\sigma^{6}s^{3}\left(15-3\right)+\text{h.o.t.}
=12\sigma^6 s^{3}+{\text{h.o.t.}}$, {where h.o.t. mean higher order terms.}
\item Finally, the contribution to $A$ due to ${\rm Var}\left(X_{s_{2}}\right)$ will generate a term of smaller order than $T^{2}K^{3}/n^{3}$. Indeed, 
	\begin{align*}
		A^{(4)}&:=4^{2}\left(2\bbe\varepsilon^{3}\right)^{2}\frac{n}{K}\frac{2}{9\sigma^{8}K^{2}T^{2}}\sum_{1\leq i<j\leq K}\left(
		{\rm Var}\left({X}_{(K+i-j)\delta_{n}}\right)+{\rm Var}\left({X}_{(j-i)\delta_{n}}\right)\right)\\
%		&:=4^{2}\left(2\bbe\varepsilon^{3}\right)^{2}\frac{n}{K}\frac{2}{9\sigma^{8}K^{2}T^{2}}\sigma^{2}\delta_{n}
%		\sum_{1\leq i<j\leq K}\left((K+i-j)+(j-i)\right)\\\
%				&:=4^{2}\left(2\bbe\varepsilon^{3}\right)^{2}\frac{n}{K}\frac{2}{\sigma^{8}K^{2}T^{2}}\sigma^{2}\delta_{n}
%		%\sum_{\ell=1}^{K-1}\ell
%		\frac{K^{2}(K-1)}{2}\\
		&=\frac{4^{2}}{9}\left(2\bbe\varepsilon^{3}\right)^{2}\frac{1}{\sigma^{6}T}.
\end{align*}

\end{itemize}
Putting together the above relationships, 
\begin{align*}
	A&=\frac{192}{5(9)}\frac{T^{2}K^{3}}{n^{3}}+\frac{(6)(4)(13)}{5(9)}\frac{T^{2}K^{3}}{n^{3}}+
\frac{6^{2}(2)}{5(9)}\frac{T^{2}K^{3}}{n^{3}}
++O\left(\frac{T^{2}K^{2}}{n^{3}}\right)
\\
&=\frac{576}{5(9)}\frac{T^{2}K^{3}}{n^{3}}+
O\left(\frac{T^{2}K^{2}}{n^{3}}\right).
\end{align*}
Now, we consider the second term in (\ref{VarExpFrmb}), which we denote $B$. Each variance {term, $B_{i}:={\rm Var}\big([\widetilde{X},\widetilde{X}]_{4}^{\mathcal{G}_{n}^{(i)}}\big)$,} of $B$ can be written as 
\begin{align*}
	B_{i}
		&=\sum_{q=0}^{n_{i}-1}{\rm Var}\left(\left|\widetilde{X}(t_{i-1+(q+1)K})-\widetilde{X}(t_{i-1+qK})\right|^{4}\right)\\
	&\quad +2\sum_{q=0}^{n_{i}-2}{\rm Cov}\left(\left|\widetilde{X}(t_{i-1+(q+1)K})-\widetilde{X}(t_{i-1+qK})\right|^{4},\left|\widetilde{X}(t_{i-1+(q+2)K})-\widetilde{X}(t_{i-1+(q+1)K})\right|^{4}\right).
\end{align*}
Next, using arguments similar to those following (\ref{CovArgm}), 
\begin{align}\label{VarCovTrm}
&{\rm Var}\left(\left|\widetilde{X}_{t+s}-\widetilde{X}_{t}\right|^{4}\right)
={\rm Var}\left(\left|{X}_{t+s}-{X}_{t}\right|^{4}\right)+\text{h.o.t.}=96\sigma^{8}s^{4}+\text{h.o.t.},\\%\sum_{i=1}^{3}a_{i}s^{i}\\
&{\rm Cov}\left(\left|\widetilde{X}(t+s_{1})-\widetilde{X}(t)\right|^{4},\left|\widetilde{X}(t+s_{1}+s_{2})-\widetilde{X}(t+s_{1})\right|^{4}\right)
%=-36\bbe X_{s_{1}}^{2}\bbe X_{s_{2}}^{2}\bbe\varepsilon^{2}+{\rm h.o.t.}
=-36\sigma^{4}\bbe\varepsilon^{2}s_{1}s_{2}+{\rm h.o.t.}\nonumber
\end{align}
valid for any $t,s_{1},s_{2}>0$ and where, again,  h.o.t.~means higher order terms. Therefore, 
$B_{i}=n_{i}\left(96\sigma^{8}\left(K\delta_{n}\right)^{4}\right)+{\rm h.o.t.}$
%\end{align*}
 and, thus,  %using that $\sum_{i=1}^{K}n_{i}=n-K+1$,
\begin{align*}
	B
%	&=\frac{96(K\delta_{n})^{4}}{K^{2}T^{2}}\sum_{i=1}^{K}n_{i}\\
%	&
%	=\frac{n-K+1}{K^{2}T^{2}\sigma^{8}}\left(114\sigma^{8}\left(K\delta_{n}\right)^{4}+O\left((K\delta_{n})^{3}\right)\right)
	=\frac{96}{9}\frac{K^{2}T^{2}}{n^{3}}+{\rm h.o.t.},
\end{align*}
which shows that $B=O(T^{2}K^{2}/n^{3})$. 
Finally, 
\begin{align*}\label{VarExpFrmb}
	{\rm Var}\left(\hat\kappa_{n,K}\right)&=\frac{576}{5(9)}\frac{T^{2}K^{3}}{n^{3}}+\frac{96}{9}\frac{K^{2}T^{2}}{n^{3}}+O\left(\frac{KT}{n^{2}}\right),%+
%114\frac{K^{2}T^{2}}{n^{3}}+O\left(\frac{KT}{n^{2}}\right),
%\\
%	&\quad+\frac{n-K+1}{K^{2}T^{2}}\left(2\sigma^{4}\left(K\delta_{n}\right)^{2}+3\sigma^{4}\kappa \left(K\delta_{n}\right)+8\sigma^{2}\bbe\left(\varepsilon^{2}\right)\left(K\delta_{n}\right)+4\bbe\left(\varepsilon^{4}\right)\right)\\
%	&\quad-\frac{2}{KT^{2}}\left(\bbe\left(\varepsilon^{4}\right)-\bbe\left(\varepsilon^{2}\right)^{2}\right).
\end{align*}
which implies the result. \hfill$\Box$

\noindent
\textbf{Proof of Theorem \ref{ThrmVarKappaEst2}.} %[{\bf Proof of Theorem \ref{ThrmVarKappaEst2}}.]
Let {$a_{K}:=\frac{K}{K-1}$, $b_{K}:=\frac{1}{3\sigma^{4}T(K-1)}$, and $c_{K}:=\frac{2}{n\sigma^{2}}$ so that}
\begin{align*}
	{\rm Var}\left(\hat{\bar{\kappa}}_{n}\right)&=a_{K}^{2}{\rm Var}\left(\hat{\kappa}_{n,K}\right)+b_{K}^{2}{\rm Var}\left([\widetilde{X},\widetilde{X}]_{4}^{\bar{\mathcal{G}}_{n}}\right)\\
	&\quad+c_{K}^{2}{\rm Var}\left([\widetilde{X},\widetilde{X}]_{2}^{\bar{\mathcal{G}}_{n}}\right)-2 a_{K}b_{K}{\rm Cov}\left(\hat\kappa_{n,K},[\widetilde{X},\widetilde{X}]_{4}^{\bar{\mathcal{G}}_{n}}\right)\\
	&\quad-2 a_{K}c_{K}{\rm Cov}\left(\hat\kappa_{n,K},[\widetilde{X},\widetilde{X}]_{2}^{\bar{\mathcal{G}}_{n}}\right)+2 b_{K}c_{K}{\rm Cov}\left([\widetilde{X},\widetilde{X}]_{4}^{\bar{\mathcal{G}}_{n}},[\widetilde{X},\widetilde{X}]_{2}^{\bar{\mathcal{G}}_{n}}\right)
\end{align*}
As in the case of the variance of $\hat{\bar{\sigma}}_{n,K}$, we are looking for the terms {having} the highest power of $K$ and the terms with the highest power of $n$ (and {the least} negative power of $K$). For ${\rm Var}\left(\hat{\kappa}_{n,K}\right)$, the highest power of $K$ is given {in Eq.~(\ref{VarKappEq})}. To find the highest power of $n$, we recall from the proof of Theorem \ref{VarKappaEst} that the variance can be decomposed into two terms, called $A$ and $B$ therein. The term with the highest power $n$ in $A$ is due to the term $4^{2}(2\bbe\varepsilon^{3})^{2}{\rm Var}(X_{s_{2}})$ in (\ref{VarForthMom}) and is of order $n^{0}$. In order to determine the term with the highest power of $n$ in $B$, note that this will be due to the constant terms of the variance and covariance in Eqs.~(\ref{VarCovTrm}). These are given by 
\begin{align}
&{\rm Var}\left(\left|\widetilde{X}_{t+s}-\widetilde{X}_{t}\right|^{4}\right)
%={\rm Var}\left((\varepsilon_{t+s}-\varepsilon_{t})^{4}\right)+\text{h.o.t.}
={\rm Var}\left((\varepsilon_{2}-\varepsilon_{1})^{4}\right)+\text{h.o.t.},\\%\sum_{i=1}^{3}a_{i}s^{i}\\
&{\rm Cov}\left(\left|\widetilde{X}(t+s_{1})-\widetilde{X}(t)\right|^{4},\left|\widetilde{X}(t+s_{1}+s_{2})-\widetilde{X}(t+s_{1})\right|^{4}\right)=
%{\rm Cov}\left(\left|\varepsilon_{t+s_{1}}-\varepsilon_{t}\right|^{4},\left|\varepsilon_{t+s_{1}+s_{2}}-\varepsilon_{t+s_{1}}\right|^{4}\right)+{\rm h.o.t.}\nonumber
{\rm Cov}\left(\left|\varepsilon_{2}-\varepsilon_{1}\right|^{4},\left|\varepsilon_{3}-\varepsilon_{2}\right|^{4}\right)+{\rm h.o.t.}\nonumber
\end{align}
where h.o.t. means higher order term (as powers of $s$, $s_{1}$, and $s_{2}$). These terms contribute to $B$ as follows:
\begin{align*}
	B&:=\frac{1}{9\sigma^{8}K^{2}T^{2}}\sum_{i=1}^{K}{\rm Var}\left([\widetilde{X},\widetilde{X}]_{4}^{\mathcal{G}_{n}^{(i)}}\right)
%	&=\frac{1}{9\sigma^{8}K^{2}T^{2}}\sum_{i=1}^{K}\left(n_{i}{\rm Var}\left((\varepsilon_{2}-\varepsilon_{1})^{4}\right)+2(n_{i}-1){\rm Cov}\left(\left|\varepsilon_{2}-\varepsilon_{1}\right|^{4},\left|\varepsilon_{3}-\varepsilon_{2}\right|^{4}\right)\right)+{\rm h.o.t.}\\
	=\frac{n}{9\sigma^{8}K^{2}T^{2}}d(\varepsilon)+{\rm h.o.t.},
\end{align*}
where $d(\varepsilon):={\rm Var}\left((\varepsilon_{2}-\varepsilon_{1})^{4}\right)+2{\rm Cov}\left(\left|\varepsilon_{2}-\varepsilon_{1}\right|^{4},\left|\varepsilon_{3}-\varepsilon_{2}\right|^{4}\right)$.
Now we consider $b_{K}^{2}{\rm Var}\left([\widetilde{X},\widetilde{X}]_{4}^{\bar{\mathcal{G}}_{n}}\right)$. As done with $B$, the term with the highest degree in $n$ is $\frac{n}{9\sigma^{8}T^{2}K^{2}}d(\varepsilon)$.
Clearly, all the terms in $c_{K}^{2}{\rm Var}\left([\widetilde{X},\widetilde{X}]_{2}^{\bar{\mathcal{G}}_{n}}\right)$ are of higher order than $n/(T^{2}K^{2})$.  
To compute ${\rm Cov}\left(\hat\kappa^{2}_{n,K},[\widetilde{X},\widetilde{X}]_{4}^{\bar{\mathcal{G}}_{n}}\right)$, let us first note that
\begin{equation}
	{\rm Cov}\left(\hat\kappa_{n,K},[\widetilde{X},\widetilde{X}]_{4}^{\bar{\mathcal{G}}_{n}}\right)=\frac{1}{3\sigma^{4}K}\sum_{i=1}^{K}{\frac{1}{T_{i}}}{\rm Cov}\left([\widetilde{X},\widetilde{X}]_{4}^{\mathcal{G}_{n}^{(i)}},[\widetilde{X},\widetilde{X}]_{4}^{\bar{\mathcal{G}}_{n}}\right).
\end{equation}
Each covariance term on the right hand side above, which is denoted $B_{i}$, is given by 
\begin{align*}
	B_{i}%&:={\rm Cov}\left([\widetilde{X},\widetilde{X}]_{2}^{\mathcal{G}_{n}^{(i)}},[\widetilde{X},\widetilde{X}]_{2}^{\bar{\mathcal{G}}_{n}}\right)\\
%	&={\rm Cov}\left(\sum_{q=0}^{n_{i}-1}\left|\widetilde{X}(t_{i-1+(q+1)K})-\widetilde{X}(t_{i-1+qK})\right|^{2},\sum_{r=0}^{n_{j}-1}\left|\widetilde{X}(t_{j-1+(r+1)K})-\widetilde{X}(t_{j-1+rK})\right|^{2}\right)\\
	&=\sum_{q=0}^{n_{i}-1}\sum_{r=0}^{n-1}{\rm Cov}\left(\left|\widetilde{X}(t_{i-1+(q+1)K})-\widetilde{X}(t_{i-1+qK})\right|^{4},\left|\widetilde{X}(t_{r+1})-\widetilde{X}(t_{r})\right|^{4}\right)\\
	&=(n_{i}-e_{i})\sum_{r=0}^{n-1}{\rm Cov}\left(\left|\widetilde{X}(t_{i-1+2K})-\widetilde{X}(t_{i-1+K})\right|^{4},\left|\widetilde{X}(t_{r+1})-\widetilde{X}(t_{r})\right|^{4}\right)\\
	&\quad+e_{i}\sum_{r=0}^{n-1}{\rm Cov}\left(\left|\widetilde{X}(t_{K})-\widetilde{X}(t_{0})\right|^{4},\left|\widetilde{X}(t_{r+1})-\widetilde{X}(t_{r})\right|^{4}\right),
\end{align*}
where above $e_{i}$ denote the number of subintervals in $\left\{[t_{i-1+qK},t_{i-1+(q+1)K}]\right\}_{q=0}^{n_{i}-1}$ which intersect the end points $0$ and $T$. Now, it turns out that
\begin{align}\label{Cov2and4}
	{\rm Cov}\left(|\widetilde{X}(v)-\widetilde{X}(u)|^{4},|\widetilde{X}(v')-\widetilde{X}(u')|^{4}\right)&\asymp n^{-1},\quad u<u'<v'<v\\
	{\rm Cov}\left(|\widetilde{X}(t)-\widetilde{X}(s)|^{4},|\widetilde{X}(u)-\widetilde{X}(t)|^{4}\right)&={\rm Cov}\left(|\varepsilon_{2}-\varepsilon_{1}|^{4},|\varepsilon_{3}-\varepsilon_{2}|^{4}\right)=:g(\varepsilon),\quad s<t<u,\nonumber
\end{align}
where here $a_{n}\asymp b_{n}$ means $\lim_{n\to{}\infty}a_{n}/b_{n}\in\bbr\backslash\{0\}$. We then conclude that $B_{i}=2n_{i}g(\varepsilon)-e_{i}g(\varepsilon)+{\rm h.o.t.}$.
%where above we used the formula
%\[
%	{\rm Cov}\left(|\widetilde{X}(v)-\widetilde{X}(u)|^{2},|\widetilde{X}(v')-\widetilde{X}(u')|^{2}\right)=
%	{2\sigma^{4}(v'-u')^{2}+3\kappa\sigma^{4}(v'-u')},
%\]
%for any $u<u'<v'<v$ (see note in handwriting for this evaluation). 
Then, it is clear that %using (\ref{SmpRel}) and that $|n_{i}-(n/K)|\leq{}2$,
\begin{align*}
	{\rm Cov}\left(\hat\kappa_{n,K},[\widetilde{X},\widetilde{X}]_{4}^{\bar{\mathcal{G}}_{n}}\right)
%	&=\frac{1}{3\sigma^{4}KT}\sum_{i=1}^{K}\left(2n_{i}g(\varepsilon)-e_{i}g(\varepsilon)\right)\\
%	&= \frac{2}{3\sigma^{4}}\frac{n-K+1}{KT}g(\varepsilon)-\frac{2}{3\sigma^{4}}\frac{1}{TK}g(\varepsilon)\\
	&=\frac{2}{3\sigma^{4}}\frac{{n}}{TK}g(\varepsilon)+{\rm {h.o.t.}}
\end{align*}
Therefore, the contribution here is $-\frac{4n}{9\sigma^{8}T^{2} K^{2}}g(\varepsilon)$.
%\[
%%	-\frac{2}{3\sigma^{4}TK}\left(\frac{2}{3\sigma^{4}}\frac{n-K}{TK}g(\varepsilon)\right)=
%	-\frac{4n}{9\sigma^{8}T^{2} K^{2}}g(\varepsilon).
%\]
Given that $c_{K}$ is of order $n^{-1}$, it is not hard to see that the term $-2 a_{K}c_{K}{\rm Cov}\left(\hat\kappa_{n,K},[\widetilde{X},\widetilde{X}]_{2}^{\bar{\mathcal{G}}_{n}}\right)$ is of an order smaller than $n$. 
Finally, consider the term corresponding to $D_{n}:={\rm Cov}\left([\widetilde{X},\widetilde{X}]_{4}^{\bar{\mathcal{G}}_{n}},[\widetilde{X},\widetilde{X}]_{2}^{\bar{\mathcal{G}}_{n}}\right)$. Note that 
\begin{align*}
	D_{n}%&:={\rm Cov}\left([\widetilde{X},\widetilde{X}]_{2}^{\mathcal{G}_{n}^{(i)}},[\widetilde{X},\widetilde{X}]_{2}^{\bar{\mathcal{G}}_{n}}\right)\\
%	&={\rm Cov}\left(\sum_{q=0}^{n_{i}-1}\left|\widetilde{X}(t_{i-1+(q+1)K})-\widetilde{X}(t_{i-1+qK})\right|^{2},\sum_{r=0}^{n_{j}-1}\left|\widetilde{X}(t_{j-1+(r+1)K})-\widetilde{X}(t_{j-1+rK})\right|^{2}\right)\\
	&=\sum_{q=0}^{n-1}\sum_{r=0}^{n-1}{\rm Cov}\left(\left|\widetilde{X}_{t_{q+1}}-\widetilde{X}_{t_{q}}\right|^{4},\left|\widetilde{X}_{t_{r+1}}-\widetilde{X}_{t_{r}}\right|^{2}\right)\\
%	&=(n-2)\left({\rm Cov}\left(\left|\widetilde{X}_{t_{1}}-\widetilde{X}_{t_{0}}\right|^{4},\left|\widetilde{X}_{t_{1}}-\widetilde{X}_{t_{0}}\right|^{2}\right)+2{\rm Cov}\left(\left|\widetilde{X}_{t_{1}}-\widetilde{X}_{t_{0}}\right|^{4},\left|\widetilde{X}_{t_{2}}-\widetilde{X}_{t_{1}}\right|^{2}\right)\right)\\
%	&\quad+2\left({\rm Cov}\left(\left|\widetilde{X}_{t_{1}}-\widetilde{X}_{t_{0}}\right|^{4},\left|\widetilde{X}_{t_{1}}-\widetilde{X}_{t_{0}}\right|^{2}\right)+{\rm Cov}\left(\left|\widetilde{X}_{t_{1}}-\widetilde{X}_{t_{0}}\right|^{4},\left|\widetilde{X}_{t_{2}}-\widetilde{X}_{t_{1}}\right|^{2}\right)\right)\\
		&=n\left({\rm Cov}\left(\left|\widetilde{X}_{t_{1}}-\widetilde{X}_{t_{0}}\right|^{4},\left|\widetilde{X}_{t_{1}}-\widetilde{X}_{t_{0}}\right|^{2}\right)+2{\rm Cov}\left(\left|\widetilde{X}_{t_{1}}-\widetilde{X}_{t_{0}}\right|^{4},\left|\widetilde{X}_{t_{2}}-\widetilde{X}_{t_{1}}\right|^{2}\right)\right)\\
	&\quad-{2}{\rm Cov}\left(\left|\widetilde{X}_{t_{1}}-\widetilde{X}_{t_{0}}\right|^{4},\left|\widetilde{X}_{t_{2}}-\widetilde{X}_{t_{1}}\right|^{2}\right)
\end{align*}
Using (\ref{Cov2and4}), it is clear that {$D_{n}\asymp n$}. Hence,
\[
2 b_{K}c_{K}{\rm Cov}\left([\widetilde{X},\widetilde{X}]_{4}^{\bar{\mathcal{G}}_{n}},[\widetilde{X},\widetilde{X}]_{2}^{\bar{\mathcal{G}}_{n}}\right)\asymp 
\frac{2}{3\sigma^{4}TK}.
\]
Finally, we obtain that 
\begin{align*}
	{\rm Var}\left(\hat{\bar{\kappa}}_{n,K}\right)&=\frac{64}{5}\frac{T^{2}K^{3}}{n^{3}}+\frac{n}{9\sigma^{8}K^{2}T^{2}}d(\varepsilon)+\frac{n}{9\sigma^{8}T^{2}K^{2}}d(\varepsilon)-\frac{4n}{9\sigma^{8}T^{2} K^{2}}g(\varepsilon)+{\rm h.o.t.}
%	\\
%	&=
%	a_{n,K}^{2}{\rm Var}\left(\hat{\sigma}^{2}_{n,K}\right)+b_{n,K}^{2}{\rm Var}\left([\widetilde{X},\widetilde{X}]_{2}^{\bar{\mathcal{G}}_{n}}\right)-2 a_{n,K}b_{n,K}{\rm Cov}\left(\hat\sigma^{2}_{n,K},[\widetilde{X},\widetilde{X}]_{2}^{\bar{\mathcal{G}}_{n}}\right)\\
%	&=\left(\frac{4\sigma^{4}K}{3n}+4\frac{n}{K^{2}T^{2}}\bbe(\varepsilon^{4})\right)+O\left(\frac{1}{n}\right)+O\left(\frac{n}{K^{3}T^{2}}\right)\\
%	&\quad + \frac{1}{T^{2}K^{2}}\left(4n\bbe\varepsilon^{4}\right)-\frac{4n}{T^{2}K^{2}}(\bbe\varepsilon^{4}-(\bbe\varepsilon^{2})^{2})+O\left(\frac{1}{TK}\right)\\
%	&=\frac{576}{5}\frac{T^{2}K^{3}}{n^{3}}+\frac{n}{9\sigma^{8}T^{2}K^{2}}e(\varepsilon)+{\rm h.o.t.},
\end{align*}
which implies the result.% since $d(\varepsilon)={\rm Var}\left((\varepsilon_{2}-\varepsilon_{1})^{4}\right)+2{\rm Cov}\left(\left|\varepsilon_{2}-\varepsilon_{1}\right|^{4},\left|\varepsilon_{3}-\varepsilon_{2}\right|^{4}\right)$ and $g(\varepsilon)={\rm Cov}\left(|\varepsilon_{2}-\varepsilon_{1}|^{4},|\varepsilon_{3}-\varepsilon_{2}|^{4}\right)$.
%where $e(\varepsilon)=2{\rm Var}\left((\varepsilon_{2}-\varepsilon_{1})^{4}\right)$.
\hfill$\Box$

\bibliographystyle{jtbnew}

\end{document}